\documentclass[12pt,a4paper]{book}

\usepackage[ansinew]{inputenc}

\usepackage[french]{babel}

\usepackage{amssymb,amsmath, amsfonts}

\usepackage [dvips] {graphicx}

\usepackage{amsthm}

\usepackage{stmaryrd}

\newtheorem{df}{Définition}
\newtheorem{prop} {Proposition}
\newtheorem{lm} [prop]{Lemme}
\newtheorem{thm} [prop] {Théorème}
\newtheorem{cor} [prop] {Corollaire}

\newtheorem{notations}{Notations}

\theoremstyle{remark}
\newtheorem{rmq}{Remarque} 
\newtheorem{exm}{Exemple} 
\newtheorem{exc}{Exercice}

\author{S. Dugowson}

\title{Introduction aux dynamiques \\ catégoriques connectives}

\begin{document}

\maketitle


\noindent \textbf{Abstract}. \textsc{Introduction To Categorical Connective Dynamics.} --- This text is a continuation to my former article ``On Connectivity Spaces". It takes into account that connectivity spaces gives rise to phenomena which are essentially dynamic. 
In a first stage, the representation of finite connectivity spaces by links (Brunn-Debrunner-Kanenobu's theorem) leads to the notion of connective representation. But examples of connective representations often come from dynamical systems. And this is even more obvious when we study the adjoint notion of connective foliation. To apply those notions to dynamics, we first need to consider dynamical systems in an unified way. This is done with a categorical point of view on temporalities and dynamics. It is then possible to define categorical connective dynamics, and to apply to them the various connective notions, specially the connectivity order of a connectivity space.

\mbox{}\linebreak

\noindent\textbf{Résumé}. Ce texte est la suite de mon article \og On connectivity Spaces\fg. Cette fois, c'est la nature essentiellement dynamique de ces espaces qui nous intéresse. Dans un premier temps, la représentation par entrelacs des espaces connectifs finis (théorème de Brunn-Debrunner-Kanenobu) conduit à la notion de représentation connective. Mais les exemples de telles représentations proviennent souvent de systèmes dynamiques. Cela est encore plus flagrant avec la notion adjointe de feuilletage connectif. Pour appliquer ces notions aux systèmes dynamiques, il nous faut d'abord considérer ces systèmes de façon unifiée. Cela est rendu possible par l'adoption d'un point de vue catégorique sur les temporalités et les dynamiques. Il est alors possible de définir les dynamiques catégoriques connectives, et de leur appliquer les diverses notions connectives, en particulier celle d'ordre connectif d'un espace connectif. 

\mbox{}\linebreak


\mbox{}\linebreak


\chapter*{Introduction}

Les espaces connectifs constituent une \og catégorie topologique\footnote{Au sens de \cite{joycats:2005}. Cela revient à peu près à dire que les structures connectives sur un ensemble donné de points s'organisent (fonctoriellement) en un treillis complet, exactement comme le font les structures de la topologie générale.} \fg, mais elles n'entrent pas pour autant dans les définitions de la topologie générale. Faisant suite à notre article \og On connectivity Spaces\fg\, \cite{Dugowson:201012}, dans lequel ces espaces et plusieurs notions fondamentales les concernant sont définis, le présent ouvrage commence par un chapitre de  courts rappels de ces notions, avec également quelques compléments qui nous serons utiles pour la suite (en particulier la notion d'ordre connectif pour des espaces connectifs quelconques).

La considération des structures connectives  les plus simples, telle celle du n\oe ud borroméen, conduisent naturellement à s'intéresser à la représentation par entrelacs de ces structures. Comme il est indiqué dans \cite{Dugowson:201012}, le théorème de Brunn-Debrunner-Kanenobu répond positivement à la question de la représen\-tabilité par entrelacs de toute structure connective finie. 

Si l'on souhaite un jour --- nous ne le ferons pas ici --- traiter cette même question mais dans le cas des espaces connectifs infinis, il nous faut commencer par préciser la notion générale de représentation d'un espace connectif dans un autre. Ceci fait, la notion de représentabilité \textit{par entrelacs} apparaît alors comme un cas particulier de la notion générale, à condition de munir l'espace où la représentation a lieu d'une structure connective particulière, définie \og par séparation \fg. Ces notions, et la définition de certaines catégories de représentations connectives, font l'objet de la première section du chapitre 2. 

Les premiers exemples de représentations par entrelacs de structures {con\-nectives} infinies sont issus de la mécanique (fibration de Hopf, tores de Liouville-Arnold pour les systèmes hamiltoniens intégrables, ...) et des systè\-mes dynamiques (système dynamique de Lorentz, de Ghrist, etc). Ce constat nous a conduit à la notion de feuilletage connectif, définie et développée dans la deuxième section du chapitre 2.

Les relations entre représentations et feuilletages connectifs s'avé\-rent fonctorielles. Plus précisément, dans certains cas, un couple de foncteurs adjoints structure ces relations. Ces questions sont traitées dans la troisième section du chapitre 2.

Afin de pouvoir appliquer les notions précédentes aux systèmes dynamiques, et en particulier la notion d'ordre connectif définie en toute généralité dans les compléments du chapitre 1, nous devons pouvoir nous appuyer sur une définition suffisamment générale des systèmes dynamiques connectifs, elle-même construite sur une conception unifiée des systèmes dynamiques. Il n'y aurait en effet guère eu de sens à faire appel à des définitions fragmentées des systèmes dynamiques, traitant en particulier de manière séparée les systèmes discrets et les systèmes continus, alors que les structures connectives que nous utilisons résultent précisément d'un souci d'unification. 

J'ai donc profité de l'occasion pour reprendre à la base la notion de système dynamique, donc aussi celles de temps, d'instants, etc., en m'appuyant sur la propriété fondamentale de ces systèmes : la transformation de la composition des écoulements temporels en composition des transitions correspondantes. 
 
En faisant des écoulements temporels les ingrédients fondamentaux de toute dynamique, et donc du temps, les instants se trouvent en quelque sorte relégués au second plan. En particulier, le fait que l'ensemble des instants puisse être ordonné, ne serait-ce que partiellement, ne constitue plus une exigence première. Sans rechercher la généralisation abstraite à tout prix et pour elle-même, le cadre qui nous a paru le mieux adapté à un point de vue unificateur est un cadre catégorique, de sorte que les écoulement temporels seront définis comme les flèches d'une petite catégorie quelconque. Sur cette base, les dynamiques --- non nécessairement déterministes, puisque là aussi il nous a paru nécessaire de ne pas nous enfermer \textit{a priori} dans un cadre dont plusieurs développements scientifiques déjà anciens ont clairement montré la nature par trop restrictive --- sont alors définies par des foncteurs. Les instants apparaissent alors comme les états de dynamiques (déterministes) particulières, notamment celles que nous baptisons \textit{existentielles} et \textit{essentielles}. Les dynamiques ainsi définies s'organisent elles-mêmes en catégories, et des foncteurs intéressants apparaissent alors entre la catégorie des petites catégories et celle des dynamiques, donnant lieu à des opérations dont nous ignorons si elles pourraient avoir un jour des applications scientifiques, mais dont nous pouvons d'ores et déjà relever le puissant parfum métaphysique : solution existentielle canonique d'une dynamique essentielle, verticalisation des temporalités, essentialisation des dynamiques. Le chapitre 3 se conclut sur deux sections consacrées à l'idée que l'on peut interpréter une dynamique par une autre, la première constituant en quelque sorte une perception réductrice, une projection de la seconde.

Finalement, le chapitre 4 introduit les dynamiques catégoriques connectives. En particulier, à toute dynamique catégorique connective est associé son ordre connectif, un ordinal qui mesure en quelque sorte la complexité connective de la dynamique considérée. 

L'objet de ce texte est essentiellement d'introduire la notion de dynamique catégorique connective et il est évident que beaucoup de questions restent ouvertes dans ce champ immense, à commencer par la question de savoir ce qui du chapitre sur les dynamiques catégoriques ensemblistes peut être transposé dans celui concernant les dynamiques connectives, ainsi que celles relatives à la détermination effective de l'ordre connectif des exemples classiques de dynamiques et à l'obtention de dynamiques d'ordre connectif donné. 

\begin{flushright}
Version 1.0, achevée le 22 décembre 2011, à Paris.
\end{flushright}

PS : On trouvera des indications sur les notations utilisées en fin d'ouvrage, juste avant les références bibliographiques et la table des matières.


\chapter{Espaces connectifs (rappels et compléments)}

Pour une présentation détaillée des espaces connectifs, nous renvoyons à notre article \cite{Dugowson:201012}. Au rappel succinct des notions de base (définition des catégories $\mathbf{Cnct}$ et $\mathbf{Cnc}$, treillis et engendrements de structures, limites et colimites, structures initiales et finales, structure induite, structure quotient...), le présent chapitre ajoute quelques compléments, à commencer par de rapides repères historiques puis la construction du quotient d'un espace connectif par une relation d'équivalence partielle (section \ref{relations equivalence partielle}), la notion de quotient structural (section \ref{quotient structural}), celle d'espace de séparation (section \ref{espaces de separation}). Enfin, la notion d'ordre connectif, présentée dans \cite{Dugowson:201012} uniquement dans le cas des espaces finis, est étendue à tous les espaces connectifs dans la section \ref{ordre connectif}.

\section{Brefs repères historiques}

\`{A} ma connaissance, le premier article décrivant des structures générales de type proprement connectif est \cite{Brunn:1892a}, publié par Hermann Brunn en 1892. Reprises par Debrunner au début des années 1960 (\cite{Debrunner:19600416,Debrunner:1964}) puis par Kanenobu \cite{Kanenobu:198504,Kanenobu:1986} en 1985, les structures connectives considérées par Brunn n'y sont pas vraiment considérées pour elles-mêmes --- elles ne surgissent qu'en relation avec les entrelacs --- elles ne comportent qu'un nombre fini de composantes et les morphismes correspondants ne sont pas définis. Indépendamment et pour la première fois, en 1981 et 1983, Reinhard Börger \cite{Borger:1981,Borger:1983} définit la catégorie des espaces connectifs. Indépendamment encore, en 1988, Georges Matheron et Jean Serra \cite{Matheron.Serra:1988, Matheron.Serra:1988k}  posent une définition identique, mais sans les morphismes. Indépendamment, en 2006, Joseph Muscat et David Buhagiar \cite{Muscat_Buhagiar:2006} introduisent une notion plus restreinte, les \textit{connective spaces}, qu'ils étudient d'un point de vue catégorique. La définition de Muscat et Buhagiar, plus proche des espaces topologiques, est trop restrictive pour ce qui nous intéresse ici puisqu'elle ne permet pas de considérer la structure connective des entrelacs\footnote{Voir plus loin la remarque \ref{rmq Muscat Buhagiar} page \pageref{rmq Muscat Buhagiar}.}.

Au début des années 2000, ignorant à l'époque ces différents travaux, j'ai moi-même posé, à l'occasion d'une réflexion sur la nature topologique ou non du jeu de go, une définition de la catégorie des espaces connectifs intègres équivalente à celle de Börger (\cite{Dugowson:2003, Dugowson:2007c, Dugowson:201012}). Par contre, la considération d'espaces connectifs non intègres, c'est-à-dire la possibilité d'avoir des singletons non connexes, qui joue un rôle important dans la notion de feuilletage connectif, est tout-à-fait nouvelle, de même que le théorème d'engendrement des structures connectives et la mise en évidence de la structure monoïdale fermée de la catégorie des espaces connectifs intègres\footnote{Voir plus loin la remarque \ref{remarque tensoriel connectif}, ainsi que la section 4 de l'article \cite{Dugowson:201012}.},  qui constituent, avec l'histoire du théorème de Brunn-Debrunner-Kanenobu et la notion d'ordre connectif pour les espaces finis, l'essentiel de mon article \cite{Dugowson:201012}.
Sont également nouvelles les principales notions du présent ouvrage, indiquées en introduction : 
la notion d'ordre connectif d'un espace connectif quelconque\footnote{Section \ref{ordre connectif} du présent chapitre.},
les notions en quelque sorte adjointes de feuilletage connectif et de représentation d'un espace connectif dans un autre\footnote{Chapitre \ref{chapitre feuilletages et representations} du présent ouvrage.},
 celles de dynamique catégorique connective, de feuilletage et d'ordre connectif d'une telle dynamique, 
 notions construite sur les concepts des dynamiques catégoriques que nous introduisons au chapitre \ref{chapitre dynamiques categoriques ens}.  

\section{Définition}\label{section definition espaces connectifs}

\begin{df} [Espaces connectifs] Un \emph{espace connectif} est un  couple $X=(\vert X\vert,\kappa(X))$ formé d'un ensemble $\vert X\vert$ et d'un ensemble non vide $\kappa(X)$ de parties de $\vert X\vert$ tel que pour toute famille $\mathcal{I}\in \mathcal{P}(\kappa(X))$, on ait
\begin{displaymath}
\bigcap_{K\in\mathcal{I}}K\ne\emptyset\Longrightarrow \bigcup_{K\in\mathcal{I}}K\in\kappa(X).
\end{displaymath}
\noindent L'ensemble $ \vert X \vert$ est le \emph{support} de $X$, l'ensemble $\kappa(X)$ est la \emph{structure connective} de $X$ et ses élé\-ments sont les \emph{parties connexes} ou \emph{parties connectées} de l'espace connectif $X$. 
Un point $x\in\vert X\vert$ est \emph{absent} s'il n'appartient à aucune partie connexe de $X$, il est \emph{présent} dans le cas contraire. On appelle \emph{composante absente} de $X$ l'ensemble des points absents de $X$. On appelle \emph{composantes connexes} de $X$ les parties connexes maximales pour l'inclusion. Pour tout point présent $x$ de $X$, on appelle composante connexe de $x$ l'unique composante connexe de $X$ contenant $x$.
Nous dirons d'un espace connectif qu'il est \emph{intègre} si tout  singleton est connecté.
Un \emph{morphisme connectif}, ou \emph{application connective}, d'un espace connectif $(\vert X\vert ,\kappa(X))$ vers un autre $(\vert Y\vert,\kappa(Y))$ est une application $f:\vert X\vert \to \vert Y\vert $ telle que :
\begin{displaymath}
\forall K\in \kappa(X), f(K)\in \kappa(Y).
\end{displaymath}
\end{df}

\begin{rmq} Selon cette définition, la partie vide de $ \vert X \vert$ est toujours connexe puisque, dans le cas où $ \vert X \vert$ est non vide, elle est l'union de la famille vide (dont l'intersection est non vide).
\end{rmq}

\begin{rmq}
Si $X$ est un espace connectif non vide sans point absent, en particulier s'il est intègre, les composantes connexes de $X$ constituent une partition de son support $\vert X\vert$.
\end{rmq}

\begin{notations}\label{notation categories espaces connectifs}
On notera $\mathbf{Cnc}$ la catégorie des espaces et morphismes connectifs, et $\mathbf{Cnct}$ la sous-catégorie pleine de $\mathbf{Cnc}$ dont les objets sont les espaces connectifs intègres.
\end{notations}

\section{Quelques exemples}

\begin{exm}[Espaces connectifs topologiques]  On définit un foncteur $U_T:\mathbf{Top}\to\mathbf{Cnct}$ en associant à tout espace topologique l'espace connectif intègre ayant les mêmes points, et dont les connexes sont les parties connexes pour la topologie de l'espace considéré. En particulier, pour tout entier $n$, nous noterons $(E_n,\tau)$ l'espace connectif associé par $U_T$ à l'espace affine $E_n\simeq\mathbf{R}^n$ muni de la topologie usuelle. Nous dirons qu'un espace connectif est topologique s'il est dans l'image objet de $U_T$.
\end{exm}

\begin{exm}[Espaces connectifs graphiques]
La connexité au sens des graphes conduit à la définition d'un foncteur $U_G:\mathbf{Grf}\rightarrow \mathbf{Cnct}$ sur la catégorie $\mathbf{Grf}$ des graphes simples non orientés dont les flèches sont les applications qui préservent cette connexité.  Nous dirons qu'un espace connectif est graphique s'il est dans l'image objet de $U_G$.
\end{exm}

\begin{exm}[Relations d'équi\-va\-lence]
\label{structure connective pour une equivalence} Une relation d'équi\-va\-lence $\sim$ sur un ensemble $E$ peut être vue comme un graphe de sommets dans $E$, comme une structure topologique sur $E$ ou encore comme une structure connective intègre $\kappa_\sim$ sur $E$, à savoir : $K\in \kappa_\sim \Leftrightarrow \forall (x,y)\in K^2, x\sim y$. Autrement dit, une partie non vide de $E$ est connexe pour la structure $\kappa_\sim$ si et seulement si elle est incluse dans une classe d'équi\-va\-lence de $\sim$. Plus généralement, toute relation d'équi\-va\-lence partielle peut être vue comme une structure connective\footnote{Voir plus loin la section \ref{relations equivalence partielle} page \pageref{relations equivalence partielle}}, mais une telle structure n'est pas intègre en général, donc n'est ni graphique ni topologique.
\end{exm}

\begin{exm}[Espace borroméen, espaces brunniens]\label{exemple espace borromeen}
Un espace connectif intègre peut être à la fois topologique et graphique, graphique sans être topologique, topologique sans être graphique\footnote{Cette possibilité n'est réalisée que pour des espaces infinis : sur un ensemble fini de points, toute structure topologique définit une structure connective déterminée par les paires connexes, donc graphique.}, ou encore ni graphique ni topologique. L'exemple le plus simple d'espace intègre ni topologique ni graphique est l'espace borroméen  $\mathcal{B}_3$, où pour tout entier $n$ on désigne par $\mathcal{B}_n$  l'espace intègre à $n$ points dont la seule partie connexe non réduite au vide ou à un singleton est la partie pleine. Pour $n>3$, $\mathcal{B}_n$ est appelé espace brunnien à $n$ points. A noter que la méthode de Newton dans le plan complexe pour les racines cubiques de l'unité fournit un exemple \og naturel \fg\, de morphisme connectif à valeur dans l'espace borroméen $\mathcal{B}_3$ (voir \cite{Dugowson:201012}, exemple 9).
\end{exm}

\begin{exm}[Espaces difféologiques]\label{exm espaces diffeologiques}
Les parties connectées d'un espace difféo\-logique (voir \cite{PIZ:2012}, chapitre 5, en particulier les sections 5.1, 5.6 et 5.9) constituent une structure connective sur l'ensemble des points de cet espace. En outre, toute application lisse préserve la connexité (\cite{PIZ:2012}, section 5.9). Autrement dit, il y a un foncteur (d'oubli) de la catégorie des espaces difféologiques dans celle des espaces connectifs. On vérifie facilement que les espaces $\mathcal{B}_n$ décrits ci-dessus dans l'exemple \ref{exemple espace borromeen} sont difféologisables. Plus généralement :
\begin{exc}
Montrer que tout espace connectif fini est difféologisable.
\end{exc}
\end{exm}

D'autres exemples sont développés dans \cite{Dugowson:201012}.

\section{Treillis et engendrement des structures}

Pour tout ensemble $E$, l'ensemble des structures connectives dont $E$ peut être muni constitue un treillis complet (pour l'inclusion). La structure la plus fine (\textit{i.e.} la plus petite pour l'ordre défini par l'inclusion), pour laquelle seul le vide est connexe, est appelée \emph{structure totalement discrète} ou encore \textit{structure désintégrée}. La moins fine est la structure \emph{grossière}, ou \emph{indiscrète}, pour laquelle toute partie de $E$ est connexe. 

De même, l'ensemble des structures connectives intègres dont $E$ peut être muni constitue un treillis complet pour l'inclusion. La plus fine d'entre elle est appelée \textit{structure discrète}, ou encore, pour éviter tout risque de confusion, \textit{structure discrète intègre}.

Ces treillis de structures admettent une expression fonctorielle, faisant entrer la catégorie des espaces connectifs (resp. intègres) dans le cadre de ce que j'ai appelé les \emph{catégories à treillis de structures} (\cite{Dugowson:201012}, §\,3.1), qui constituent, comme par exemple la catégorie des espaces topologiques, un cas particulier de \textit{catégories topologiques} au sens de Adamek, Herrlich, et Strecker \cite{joycats:2005}. 

Comme dans le cas des espaces topologiques, une conséquence de l'organisation en treillis complets des structures connectives sur un ensemble $E$ est la notion de structure connective la plus fine contenant un ensemble donné quelconque $\mathcal{A}$ de parties de $E$. La construction d'une telle structure, dès lors appelée structure connective \emph{engendrée} par $\mathcal{A}$, est donnée dans \cite{Dugowson:201012}, §\,2. On note \label{notation structure connective engendree} $[\mathcal{A}]_0$ la structure connective engendrée par $\mathcal{A}$, tandis que la structure connective \emph{intègre} engendrée par $\mathcal{A}$ sera notée $[\mathcal{A}]_1$ ou $[\mathcal{A}]$, de sorte que \[[\mathcal{A}]=[\mathcal{A}]_0 \cup \{\{a\}, a\in E\}.\]

\section{Structures initiales, structures finales}

En conséquence du fait que $\mathbf{Cnc}$ et $\mathbf{Cnct}$ sont des catégories à treillis de structures, elles sont complètes et cocomplètes : toutes les limites (produits cartésiens, fibrés, etc.) et colimites (unions disjointes, amalgamées, etc.) existent (voir \cite{Dugowson:201012}, §\,3.2). \`{A} la base de toutes ces constructions on trouve les notions de structures initiales et finales. Soit $f:E\to E'$ une application entre deux ensembles $E$ et $E'$. Si $\mathcal{K}$ est une structure connective sur $E$, la structure finale sur $E'$ associée à $\mathcal{K}$ par $f$, notée $f_!(\mathcal{K})$ est la plus fine des structures connectives sur $E'$ qui fasse de $f$ un morphisme connectif. Elle est donnée par la formule

\begin{equation} \label{expression of f!}
f_!(\mathcal{K})=[\{f(K), K\in \mathcal{K}\}]_0,
\end{equation}

\noindent ou, dans le cas où on se limite aux structures connectives \textit{intègres}, 

\begin{equation} \label{expression of f!}
f_!(\mathcal{K})=[\{f(K), K\in \mathcal{K}\}].
\end{equation}

De même, si $E'$
 est muni d'une structure connective $\mathcal{K}'$,
  la structure la moins fine sur $E$ qui fasse de $f$ un morphisme connectif, appelée structure initiale associée à $\mathcal{K}'$ par $f$, est notée $f^*(\mathcal{K}')$ et vérifie
 
\begin{equation}  \label{expression of f*}
f^*(\mathcal{K}')=\{K\in \mathcal{P}_X, f(K)\in \mathcal{K}'\}.
\end{equation}

La section §\,3.3 de \cite{Dugowson:201012} développe en particulier les notions de quotient d'un espace connectif par une relation d'équi\-va\-lence et de structure connective induite sur une partie d'un espace connectif. \'{E}tant donnée l'importance de ces constructions, nous faisons quelques rappels et compléments à leur sujet dans les sections qui suivent.

\section{Structure connective induite sur une partie}

Dans toute catégorie d'ensembles à treillis de structures, la structure induite sur une partie d'un espace est la structure initiale associée à la structure de cet espace par l'injection canonique de la partie considérée. Dans le cas des espaces connectifs, on obtient la construction suivante : étant donné un espace connectif $X$ et $A$ une partie de l'ensemble $\vert X\vert$, la \textit{structure connective induite} par $X$ sur $A$ est la structure connective initiale $i^*(\kappa(X))$ sur $A$ associée à $\kappa(X)$ par l'injection canonique $i:A\hookrightarrow \vert X\vert$, autrement dit c'est la structure connective la moins fine sur $A$ qui fasse de l'injection canonique $i$ un morphisme connectif, de sorte que, notant $X_{\vert A}$ l'espace induit par $X$ sur $A$, on a
\[ 
\kappa(X_{\vert A})=i^*(\kappa(X))=\{K\in\mathcal{P}_A, A\in\kappa(X)\}=\kappa(X)\cap \mathcal{P}(A).
 \]
 
Il est remarquable que s'agissant des structures induites, il n'y aura pas d'ambiguïté concernant les structures connectives associées aux topologies, du fait du résultat suivant.

\begin{prop} Soit $(\vert X\vert,\tau)$ un espace topologique, $A$ une partie de l'ensemble $\vert X\vert$, $\kappa$ l'ensemble des parties connexes de l'espace topologique $X$, $\tau_{|A}$ la structure topologique induite sur $A$ par $(\vert X\vert,\tau)$ et $\kappa_{|A}$ la structure connective induite sur $A$ par $(\vert X\vert,\kappa)$. Alors $(A,\kappa_{|A})=U_T(A,\tau_{|A})$ : la structure connective associée à la structure topologique induite sur $A$ par $\tau$ coïncide avec la structure connective induite sur $A$ par la structure connective sur $X$ associée à $\tau$.
\end{prop}

\noindent \textbf{Preuve}. Soit $K\in\kappa_{|A}$. On a $K\in\kappa$,  donc $K$ ne peut être séparé par deux ouverts de $X$, au sens où il n'existe pas deux ouverts $U$ et  $V$ de $X$ tels que $U\cap K\neq \emptyset \neq V\cap K$ et $U\cup V \supset K$ et $U \cap K \cap V =\emptyset$. Alors $K$ ne peut pas non plus être séparé par deux ouverts de $(A,\tau_{|A})$, puisque ceux-ci sont la trace sur $A$ des ouverts de $X$, de sorte que $K$ est connexe dans l'espace topologique $(A,\tau_{|A})$. Réciproquement, si $K\in U_T(A,\tau_{|A})$, les mêmes arguments montrent que $K$ est connexe dans $(\vert X\vert,\tau)$, d'où $K\in\kappa_{|A}$.
\begin{flushright}$\square$\end{flushright} 
\pagebreak[3]

\section{Quotient par une relation d'équi\-va\-lence}

Dans toute catégorie d'ensembles à treillis de structures, l'espace quotient $X/\sim$ d'un espace $X$ par une relation d'équi\-va\-lence $\sim$  sur l'ensemble des points de l'espace est défini en munissant l'ensemble quotient $\vert X \vert/\sim$ de la structure finale  associée à celle de $X$ par la surjection canonique. Dans le cas des espaces connectifs, on obtient la construction suivante : étant donné un espace connectif $(\vert X \vert,\kappa(X))$ et une relation d'équi\-va\-lence $\sim$ sur l'ensemble $\vert X \vert$, le support de l'espace quotient est l'ensemble quotient $\vert X \vert/\sim$, et la structure connective de l'espace quotient est la structure connective 
\[[\{s(K), K\in \kappa(X)\}]_0\] 
\noindent \textit{engendrée} par l'ensemble $\{s(K), K\in \kappa(X)\}$ des images par la surjection canonique $s:\vert X \vert \rightarrow \vert X \vert/\sim $ des connexes de $X$.

Contrairement aux structures induites, le foncteur $U_T$ ne respecte pas les structures quotients, de sorte que le point de vue connectif peut de ce point de vue se révéler parfois plus fin que le point de vue topologique. Par exemple, le quotient de l'espace topologique usuel $\mathbf{R}$ par la relation d'équi\-va\-lence
\[ x\sim y \Longleftrightarrow x-y \in \mathbf{Q} \]
conduit à un espace grossier, puisque l'union d'une famille non vide de classes d'équi\-va\-lence ne peut être un ouvert de $\mathbf{R}$ que si elle les contient toutes, tandis que le quotient de l'espace connectif usuel $U_T(\mathbf{R})$ par la même relation d'équi\-va\-lence conduit sur l'ensemble $\mathbf{R}/\sim$ à une structure brunnienne : outre le vide et les singletons, seul l'espace entier est connexe, du fait que tout intervalle de $\mathbf{R}$ non vide et non réduit à un point est d'intérieur non vide et que toutes les classes d'équi\-va\-lence sont denses dans $\mathbf{R}$.

Dans le cas particulier où toutes les classes d'équi\-va\-lence sont connexes, la structure de l'espace quotient peut s'exprimer plus simplement. Pour toute partie $A \subset\vert X\vert$, notons $\hat{A}$ la partie de $\vert X\vert$ définie par 
\begin{equation}\label{notation A chapeau}
\hat{A}=\bigcup_{x\in A} s(x).
\end{equation}

\begin{prop}[Cas de classes connexes]
Si toutes les classes d'équi\-va\-lence sont connexes, la structure connective de l'espace quotient $X/\sim$ est l'ensemble $\{s(\hat{K}), K\in \kappa(X)\}$.
\end{prop}

\noindent \textbf{Preuve}. Pour tout connexe non vide $K$, $\hat{K}=\bigcup_{x\in K} s(x)=K\cup \bigcup_{x\in K}s(x)
=\bigcup_{x\in K} (K\cup s(x))$ 
est l'union de connexes d'intersection non vide, donc est connexe. 
Donc $\{s(\hat{K}), K\in \kappa(X)\}\subset [\{s(K), K\in \kappa(X)\}]_0$. 
Pour prouver l'inclusion réciproque, il suffit de montrer que 
$\{s(\hat{K}), K\in \kappa(X)\}$ 
est une structure connective de l'espace quotient, contenant $\{s(K), K\in \kappa(X)\}$. Par construction, on a trivialement $s(A)=s(\hat{A})$. 
Par ailleurs, si $\mathcal{A}$ est une famille de connexes de $X$ telle que $\bigcap_{K\in\mathcal{A}}s(\hat{K})\neq\emptyset$, 
alors  $\bigcup_{K\in\mathcal{A}}s(\hat{K})=s(\hat{L})$ où
$\hat{L}=L=\bigcup_{K\in\mathcal{A}}\hat{K}$ 
est une partie connexe de $X$ comme union de connexes d'intersection non vide 
(puisque $c\in s(\hat{A})\Rightarrow \emptyset\neq c\subset \hat{A}$).
\begin{flushright}$\square$\end{flushright} 
\pagebreak[3]

\begin{cor}\label{corollaire parties non connexes du quotient} Si toutes les classes d'équi\-va\-lence 
sont connexes, et si $A$ est telle que $\hat{A}
$ soit une partie non connexe de $X$, alors $s(A)
$ est une partie non connexe de l'espace quotient $X/\sim$.
\end{cor}
\noindent \textbf{Preuve}. Si par absurde $s(A)=\{s(x), x\in A\}$ était une partie connexe de $X/\sim$, on aurait $s(A)=s(\hat{K})$ pour une certaine partie connexe $K$ de $X$, d'où $\hat{A}=\hat{\hat{K}}=\hat{K}$ et $\hat{A}$ serait connexe.
\begin{flushright}$\square$\end{flushright} 
\pagebreak[3]

\section{Relations d'équi\-va\-lence partielle}\label{relations equivalence partielle}

Une \textit{relation d'équi\-va\-lence partielle}  sur un ensemble est une relation binaire sur cet ensemble, symétrique et transitive, mais non né\-cessai\-rement réflexive. On appelle \textit{partie présente} d'une telle relation l'ensemble des élé\-ments en relation avec eux-mêmes. Le complémentaire de la partie présente est la \textit{partie absente}.

Une  relation d'équi\-va\-lence partielle est entièrement définie par la donnée de ses classes d'équi\-va\-lence, c'est-à-dire les classes d'équi\-va\-lence de sa restriction à la partie présente.

Comme pour les relations d'équi\-va\-lence\footnote{Voir l'exemple \ref{structure connective pour une equivalence} page \pageref{structure connective pour une equivalence}.}, toute relation d'équi\-va\-lence partielle $\sim$ peut être identifiée à la structure connective $\kappa_\sim$ pour laquelle les parties connexes sont les parties des classes de la relation en question. La différence est que cette fois la structure $\kappa_\sim$ n'est pas intègre en général.

\begin{df}[Saturation d'une structure connective] \'{E}tant donnée une structure connective $\mathcal{K}$ sur un ensemble $E$, on appelle \emph{saturation} de $\mathcal{K}$, et l'on note  $\widetilde{\mathcal{K}}$, la structure connective sur $E$ telle qu'une partie quelconque $A$ de $E$ est $\widetilde{\mathcal{K}}$-connexe si et seulement s'il existe une partie $\mathcal{K}$-connexe $K$ de $E$ telle que $K\supset A$. Par ailleurs, on appelle relation d'équi\-va\-lence partielle \emph{engendrée} par  $\mathcal{K}$ la relation d'équi\-va\-lence partielle ${\chi[\mathcal{K}]}$ dont les classes sont les composantes connexes de $\mathcal{K}$.
\end{df}

Bien entendu, toute structure connective $\mathcal{K}$ sur $E$ a même partie absente que la structure saturée $\widetilde{\mathcal{K}}$. La proposition suivante est également évidente.

\begin{prop} La structure connective $\kappa_{\chi[\mathcal{K}]}$ associée à la relation d'équi\-va\-lence partielle ${\chi[\mathcal{K}]}$ engendrée par $\mathcal{K}$ coïncide avec la saturée $\widetilde{\mathcal{K}}$ de $\mathcal{K}$.
\end{prop}

\begin{df}[quotient par une relation d'équi\-va\-lence partielle]\label{def quotient par une relation partielle}  Soit $X$ un espace connectif, et $\sim$ une relation d'équi\-va\-lence partielle sur $\vert X\vert$. On définit le quotient $X/\sim$ comme étant l'espace connectif quotient $X'/\sim'$, où $X'$ est l'espace connectif induit par $X$ sur la partie présente de $\sim$, et $\sim'$ est la relation d'équi\-va\-lence induite par $\sim$ sur $\vert X'\vert$.
\end{df}

\section{Quotient structural}\label{quotient structural}

%

Dans cet ouvrage, nous aurons en outre besoin d'une construction voisine de celle du quotient d'un espace $X$ par une relation d'équi\-va\-lence $\sim$, mais consistant à ne modifier que la structure connective et non le support de l'espace considéré, de sorte que les points situés dans une même classe d'équi\-va\-lence partagent les mêmes relations connectives, au sens où s'il existe une bijection $\varphi$ entre deux parties $A$ et $B$ de l'espace telle que pour tout $a\in A$, $a\sim\varphi(a)$, alors $A$ est connexe si et seulement si $B$ l'est. Plus précisément, on pose la définition suivante.

\begin{df}[Quotient structural] Soit $X$ un espace connectif et $\sim$ une relation d'équi\-va\-lence sur $\vert X\vert$. Le \emph{quotient structural} de $X$ par $\sim$ est l'espace connectif $\dfrac{X}{\sim}$ défini par 
\[ 
\dfrac{X}{\sim}=(\vert X\vert, s^*(\kappa(X/\sim))),
 \]
où $s$ désigne la surjection canonique $s:\vert X\vert \twoheadrightarrow \vert X/\sim \vert$.
Autrement dit, on a $\vert\dfrac{X}{\sim}\vert=\vert X \vert$ et
$\kappa(\dfrac{X}{\sim})=\{A\in \mathcal{P}_{\vert X \vert}, s(A)\in [s(K), K\in \kappa(X)]_0\}$.
\end{df}

\begin{rmq} Plus généralement, on pourrait définit le quotient structural d'un espace connectif $X$ par une relation d'équivalence partielle $\sim$ comme le quotient structural de l'espace $X'$ par la relation $\sim'$, où $X'$ désigne l'espace induit par $X$ sur la partie présente de $\sim$, et $\sim'$ la restriction à $X'$ de la relation $\sim$. Avec cette définition, le quotient structural d'un espace $X$ par la relation d'équivalence partielle $\sim[\kappa(X)]$ coïncide avec sa partie présente munie de la structure saturée $\widetilde{\kappa(X)}$.


\end{rmq}

\begin{exm} Pour la relation d'équi\-va\-lence $x\sim y \Leftrightarrow x-y\in \mathbf{Q}$, l'espace $\dfrac{\mathbf{R}}{\sim}$
a pour connexes d'une part toute partie de toute classe d'équi\-va\-lence, et d'autre part tout ensemble de réels contenant au moins un représentant de chaque classe, comme par exemple les intervalles d'intérieur non vide.
\end{exm}

\begin{exm}[Sphère ligaturée] \label{Sphere ligaturee} 

Soit $\mathcal{S}_2$  la sphère unité  de l'espace euclidien $\mathbf{R}^3$, muni de la structure connective usuelle (topologique), et soit $\sim$ la relation d'équi\-va\-lence telle que deux points distincts $(x,y,z)$ et $(x',y',z')$ sur la sphère sont équivalents si et seulement s'ils ont même hauteur non rationnelle $z=z'\notin  \mathbf{Q}$. Alors les parties connexes de 
$\dfrac{\mathcal{S}_2}{\sim}$
sont d'une part les arcs de cercles de hauteur constante $z\in\mathbf{Q}$, et d'autre part toute partie $A$ de la sphère dont la projection orthogonale $p(A)$ sur l'axe des $z$ soit un intervalle non réduit à un singleton rationnel. 

En effet, notant $\mathcal{C}_z$ le cercle de hauteur $z$ tracé sur $\mathcal{S}_2$, l'image par la surjection canonique $s:\mathcal{S}_2\rightarrow \mathcal{S}_2/\sim$ d'une telle partie $A$ coïncide avec l'image par $s$ de l'ensemble $A'=A\cup \bigcup_{z\in p(A)\setminus \mathbf{Q}}{\mathcal{C}_z}$. 

Or, $A'$ est une partie connexe de la sphère. Soit en effet $U$ et $V$ deux ouverts disjoints de $\mathcal{S}_2$ recouvrant $A'$. Il n'existe pas de $z\in p(A)$ tel que $U\cap \mathcal{C}_z \neq \emptyset \neq V\cap \mathcal{C}_z$. En effet, cela est impossible si $z\notin \mathbf{Q}$, car chaque $\mathcal{C}_z$ est connexe, tandis que si l'on suppose que $U$ et $V$ rencontrent tous deux $\mathcal{C}_z$ avec $z\in \mathbf{Q}\cap p(A)$, alors $U$ et $V$ devraient également rencontrer un cercle $\mathcal{C}_{z'}$ avec $z'\in  p(A)\setminus \mathbf{Q}$, ce qui est impossible. Donc, la projection $p$ étant une application ouverte, $p(U)$ et $p(V)$ sont des ouverts de $\mathbf{R}$ recouvrant l'intervalle $p(A)$ mais tels que $p(U)\cap p(V)\cap p(A)=\emptyset$, donc $U\cap A$ ou $V\cap A$ est en fait vide. 

Réciproquement, il est clair que $\mathcal{S}_2[\sim]$ ne peut contenir d'autres parties connexes que celles indiquées. 
\end{exm}

\section{Espaces de séparation}\label{espaces de separation}

Cette dernière section est consacrée à une manière particulière de spécifier les espaces connectifs intègres, grâce à ce que j'appelle des dispositifs de séparation. Ce procédé jouera un rôle dans la notion de représentation connective, en particulier dans la représentation par entrelacs (voir l'exemple \ref{representations par entrelacs comme representations} page \pageref{representations par entrelacs comme representations}).

Soit $E$ un ensemble. 

\begin{df}
On appelle \emph{dispositif de séparation} sur $E$ tout ensemble  $\mathcal{S}$ de paires $\{S,T\}$ de parties non vides et disjointes de $E$ : 
\[ S \neq \emptyset \neq T\quad  \mathrm{ et }\quad S\cap T = \emptyset. \]
 Les paires $\{S,T\}$ d'un tel dispositif sont appelés \emph{paires séparatrices}, ou \emph{paires de séparateurs}.
\end{df}

\begin{df} Soit $\mathcal{S}$ un dispositif de séparation sur $E$. On dit qu'une partie $A$ de $E$ est \textit{séparée} par $\mathcal{S}$, et l'on note $(\mathcal{S}:A)$,
s'il existe dans $\mathcal{S}$ une paire de séparateurs $\{S,T\}$ qui recouvre $A$ et dont chaque membre rencontre $A$, 
\[A\subset S\cup T, \quad A\cap S\neq \emptyset \quad\mathrm{et}\quad A\cap T\neq \emptyset.\]
\end{df}

\begin{rmq} Pour tout sous-groupe $G$ du groupe des permutations de $E$ et tout dispositif de sé\-para\-tion $\mathcal{S}$ sur $E$, l'ensemble $G\mathcal{S}=\{\{\varphi(S),\varphi(T)\}, \varphi\in G, (S,T)\in\mathcal{S}\}$  est encore un dispositif de séparation sur $E$. On a alors $(G\mathcal{S}:A)$ si et seulement s'il existe $\varphi\in G$ tel que $(\mathcal{S}:\varphi(A))$.
\end{rmq}

\begin{df} [Espace connectif défini par séparation] Soit $\mathcal{S}$ un dispositif de séparation sur $E$. L'ensemble des parties de $E$ qui ne sont pas séparées par $\mathcal{S}$ constitue une structure connective intègre sur $E$, que nous noterons $\kappa(\mathcal{S})$ :
\[\kappa(\mathcal{S})=\{K\in \mathcal{P}(E), \neg (\mathcal{S}:K)\}.\] 
L'espace connectif $(E,\kappa(\mathcal{S}))$ est \emph{l'espace
 connectif défini sur $E$ par le dispositif de séparation
  $\mathcal{S}$}. On le notera $E[\mathcal{S}]$, de sorte que
   $\vert E[\mathcal{S}] \vert=E$ et 
   $\kappa (E[\mathcal{S}])=\kappa(\mathcal{S})$.
\end{df}

\begin{thm} Tout espace connectif intègre peut être défini par un dispositif de séparation.
\end{thm}

\noindent \textbf{Preuve}. On forme un dispositif de séparation adéquat en prenant tous les couples de parties disjointes non vides $(A,B)$ telles que toute composante connexe de $A\cup B$ soit contenue dans $A$ ou bien dans $B$.
\begin{flushright}$\square$\end{flushright} 
\pagebreak[3]

\begin{exm} Le dispositif de séparation le plus faible sur l'ensemble $E$ est obtenu en prenant $\mathcal{S}=\emptyset$. Aucune partie ne peut alors être séparée par $\mathcal{S}$, et $\kappa(\mathcal{S})$ est la structure connective grossière sur $E$ : l'espace obtenu est totalement connecté.
\end{exm}

\begin{exm} Le dispositif de séparation le plus fort sur l'ensemble $E$ est obtenu en prenant pour $\mathcal{S}$ l'ensemble de toutes les paires de parties non vides disjointes de $E$. Toute partie ayant au moins deux éléments est séparée par $\mathcal{S}$, et $\kappa(\mathcal{S})$ est alors la structure connective intègre la plus fine sur $E$, pour laquelle seuls le vide et les singletons sont connectés.
\end{exm}

\begin{exm} [Foncteur $V_T$]
A tout espace topologique, on associe l'espace connectif défini sur le même ensemble de points en prenant pour dispositif de séparation l'ensemble des paires d'ouverts non vides disjoints. En associant en outre elle-même à toute application continue, on définit ainsi un nouveau foncteur, notons-le $V_T:\mathbf{Top}\to\mathbf{Cnct}$, de la catégorie des espaces topologiques dans celle des espaces connectifs intègres : en effet, une application continue transforme né\-ces\-saire\-ment une partie  non séparable par ouverts disjoints de l'espace de départ en partie non séparables par ouverts disjoints de l'espace d'arrivée (puisque dans le cas contraire l'image réciproque des ouverts séparateurs de l'espace d'arrivée permettrait de séparer la partie considérée dans l'espace de départ).

Remarquons que toute partie connexe au sens topologique est néces\-sai\-re\-ment connexe au sens de cette nouvelle structure connective. Autrement dit, le foncteur $U_T$ est connectivement plus fin que le foncteur $V_T$. Si par exemple on prend pour espace topologique $X$ l'ensemble de points $\{1,2,3\}$ avec pour ouverts non triviaux $\{1,2\}$ et $\{1,3\}$, alors $\{2,3\}$ est une partie non connexe de l'espace $U_T(X)$, mais est une partie connexe de $V_T(X)$. 

Par contre, il est connu que dans un espace métrique la connexité d'une partie est équivalente à l'impossibilité de la séparer par des ouverts \textit{disjoints} :

\begin{prop}Pour tout espace métrique (ou métrisable) $X$, on a $U_T(X)=V_T(X)$.
\end{prop}

\noindent \textbf{Preuve}. Montrons que toute partie non connexe $W$ d'un espace métrique, dont nous noterons $\delta$ la distance, est effectivement séparable par deux ouverts \textit{disjoints}. Par définition, il existe deux ouverts $A$ et $B$ de $E$ vérifiant 
\[ 
\left\lbrace 
\begin{array}{l}
W\subseteq A \cup B\\
A\cap B \cap W =\emptyset\\
A\cap W \neq \emptyset\\
B\cap W \neq \emptyset
\end{array}
\right.
 \]
Désignant par $\mathcal{B}(a,\rho)$ la boule ouverte de centre $a$ et de rayon $\rho$, l'ensemble $\{\rho\in \mathbf{R}_+, \mathcal{B}(a,\rho)\subseteq A\}$ est, pour tout $a\in A\cap W$, de la forme $[0,R_a]$, avec $R_a$ un réel strictement positif. On pose de même, pour tout $b\in B\cap W$, $S_b=\max\{\rho\in \mathbf{R}_+, \mathcal{B}(b,\rho)\subseteq B\}$. Pour tout couple $(a,b)\in (A\cap W)\times (B\cap W)$, on a $\delta(a,b)\geq R_a$ et $\delta(a,b)\geq S_b$, donc $\mathcal{B}(a,\frac{1}{2}R_a)\cap \mathcal{B}(b,\frac{1}{2}S_b)=\emptyset$, de sorte que les ouverts $A'=\bigcup_{a\in A\cap W}\mathcal{B}(a,\frac{1}{2}R_a)$ et $B'=\bigcup_{b\in B\cap W}\mathcal{B}(b,\frac{1}{2}S_b)$ sont disjoints et séparent $W$ :
\[ 
\left\lbrace 
\begin{array}{l}
W\subseteq A' \cup B'\\
A'\cap B' =\emptyset\\
A'\cap W \neq \emptyset\\
B'\cap W \neq \emptyset
\end{array}
\right.
 \]
\begin{flushright}$\square$\end{flushright} 
\pagebreak[3]

\end{exm}

\begin{exm} \label{espace affine de separation} Tout espace affine réel $E_n$ de dimension $n$ est muni d'une structure connective notée $\sigma_n$ ou $\sigma$, et appelée la \textit{structure connective usuelle de séparation} sur $E_n$ en prenant pour dispositif de séparation $G\mathcal{S}$ avec $G$ le groupe des homéo\-mor\-phismes de l'espace topologique $E_n \simeq\mathbf{R}^n$ et pour $\mathcal{S}$ le singleton $\{\{S,T\}\}$, avec $S$ et $T$ les deux demi-espaces ouverts définis par un hyperplan quelconque. 

\begin{df} L'espace connectif $(E_n,\sigma)$ est appelé l'espace usuel de séparation $n$-dimensionnelle.
\end{df}

A noter que $(E_n,\sigma)$ est un espace connectif moins fin (il a plus de connexes) que l'espace connectif $(E_n,\tau)$ associé par $U_T$ à l'espace topologique usuel $E_n$.

On munit de même la sphère $n$-dimensionnelle réelle d'une structure connective, notée également $\sigma$, en prenant pour $G$ le groupe des homéo\-mor\-phismes de l'espace topologique $S^n$ et pour paire séparatrice de base les deux demi-espaces séparés par une sphère $(n-1)$-dimensionnelle plongée dans $S^n$.

\begin{df} L'espace connectif $(S^n,\sigma)$ est appelé la sphère usuelle de séparation $n$-dimensionnelle.
\end{df}

Notons que d'autres structures connectives, moins fines, sont obtenues en remplaçant les homéo\-mor\-phismes par des difféomorphismes de classe $C^k$, ou en remplaçant les demi-espaces séparés par un hyperplan par des demi-espaces séparés par un voisinage tubulaire d'un tel hyperplan.

\begin{prop}
La structure connective de l'espace usuel de séparation $(E_n,\sigma)$ n'est pas celle d'un espace topologique.
\end{prop}

\noindent \textbf{Preuve}. On vérifie facilement que dans tout espace topologique, si $A$ et $B$ sont deux parties connexes non vides et que $x$ est un point de l'espace tel que $x \notin A\cup B$ et que $A\cup\{x\}$ et $B\cup\{x\}$ soient non connexes, alors $A\cup B \cup\{x\}$ est encore non connexe. Or, dans $(E_n,\sigma)$, si l'on prend par exemple pour $A$ une demi-sphère,  pour $B$ la demi-sphère complémentaire et pour $x$ le centre de la sphère $A\cup B$, la propriété précédente est contredite.
\begin{flushright}$\square$\end{flushright} 
\pagebreak[3]
\end{exm}

\begin{rmq}\label{rmq Muscat Buhagiar} Cette démonstration repose sur le fait qu'une certaine propriété des structures connectives topologiques n'est pas vérifiée par la structure connective considérée ici. Cette propriété fait partie de celles incorporées par Muscat et Buhagiar dans leur définition des \textit{connective spaces} \cite{Muscat_Buhagiar:2006}. Ainsi, en excluant les structures associées à la séparation par les hyperplans (transformés par les homéo\-mor\-phismes de l'espace) les espaces de Muscat et Buhagiar ne permettent pas de rendre compte de la structure connective des entrelacs.
\end{rmq}

\section{Ordre d'un espace connectif}\label{ordre connectif}


On rappelle qu'un ordinal est un ensemble transitif bien ordonné par l'appartenance. Si $\alpha$ et $\beta$ sont des ordinaux, on a l'équivalence
\[ 
\beta\in\alpha \Leftrightarrow \beta \subsetneq \alpha \Leftrightarrow \beta< \alpha
 \]

Les premiers ordinaux sont $0=\emptyset$, $1=\{0\}$, $2=\{0,1\}$, etc.  On note $Ord$ la classe des ordinaux. On désignera par $\omega_0$ ou $\aleph_0$ le plus petit ordinal infini, c'est-à-dire l'ensemble des ordinaux finis, et par $\aleph_1$ le plus petit ordinal non dénombrable, c'est-à-dire l'ensemble des ordinaux dénombrables.

Pour tout ordinal $\alpha$, nous notons $\alpha^-$ l'ordinal défini par
\begin{itemize}
\item si $\alpha$ a un prédécesseur $\beta$, $\alpha^-=\beta$,
\item si $\alpha$ n'a pas de prédécesseur, $\alpha^-=\alpha$.
\end{itemize}

\begin{df} Soit $\alpha\in Ord$ un ordinal. Un ensemble (partiellement) ordonné $(R,\preceq)$ est dit \emph{supérieur ou égal} à $\alpha$, et l'on note $\alpha\leq R$, s'il existe une application strictement croissante de l'ensemble ordonné $\alpha$ dans l'ensemble ordonné $(R,\preceq)$.
\end{df}

\begin{rmq} La définition précédente est compatible avec la relation d'ordre entre ordinaux : $\alpha\leq\beta$  si et seulement si $\alpha\in\{0, ..., \beta\}=\beta+1$.
\end{rmq}

Soit maintenant $(R,\preceq)$ un ensemble ordonné. La classe des ordinaux $\alpha$ tels que $\alpha\leq R$ est bornée (en fonction du cardinal de $R$), c'est donc un ensemble, et c'est un ordinal puisque $\alpha\leq R\Rightarrow \beta\leq R$ pour tout $\beta\leq\alpha$. 

\begin{df} On appelle \emph{hauteur} de l'ensemble partiellement ordonné ${(R,\preceq)}$, et on note $\Gamma(R)$, l'ordinal
\[\Gamma(R)=\{\alpha\in Ord, \alpha\leq R\} \]
\end{df}

\begin{exm} Soit $(\mathbf{R},\leqslant)$ 
la droite réelle munie de l'ordre usuel. 
Alors $\Gamma(\mathbf{R})=\aleph_1$. 
En effet, on a d'un coté $\Gamma(\mathbf{R})\leqslant\aleph_1$
car tout ordinal $\alpha$ qui appartient à $\Gamma(\mathbf{R})$
est nécessairement dénombrable : dans l'hypothèse contraire, l'existence d'une injection strictement croissante 
$\phi:\alpha\rightarrow \mathbf{R}$ 
entrainerait celle d'une famille non dénombrable ${(]\phi_i,\phi_{i+1}[)}_{i\in\alpha}$ d'intervalles ouverts non vides disjoints de $\mathbf{R}$, ce qui ne se peut. 
D'un autre coté, montrons que $\Gamma(\mathbf{R})\geqslant\aleph_1$. Il suffit pour cela d'établir que $\Gamma(\mathbf{R})$ n'est pas dénombrable. Supposons le contraire. Dans cette hypothèse, il existe une bijection $\Psi:\mathbf{N}\rightarrow\Gamma(\mathbf{R})$. Pour tout entier $n\in\mathbf{N}$, $\Psi_n\leq \mathbf{R}$, 
d'où, par une composition évidente, $\Psi_n\leq [n,n+1[$,
 autrement dit il existe une application strictement croissante 
 $\phi_n:\Psi_n\rightarrow[n,n+1[$. 
Alors l'ensemble $A=\bigcup_{n\in\mathbf{N}}{\phi_n(\Psi_n)}$ est une partie de $\mathbf{R}$ bien ordonnée 
--- car pour toute partie non vide $B$ de $A$, il existe un plus petit entier $n$ tel que $B\cap [n,n+1[ \neq \emptyset$, et $B\cap [n,n+1[$ admet à son tour un plus petit élément $b$ puisque $B\cap [n,n+1[$ 
est une partie non vide de l'ensemble bien ordonné $\phi_n(\Psi_n)$, et $b$ est alors le plus petit élément de $B$ --- 
donc $A$ est en bijection croissante avec un ordinal $\gamma\leq \mathbf{R}$, qui vérifie donc $\gamma\in \Gamma(\mathbf{R})$. Or, par construction, cet ordinal est strictement supérieur  à tout élément de $\Gamma(\mathbf{R})$, puisque par exemple on définit facilement, pour tout $n$, une application strictement croissante $\Psi_n+1\rightarrow \gamma$. Donc $\gamma$ est strictement supérieur à lui-même, ce qui est absurde.

\end{exm}


\begin{df}[Connexes irréductibles]  Soit $X=(\vert X\vert,\kappa_X)$ un espace connectif. Une partie connexe $K\in\kappa_X$ est dite \emph{irréductible} si et seulement si elle n'appartient pas à la structure connective engendrée par les autres parties : $K\notin [ \kappa_X\setminus \{K\} ]_0$. 
\end{df}

\begin{notations} Pour tout espace connectif $X$, on note $(G_X,\subset)$ l'ensemble ordonné par l'inclusion des parties connexes irréductibles de $X$.
\end{notations}

\begin{rmq} Un morphisme connectif transformant toute partie connexe irréductible en une partie connexe irréductible est dit \emph{distingué}. Un espace connectif dont toutes les parties connexes sont irréductibles est dit \emph{distingué}. 
\end{rmq}

\begin{rmq} La partie vide n'est pas irréductible, tandis que les singletons connexes le sont.
\end{rmq}

\begin{rmq} [Parties connexes irréductibles d'un espace fini] Soit $X$ un espace connectif fini. Une partie connexe non vide $K\in\kappa(X)$ est irréductible si et seulement s'il n'existe pas deux parties propres connexes $A\subsetneqq K$ et $B\subsetneqq K$  telles que
\begin{displaymath}
K=A\cup B \textrm{ et }
A\cap B \neq \emptyset.
\end{displaymath}
\end{rmq}


\begin{df}[Ordre connectif d'un espace] Soit $X$ un espace connectif. On appelle \emph{ordre connectif} de $X$ l'ordinal $\Omega(X)=\Gamma(G_X)^{--}$.  
\end{df}

\begin{prop} $\Omega(X)=\{\alpha\in Ord, \alpha+2\leq G_X\}$. 
\end{prop}

\begin{exm} 
Pour un espace $X$ sans connexe irréductible, $G_X$ n'est supérieur qu'à $0$, de sorte que $\Gamma(G_X)=1$ et $\Omega_0(X)=\emptyset=0$. C'est par exemple le cas de la droite réelle lorsqu'on la munit de la structure connective pour laquelle les connexes non vides sont les intervalles \textit{non réduits à un point}.
\end{exm}

\begin{exm} Pour un espace $X$ comportant au moins un connexe irréductible mais  aucun couple de connexes irréductibles emboîtés. Alors $\Gamma(G(X))=2$, et $\Omega_0(X)=\emptyset=0$. C'est par exemple le cas de la droite réelle munie de la structure connective topologique usuelle, qui admet les singletons pour seuls irréductibles.
\end{exm}

\begin{exm} L'ordre connectif $\Omega(X)$ d'un espace connectif fini intègre $X$ coïncide avec l'ordre connectif défini dans \cite{Dugowson:201012}, à savoir la hauteur du graphe orienté acyclique $G_X$ constitué des connexes irréductibles (points génériques) munis de la relation d'inclusion. Par exemple, un graphe fini comportant au moins une arête est d'ordre  $1$, de même l'espace borroméen. A noter qu'un espace connectif fini est entièrement caractérisé par la donnée de ses connexes irréductibles.
\end{exm}

\begin{exm}
On munit l'intervalle $]0,1[$ de la structure connective pour laquelle les connexes sont exactement les intervalles de la forme $]0,x[$, avec $x\in[0,1]$. Il y a une infinité continue de connexes emboîtés, et ils sont tous irréductibles. On a $\Gamma(G_X)=\Gamma(]0,1[)=\Gamma(\mathbf{R})=\aleph_1$, d'où  $\Omega(X)=\aleph_1$.
\end{exm}

\begin{exm} Soit $(X_i)_{i\in \mathbf{N}}$ une famille d'espace connectifs intègres finis tel que tout espace connectif intègre fini soit isomorphe à l'un et un seul des $X_i$ (une telle famille existe puisque pour tout cardinal fini donné il n'y a qu'un nombre fini d'espaces connectifs (à isomorphismes près) ayant ce cardinal pour nombre de points), et soit $Z$ l'union disjointe de tous les $X_i$. Alors $Z$ est un espace connectif intègre, non connexe, dénombrable, d'ordre connectif $\Omega(Z)=\omega_0$.
\end{exm}

\begin{exm} On munit l'ensemble $\{0,1,2\}^{\mathbf{Z}}$, 
c'est-à-dire l'ensemble de toutes les suites $u=(u_n)_{n\in \mathbf{Z}}$
 à valeur dans $\{0,1,2\}$, 
 de la structure connective pour laquelle les connexes sont de la forme $K_{(k,v)}=\{u\in \{0,1,2\}^{\mathbf{Z}}, \forall n\leq k, u_n=v_n\}$
  avec $k\in \mathbf{Z}$ et $v\in \{0,1,2\}^{\mathbf{Z}}$. 
  Cet espace est non connexe et non intègre, il est d'ordre connectif $\omega_0$.
\end{exm}

\begin{exm} Pour tout ordinal $\alpha$, il existe un espace connectif
 d'ordre  $\alpha$. Un exemple  d'un tel espace est obtenu en munissant l'ensemble $\alpha+2$ de la structure connective pour laquelle
   les connexes non vides sont les sections commençantes $k+1=\{0,1,...,k\}$ avec $k\in\alpha$. 
\end{exm}

\chapter{Feuilletages et représentations}\label{chapitre feuilletages et representations}

\section{Représentations connectives}

\subsection{Théorème de Brunn-Debrunner-Kanenobu}

L'espace borroméen $\mathcal{B}_3$ (voir exemple \ref{exemple espace borromeen} page \pageref{exemple espace borromeen}), plus simple des espaces connectifs intègres dont la structure ne soit ni celle d'un espace topologique ni celle d'un graphe, tient son nom de la possibilité de le représenter par l'entrelacs borroméen. Plus généralement, on peut se poser la question de savoir si tout espace connectif fini peut être représenté par un entrelacs dans $\mathbf{R}^3$. Il s'avère que la question avait déjà été posée et en partie résolue dès 1892 par Hermann Brunn \cite{Brunn:1892a, Brunn:1982a, Brunn:1982b}, sans que le concept d'espace connectif (fini) ait toutefois été clairement dégagé par lui : d'une certaine façon, les structures connectives considérées par Brunn restent attachées aux entrelacs. Quoi qu'il en soit, pour Brunn, toute structure connective finie peut effectivement être représentée par un entrelacs, et il donne une idée de preuve, utilisant comme \og briques élémentaires\fg\, de sa construction certains entrelacs particuliers, appelés aujourd'hui, depuis Rolfsen \cite{Rolfsen:1976}, les \og entrelacs brunniens \fg. En 1964, 
Debrunner \cite{Debrunner:1964} affirme que la preuve de Brunn est insuffisante et il en propose une autre mais valable seulement pour les entrelacs $n$-dimensionnels plongés dans l'espace de dimension $n+2$ avec $n\geq 2$. En 1985 et 1986, Kanenobu \cite{Kanenobu:198504,Kanenobu:1986} publie finalement deux preuves de la possibilité de représenter toute structure connective finie par un entrelacs classique. L'idée essentielle de ces différentes construction  se trouve déjà dans l'article de Brunn  : plusieurs entrelacs peuvent toujours être globalement connectés en nouant leurs brins à la façon des \og entrelacs brunniens \fg.

\begin{thm} [Brunn-Debrunner-Kanenobu] Toute structure connective finie est la structure connective d'un entrelacs plongé dans $\mathbf{R}^3$ (ou dans $S^3$).
\end{thm}

On peut alors se demander ce qu'il en est de la représentabilité par entrelacs des espaces connectifs infinis. Il est clair qu'un espace connectif possédant davantage de parties connexes qu'il n'y a de courbes fermées dans $\mathbf{R}^3$, autrement dit un espace connectif dont l'ensemble des parties connexes a un cardinal strictement supérieur à la puissance du continu, ne saurait être représentable par entrelacs. Mais qu'en est-il pour les structures connectives de cardinal inférieur ou égal au continu ? Avant de prétendre aborder cette question, ce que nous ne ferons pas dans le présent ouvrage, il est nécessaire de préciser ce que l'on entend par \og représentation par entrelacs\fg : dira-t-on, par exemple, que le cylindre engendré par le déplacement d'un cercle dans une direction extérieure au plan qui le contient constitue une représentation \textit{par entrelacs} d'un intervalle réel ? Les cercles obtenus ne sont-ils pas, en effet, davantage \og collés\fg\, par continuité qu'\textit{entrelacés} ? Le point de vue que nous soutenons ici est que la représentation d'un espace connectif fini par entrelacs doit être comprise comme un cas particulier de représentation d'un espace connectif dans un autre. C'est cette dernière notion qui constitue l'objet de la présente section.

\subsection{L'espace des parties}

L'idée essentielle de la représentation connective est d'associer à tout \textit{point} de l'espace connectif représenté une \textit{partie} non vide de l'espace connectif dans lequel a lieu la représentation. D'où l'appel aux foncteurs \textit{puissances}, avec une variante selon que l'on souhaite  ou non représenter uniquement des espaces connectifs intègres. Définissons d'abord le foncteur \emph{puissance connective générale}, que nous noterons $\mathcal{P}^*$ comme le foncteur ensembliste des parties non vide dont il constitue en quelque sorte un prolongement aux espaces connectifs.

\begin{df} On définit un endofoncteur $\mathcal{P}^*$ de la catégorie $\mathbf{Cnc}$, appelé \emph{puissance connective générale} ou \emph{espace connectif des parties non vides}, en associant à tout espace connectif $X$ l'espace connectif, noté $\mathcal{P}^*_X$ (ou $\mathcal{P}^*(X)$, ou $\mathcal{P}^*X$) défini par
\begin{itemize}
\item son support $\vert \mathcal{P}^*_X\vert =\mathcal{P}^*_{\vert X\vert}$,
\item et sa structure connective $\kappa(\mathcal{P}^*_X) =\{\mathcal{A}\in \mathcal{P}\mathcal{P}^*_{\vert X\vert}, \bigcup \mathcal{A} \in \kappa(X)\}$,
\end{itemize}
et en associant à tout morphisme connectif $f:X\rightarrow Y$, le morphisme connectif, encore noté $f$, défini pour toute partie non vide $A$ de $\vert X\vert$ par $f(A)=\{f(a), a\in A\}$.  
\end{df}

\begin{rmq} L'espace connectif $\mathcal{P}^*X$ n'est pas intègre en général, même lorsque $X$ l'est. Plus précisément, $\mathcal{P}^*X$ est un espace intègre \textit{uniquement} si $X$ est un espace grossier. Nous verrons d'autres situations, en particulier en relation avec les feuilletages, où les espaces connectifs non intègres apparaissent assez \og naturellement\fg.
\end{rmq}

\subsection{L'espace des parties connexes}

\begin{df} On définit un endofoncteur $\mathcal{K}^*$ de la catégorie $\mathbf{Cnct}$, appelé \emph{puissance connective intègre} ou \emph{espace des parties connexes non vides}, en associant à tout espace connectif intègre $X$ l'espace connectif $\mathcal{K}^*_X$ défini par
\begin{itemize}
\item son support $\vert \mathcal{K}^*_X\vert =\kappa(X)\setminus \{\emptyset\}$,
\item et sa structure connective $\kappa(\mathcal{K}^*_X) =\{\mathcal{A}\in \mathcal{P}(\vert \mathcal{K}^*_X\vert), \bigcup \mathcal{A} \in \kappa(X)\}$,
\end{itemize}
et en associant à tout morphisme connectif $f:X\rightarrow Y$, le morphisme connectif, encore noté $f$, défini pour toute partie connexe non vide $A$ de $\vert X\vert$ par $f(A)=\{f(a), a\in A\}$.  
\end{df}

\begin{rmq} En appliquant cette définition aux espaces connectifs non né\-ces\-saire\-ment intègres, on prolonge naturellement l'endofoncteur $\mathcal{K}^*$ en un foncteur $\mathbf{Cnc}\to \mathbf{Cnct}$ puis en un endofoncteur $\mathbf{Cnc}\to \mathbf{Cnc}$. Ces divers foncteurs seront tous notés $\mathcal{K}^*$.
\end{rmq}



\subsection{Représentation connective}

\begin{df} [Représentation connective] On appelle \emph{représentation} connective d'un espace connectif $X$ dans un espace connectif $Y$ tout morphisme connectif de $X$ dans l'espace connectif $\mathcal{P}^*(Y)$. On écrira  $\rho:X\rightsquigarrow Y$ pour exprimer que $\rho$ est une représentation de $X$ dans $Y$. Etant donnée $\rho$ une telle représentation, $Y$ sera appelé l'\emph{espace} de $\rho$, et sera noté $Y=sp(\rho)$; $X$ sera appelé l'\emph{objet} de $\rho$, et sera noté $ob(\rho)$.
\end{df}

Dans le cas où $X$ est intègre, une représentation $\rho:X\rightsquigarrow Y$ s'identifie à un morphisme connectif de $X$ dans l'espace intègre $\mathcal{K}^*(Y)$.

\begin{df} On dit qu'une représentation $f:X\rightsquigarrow Y$ est \emph{intègre} si son objet et son espace sont tous deux intègres.
\end{df}

\subsection{Composition des représentations}\label{composition des representations}

Soit $\epsilon$ la transformation naturelle $Id_{\mathbf{Cnc}}\rightarrow \mathcal{P}^*$ définie pour tout espace connectif $X$ par $\forall x\in \vert X\vert, \epsilon_X(x)=\{x\}$, et $\mu$ la transformation naturelle $\mathcal{Q}^*\rightarrow \mathcal{P}^*$ définie par $\forall \mathcal{A}\in \mathcal{Q}^*_{\vert X\vert}, \mu_X(\mathcal{A})=\bigcup\mathcal{A}$. Le triplet $(\mathcal{P}^*,\epsilon,\mu)$ constitue alors une monade sur $\mathbf{Cnc}$. La catégorie de Kleisli associée à cette monade a pour objets les espaces connectifs, et pour morphismes les re\-pré\-sen\-ta\-tions, dont la composition est définie pour $\rho:X\rightsquigarrow Y$ et $\tau:Y\rightsquigarrow Z$ par
\[ 
\begin{array}{c}
\tau\odot \rho :X\rightsquigarrow Z \\ 
x\mapsto \mu_Z(\tau^\mathcal{P} (\rho(x))).
\end{array} 
 \]

Étant donnée une représentation $\rho:X\rightsquigarrow Y$, on notera 
${^\mu{\rho}}$ l'application de $\mathcal{P}^*X$ dans $\mathcal{P}^*Y$ définie par ${^\mu{\rho}}=\mu_Y\circ{\rho^\mathcal{P}}$. Une représentation de $X$ dans $Y$ est donc une application $\rho$ de $X$ dans $\mathcal{P}^*_Y$ telle que $^\mu \rho$ transforme toute partie connexe non vide de $X$ en une partie connexe non vide de $Y$.

Pour toute partie non vide $A$ de $X$, on a donc
\[ 
{^\mu{\rho}}(A)=\mu_Y(\rho^\mathcal{P}(A))=\mu_Y(\{\rho(a), a\in A\})=\bigcup_{a\in A}{\rho(a)}\subset Y
\]
tandis que la composée de deux représentations s'écrit
\[ 
\tau\odot \rho = {^\mu{\tau}}\circ \rho.
\]

On définit une catégorie de Kleisli analogue pour les espaces connectifs intègres : notant encore  $\epsilon$ la transformation naturelle $Id_{\mathbf{Cnct}}\rightarrow \mathcal{K}^*$ définie pour tout espace connectif intègre $X$ par $\forall x\in \vert X\vert, \epsilon_X(x)=\{x\}$, et $\mu$ la transformation naturelle $\mathcal{K}^*\mathcal{K}^*\rightarrow \mathcal{K}^*$ définie par $\forall \mathcal{A}\in \vert \mathcal{K}^*\mathcal{K}^* X\vert, \mu_X(\mathcal{A})=\bigcup\mathcal{A}$, le triplet $(\mathcal{K}^*,\epsilon,\mu)$ constitue une monade sur $\mathbf{Cnct}$, dont la catégorie de Kleisli associée a pour objets les espaces connectifs intègres, et pour morphismes les re\-pré\-sen\-ta\-tions intègres, avec la composition des re\-pré\-sen\-ta\-tions définie comme pour le cas général.

%





\subsection{Représentations claires, représentations distinctes}

\begin{df}
Soit $\rho:X\rightsquigarrow Y$ une représentation d'un espace $X$ dans un espace $Y$. On dit que $\rho$ est \emph{claire} si 
\[\forall A\in \mathcal{P}_{\vert X\vert}, A\notin \kappa(X)\Rightarrow ^\mu \rho(A)\notin \kappa(Y).\]

On dit que $\rho$ est \emph{distincte} si
\[ \forall(x,y)\in X^2, x\neq y\Rightarrow \rho(x)\cap \rho(y)=\emptyset. \]
\end{df}

\begin{exm} Soit $X$ l'espace connectif non intègre de support $\vert X\vert=\{a,b,c\}$ et de structure connective $\kappa(X)=\{\emptyset, \{a\},\{b\},\{a,b,c\}\}$. Une représentation claire mais non distincte de $X$ dans la droite connective usuelle $(\mathbf{R},\tau)$ est donnée par 

\begin{itemize}
\item $\rho(a)=[0,1[$,
\item $\rho(b)=]1,2]$,
\item $\rho(c)=\{0,1,2\}$.
\end{itemize}

On obtient une représentation claire et distincte de l'espace non intègre $X$ en remplaçant ci-dessus les intervalles par leur intérieur.
\end{exm}

\begin{exm} Une représentation claire et distincte de l'espace borroméen $\mathcal{B}_3$ est obtenue en associant à chacun des trois éléments de cet espace une des trois composantes d'un \textit{noeud borroméen} plongé dans l'espace tridimensionnel usuel de séparation $(E_3,\sigma_3)$. Plus généralement, les entrelacs brunniens constituent, dans l'espace de séparation $(E_3,\sigma_3)$, des représentations claires et distinctes  des espaces connectifs brunniens. 

A noter qu'en supprimant un point sur chacune des composantes des entrelacs considérés on obtient encore des représentations claires et distinctes des espaces brunniens.
\end{exm}

\begin{thm} \label{representation canonique} Tout espace connectif admet une représentation claire et distincte dans un espace intègre. En particulier, tout espace connectif fini admet une représentation par entrelacs, les points non connexes étant représentés par deux ou plusieurs composantes séparables de tels entrelaces.
\end{thm}

\noindent \textbf{Preuve}. Soit $X$ un espace connectif. On pose $\vert X'\vert=\vert X\vert \times \{0,1\}\simeq \vert X\vert \sqcup \vert X\vert$, et soit $\rho:\vert X\vert\to\mathcal{P}(\vert X'\vert)$ l'application définie par 
\[ 
\rho(x)=\{(x,0),(x,1)\}.
 \]
On munit l'ensemble $\vert X'\vert$ de l'unique structure connective intègre dont les connexes non triviaux sont donnés par
\[ \label{bullet}
 \kappa(X')^\bullet=\{^\mu\rho(K), K\in \kappa(X)\}.
\]
Alors $\rho$ est une représentation claire et distincte de $X$ dans l'espace intègre $X'$. 

En particulier, tout espace connectif fini admet une représentation claire et distincte dans un espace connectif fini intègre, qui admet à son tour une représentation claire et distincte par entrelacs (Brunn-Debrunner-Kanenobu). En composant ces deux représentations, on obtient le résultat annoncé
\begin{flushright}$\square$\end{flushright} 
\pagebreak[3]


\begin{df} [Représentations de type $\mathcal{S}$]
Soit  $\mathcal{S}$ un dispositif de séparation sur un ensemble $Y$. On appelle \emph{représentation de type $\mathcal{S}$} toute représentation claire et distincte d'un espace connectif $X$ dans l'espace $Y[\mathcal{S}]$.
\end{df}

\begin{exm}\label{representations par entrelacs comme representations} Soit $\mathcal{S}$ un dispositif de séparation engendrant l'espace de séparation usuel $(E_3,\sigma_3)$. Une représentation par entrelacs d'un espace connectif fini est une représentation de type $\mathcal{S}$.
\end{exm}

\subsection{Catégories de représentations}

On a vu (section \ref{composition des representations} ci-dessus) que les représentations étaient les morphismes entre espaces connectifs dans certaines catégories de Kleisli. Mais l'on peut prendre à leur tour les représentations comme objets, dès lors qu'on aura défini les morphismes entre représentations. 

\begin{df}\label{categorie des representations et rcd} On définit une catégorie $\mathbf{RC}$, dite \emph{catégorie des re\-pré\-sen\-ta\-tions connectives} en prenant pour objets les re\-pré\-sen\-ta\-tions connectives, et pour morphismes  d'une re\-pré\-sen\-ta\-tion $\rho:A\rightsquigarrow B$ vers une re\-pré\-sen\-ta\-tion $\rho':A'\rightsquigarrow B'$ les couples $(\alpha,\beta)$ où $\alpha:A\to A'$ et $\beta:B\to B'$ sont des morphismes connectifs tels que
\[ 
\beta^{\mathcal{P}}\circ \rho \subset \rho'\circ\alpha,
 \]
 au sens où, pour tout $a\in A$, $\beta^{\mathcal{P}}(\rho(a))\subset \rho'(\alpha(a))$.
 
 La sous-catégorie pleine de $\mathbf{RC}$ admettant pour objets les représentations claires et distinctes sera notée $\mathbf{RCD}$.
\end{df}

\begin{exm}[points d'une représentation] La catégorie $\mathbf{RC}$ admet comme objet final l'unique représentation $\mathbf{1}_\mathbf{RC}:\bullet\mapsto\{\bullet\}$ d'un singleton connecté dans lui-même. Un point d'une représentation $\rho:A\rightsquigarrow B$ est alors un morphisme $\mathbf{1}_\mathbf{RC}\to\rho$, c'est-à-dire la donnée d'un point connecté $p$ de $A$ et d'un point connecté $q$ de $\rho(p)\subset B$. En particulier, si l'objet ou l'espace d'une représentation ne possède pas de point intègre, celle-ci n'a pas de point.
\end{exm}

\begin{exm} On définit un foncteur $RC:\mathbf{Cnc}\to\mathbf{RC}$, appelé \textit{représentation canonique}, en associant à tout espace connectif $X$ sa représentation canonique dans un espace intègre $X'$ tel que $\vert X'\vert=\vert X\vert \times \{0,1\}$ (voir le théorème \ref{representation canonique} page \pageref{representation canonique}), et à tout morphisme d'espaces connectifs $f:X\to Y$ le morphisme $(\alpha,\beta):X'\to Y'$, où $\alpha=f$ et $\beta((x,i))=(f(x),i)$ pour tout $x\in X$ et $i\in \{0,1\}$.
\end{exm}

\section{Feuilletages connectifs}


\begin{df} [Feuilletage connectif]
Un \emph{feuilletage connectif}  est un triplet $(E,\kappa_0, \kappa_1)$ constitué d'un ensemble $E$ appelé le \emph{support} du feuilletage, et d'un couple $(\kappa_0, \kappa_1)$ de structures connectives sur $E$, la première, $\kappa_0$, étant dite structure connective \emph{interne}, et la seconde, $\kappa_1$, structure connective \emph{externe}. Lorsque que $\kappa_0\subset\kappa_1$, le feuilletage est dit \emph{régulier}.
\end{df}

Lorsqu'une partie de $E$ est connexe pour $\kappa_0$ (resp. $\kappa_1$), on dit aussi qu'elle est $\kappa_0$-connexe, ou encore qu'elle est \textit{connexe interne}, ou encore \textit{intérieurement connexe} (resp. $\kappa_1$-connexe, ou \textit{connexe externe}, ou encore \textit{extérieurement connexe}). Étant donné un feuilletage connectif $Z$, on notera $\vert Z\vert$ son support, $\kappa_0(Z)$ sa structure connective interne et $\kappa_1(Z)$ sa structure connective externe, de sorte que $Z=(\vert Z\vert, \kappa_0(Z), \kappa_1(Z))$. Souvent, on notera  $Z_0$ l'espace connectif intérieur $Z_0=(\vert Z\vert, \kappa_0(Z))$, et $Z_1$ l'espace connectif extérieur $Z_1=(\vert Z\vert, \kappa_1(Z))$.

\begin{df}\label{categorie FC} La \emph{catégorie des feuilletages connectifs} $\mathbf{FC}$ a pour objets les feuilletages connectifs, et pour morphismes d'un feuilletage $Z$ vers un feuilletage $Z'$ les applications $\vert Z\vert\to\vert Z'\vert$ qui sont connectives de $Z_i=(\vert Z\vert,\kappa_i(Z))$ vers $Z'_i=(\vert Z'\vert,\kappa_i(Z'))$ pour chacun des deux indices $i\in\{0,1\}$.
\end{df}

%
%
%

\begin{df} [Feuilles] Soit $Z$ un feuilletage. On appelle \emph{domaine} de $Z$, et on note $dom(Z)$, la partie présente de la structure interne $\kappa_0(Z)$. On appelle \emph{feuilles} de $Z$ les composantes connexes non vides de la structure interne $\kappa_0(Z)$. La \emph{structure interne} d'une feuille $F$ est la structure connective induite sur $F$ par $\kappa_0(Z)$. La structure externe de $F$ est la structure induite sur $F$ par $\kappa_1(Z)$. 
\end{df}

Pour tout feuilletage $Z$, on note $\mathcal{F}(Z)$ l'ensemble des feuilles de $Z$. Si $dom(Z)$ est non vide, $\mathcal{F}(Z)$ en constitue une partition.  

\begin{rmq} Il ne peut exister de composante connexe intérieurement vide que dans le cas où le domaine du feuilletage est lui-même vide, autrement dit lorsque la structure interne est la structure désintégrée. Et dans ce cas, il n'y a pas de feuilles : $\mathcal{F}(Z)=\emptyset$.
\end{rmq}

\begin{rmq}Par définition, chaque feuille est intérieurement connexe. Par contre, si le feuilletage n'est pas régulier, une feuille peut ne pas être extérieurement connexe. 
\end{rmq}


\begin{df} On dira qu'un morphisme de feuilletages $\phi:Z\to Z'$ est \emph{strict} si $\phi^\mathcal{P}$ transforme toute feuille de $Z$ en une feuille de $Z'$. La catégorie ayant pour objets les feuilletages connectifs et pour morphismes les morphismes de feuilletages stricts sera notée $\mathbf{FS}$.
\end{df}


\begin{exm} Un espace topologique $Y$ muni d'une relation d'equivalence $\rho$ définit un feuilletage connectif, en prenant $(\vert Z\vert,\kappa_1(Z))=U_T(Y)$ et $\kappa_0(Z)=\rho$, la structure connective associée à la relation d'équivalence $\rho$. 
\end{exm}





\begin{exm} Un exemple fondamental de feuilletages connectifs est issu des feuilletages au sens des variétés : une variété feuilletée est en effet munie de deux structures de variété de dimensions différentes, la structure de plus faible dimension étant construite sur une topologie plus fine que celle sur laquelle est construite la structure de variété de dimension plus grande. On prendra pour structure connective interne celle associée à la topologie la plus fine, et celle associée à l'autre pour structure connective externe. A noter que dans cet exemple, les composantes connexes pour la structure interne (les feuilles) sont nécessairement connexes également pour la structure externe.
\end{exm}

\subsection{Espace induit des feuilles d'un feuilletage}

\begin{df}[Espace induit des feuilles]\label{espace induit des feuilles} Soit $Z=(\vert Z\vert,\kappa_0(Z), \kappa_1(Z))$ un feuilletage connectif. 
L'\emph{espace induit} des feuilles de $Z$ est l'espace connectif noté $\mathcal{F}^\downarrow(Z)$, ou plus simplement $Z^\downarrow$,  de support $\vert Z^\downarrow\vert$ l'ensemble $\mathcal{F}(Z)$ des feuilles de $Z$, 
et de structure connective celle qui y est induite par l'espace connectif des parties non vides $\mathcal{P}^*(Z_1)$, 
où $Z_1=(\vert Z\vert, \kappa_1(Z))$, de sorte qu'un ensemble $\mathcal{A}$ de feuilles 
est $\kappa(Z^\downarrow)$-connexe si et seulement si $\bigcup_{F\in\mathcal{A}}{F}\in \kappa_1(Z)$.
\end{df}

\begin{rmq} Si une $\kappa_0$-composante connexe n'est pas $\kappa_1$-connexe, elle définit un point non connexe de l'espace $Z^\downarrow$. Ainsi, l'espace induit des feuilles d'un feuilletage $Z$ est-il intègre si et seulement si toute composante connexe de la structure interne de $Z$ est extérieurement connexe.
\end{rmq}

\begin{rmq} Nous appellerons également \emph{entrant} l'espace induit des feuilles d'un feuilletage, pour exprimer ce qu'on pourrait appeler la suprématie qui y est accordée à la structure externe du feuilletage sur la structure interne. 
\end{rmq}

\subsection{Espace quotient des feuilles}

Il existe une autre façon, qui pourrait d'ailleurs sembler plus naturelle, de munir l'espace des feuilles d'une structure connective. En effet, les feuilles d'un feuilletage $Z$ étant
les composantes connexes (non vides) de sa structure interne $\kappa_0(Z)$, elles sont également les classes de la relation d'équivalence partielle $\chi[\kappa_0(Z)]$, d'où très naturellement la définition suivante de l'espace \textit{quotient} des feuilles, fondée sur la notion de quotient par une relation d'équivalence partielle (voir la définition \ref{def quotient par une relation partielle} page \pageref{def quotient par une relation partielle}).

\begin{df} [Espace quotient des feuilles] 
Soit $Z$ un feuilletage. L'\emph{espace quotient des feuilles} de $Z$, noté $\mathcal{F}^\uparrow(Z)$ ou plus simplement $Z^\uparrow$, est le quotient de l'espace connectif externe $Z_1=(\vert Z\vert, \kappa_1(Z))$ par la relation d'équivalence partielle $\chi[\kappa_0(Z)]$
 \[ 
Z^\uparrow=Z_1/{\chi[\kappa_0(Z)]}.
\]
\end{df}

\begin{rmq} Nous appellerons également \emph{sortant} l'espace quotient des feuilles d'un feuilletage, premièrement en relation avec le fait que, sous des conditions en pratique souvent vérifiées (par exemple que l'espace externe soit intègre), l'espace quotient est lui-même intègre, ce qui revient à affirmer en un sens la suprématie  du point vue interne sur le point de vue externe, deuxièmement parce que cet espace est destiné à être représenté dans un espace dont la structure dépendra non seulement de la structure externe du feuilletage, mais aussi de sa structure interne, troisièmement par opposition à l'espace entrant des feuilles.
\end{rmq}

\begin{prop} Pour qu'un ensemble de feuilles soit une partie connexe de l'espace quotient, il \textit{suffit} qu'il existe inclus dans l'union de ces feuilles un connexe externe rencontrant chacune des feuilles. Cependant, cette condition n'est pas, en général, nécessaire.
\end{prop}
\noindent \textbf{Preuve}. La condition est évidemment suffisante, puisque l'image par la surjection canonique d'un tel connexe externe est égal à l'ensemble des feuilles considérées, qui est donc connexe dans l'espace quotient. 

Montrons que cette condition n'est pas nécessaire. Soit $Z$ le feuilletage défini par $\vert Z\vert=\{a,a',b,b',c,c'\}$, $\kappa_0(Z)=\{\emptyset,\{a,a'\},\{b,b'\},\{c,c'\}\}$ et $\kappa_1(Z)=\{\emptyset,\{a,b\},\{b',c'\}\}$.  
Il y a trois feuilles, notons-les $A, B, C$ (avec $A=\{a,a'\}$, etc.). 
La structure de l'espace quotient est \[\kappa(Z^\uparrow)=\{\emptyset, \{A,B\}, \{B,C\}, \{A,B,C\}\}.\]

C'est donc un espace connexe, même s'il n'existe pas de connexe dans l'espace externe rencontrant chacune de ces trois feuilles (à noter que l'espace induit des feuilles est, lui, désintégré).
\begin{flushright}$\square$\end{flushright} 
\pagebreak[3]

\begin{cor} Pour tout feuilletage, l'espace entrant (induit) des feuilles est plus fin que l'espace sortant (quotient) des feuilles.
\end{cor}

Néanmoins, même dans le cas de la structure quotient, l'espace des feuilles n'est pas nécessairement intègre (si une feuille ne contient aucun connexe externe, le singleton qu'elle forme sera non connexe).

\section{Relations fonctorielles entre feuilletages et représentations}

\subsection{Foncteurs $\mathbf{RC} \rightarrow \mathbf{FC}$}

\`{A} toute représentation connective $\rho:ob(\rho)\rightsquigarrow sp(\rho)$ on souhaite associer fonctoriellement un feuilletage $\Phi(\rho)$. Pour la structure externe de ce feuilletage, on prendra naturellement la structure de l'espace $sp(\rho)$ de la représentation. Par contre, il y a plusieurs choix possibles, \textit{a priori} légitimes,  pour la structure interne du feuilletage. 

Pour préciser ces choix, nous aurons besoin de faire appel à ce que nous appellerons des structures connectives fonctorielles :

\begin{df} Une \emph{structure connective fonctorielle} est une application $\gamma$ définie sur la classe des espaces connectifs et qui à tout espace connectif $B$ associe une structure connective $\gamma(B)$ sur $\vert B\vert$ qui soit fonctorielle au sens où il existe un endofoncteur $\Gamma$ de $\mathbf{Cnc}$ défini 
\begin{itemize}
\item sur les objets de $\mathbf{Cnc}$ par $\Gamma(B)=(\vert B\vert,\gamma(B))$, 
\item sur les flèches par $\Gamma(f)=f$.
\end{itemize}

Nous dirons qu'une structure connective fonctorielle $\gamma$ est plus fine qu'une autre, $\gamma'$, et l'on notera $\gamma\subset \gamma'$, si pour tout espace connectif $B$ on a $\gamma(B)\subset \gamma'(B)$.
\end{df}

Par exemple, notons respectivement $\kappa_D(B)$ et $\kappa_G(B)$ la structure connective désintégrée et la structure connective grossière sur $\vert B\vert$.  Alors  $\kappa_D$ et  $\kappa_G$ sont des structures connectives fonctorielles. De même, l'application notée simplement $\kappa$ qui à tout espace connectif $B$ associe sa structure connective $\kappa(B)$ est une structure connective fonctorielle, et l'on a :
\[ 
\kappa_D\subset \kappa \subset \kappa_G.
 \]

Soit $(\gamma_0,\gamma_1)$ 
un couple de structures connectives fonctorielles tel que 
$\gamma_0 \subset \gamma_1$. On va définir à partir de ce couple la structure interne du feuilletage associé à une représentation $\rho$ comme la structure connective engendrée par les parties des $\rho(a)$ qui sont connexes pour l'une ou l'autre des structures $\gamma_i$ selon que $a$ est un point connexe ou non de l'objet de la représentation $\rho$. Plus précisément, à toute représentation connective 
$\rho : A \rightsquigarrow B$,
on associe le feuilletage $Z=\Phi_{(\gamma_0,\gamma_1)}(\rho)= \Phi(\rho)$ de support 
$\vert Z\vert=\vert B\vert$,
de structure externe $\kappa_1(Z)=\kappa(B)$ et de structure interne
\[ 
\kappa_0(Z)=[\bigcup_{i\in\{0,1\}}\bigcup_{a\in A_i}{( \gamma_i(B)\cap \mathcal{P}_{\rho(a)})}]_0,
 \]
où $A_0$ désigne la partie absente de $A$ et $A_1$ sa partie présente.

\begin{prop}\label{prop Phi foncteur} Soit $(\alpha,\beta):\rho\to\rho'$  un morphisme de représentations connectives. Alors l'application $\beta:\vert sp(\rho)\vert \to \vert sp(\rho')\vert$ est un morphisme de feuilletages $\beta:\Phi(\rho)\to\Phi(\rho')$.
\end{prop}
\noindent \textbf{Preuve}. Posons $A=ob(\rho)$, $B=sp(\rho)$, $A'=ob(\rho')$ et $B'=sp(\rho')$.

Si $K$ est une partie extérieurement connexe du feuilletage $Z=\Phi(\rho)$, alors $\beta^\mathcal{P}(K)$ est une partie extérieurement connexe de $Z'=\Phi(\rho')$ puisque les structures extérieures des feuilletages coïncident avec les structures des espaces de représentation que respecte $\beta$. 

Soit maintenant $K$ un connexe de base\footnote{C'est-à-dire un connexe appartenant à la famille qui, dans la définition qui en est donnée ci-dessus, engendre la structure interne de $Z$.} pour la structure interne de $Z$. On veut montrer que $\beta^\mathcal{P}(K)$ est intérieurement connexe dans $Z'$.

Si $K\in\gamma_0(B)\cap \mathcal{P}_{\rho(a)}$  avec $a\in A_0$, alors $\beta(K)\in\gamma_0(B')$ puisque $\gamma_0$ est fonctoriel. Et $K\subset \rho(a) \Longrightarrow \beta^\mathcal{P}(K)\subset \beta^\mathcal{P}(\rho(a))\subset \rho'(\alpha(a))$. Si $a'=\alpha(a)\in A'_0$, on a alors $\beta^\mathcal{P}(K)\in \gamma_0(B')\cap \mathcal{P}_{\rho(a')}\subset\kappa_0(Z')$, tandis que si $a'\in A'_1$, on a $\beta^\mathcal{P}(K)\in \gamma_1(B')\cap \mathcal{P}_{\rho(a')}\subset\kappa_0(Z')$,
puisque $\gamma_0\subset\gamma_1$.

Si $K\in\gamma_1(B)\cap \mathcal{P}_{\rho(a)}$ avec $a\in A_1$, on a nécessairement $a'=\alpha(a)\in A'_1$, et comme précédemment le fait que $\beta^\mathcal{P}(\rho(a))\subset \rho'(a')$ permet de conclure que $\beta^\mathcal{P}(K)\in \gamma_1(B')\cap \mathcal{P}_{\rho(a')}\subset\kappa_0(Z')$,
puisque $\gamma_1$ est fonctoriel.
\begin{flushright}$\square$\end{flushright} 
\pagebreak[3]

On en déduit que l'application $\Phi=\Phi_{(\gamma_0,\gamma_1)}$ qui à toute représentation $\rho$ associe le feuilletage $\Phi(\rho)$ et à tout morphisme de représentations $(\alpha,\beta)$ associe $\beta$ est un  foncteur de la catégorie des représentations dans celle des feuilletages. Dans le cas où $\gamma_0=\gamma_1=\gamma$, on le notera simplement $\Phi_\gamma$. Lorsque $\gamma=\kappa_G$ (resp. $\gamma=\kappa_D$), on notera simplement $\Phi_G$ (resp. $\Phi_D$) le foncteur $\Phi_\gamma$\label{notation PHIG}.

\begin{exm} Soit $\rho$ la représentation claire et distincte, dans la droite $\mathbf{R}$ usuelle, de l'objet $X=(\{a,b,c\},\{\emptyset,\{a\},\{b\},\{a,b,c\} \})$ définie par $\rho(a)=]-2,-1[$, $\rho(b)=]1,2[$ et $\rho(c)=[-3,-2]\cup[-1,1]\cup[2,3]$.

Les feuilletages $\Phi_{(\kappa_D, \kappa)}(\rho)$ et $\Phi_{(\kappa_D, \kappa_G)}(\rho)$ possèdent chacun deux feuilles, et celles-ci sont extérieurement connexes. $\Phi_{\kappa}(\rho)$ et $\Phi_{(\kappa, \kappa_G)}(\rho)$ ont chacun cinq feuilles  extérieurement connexes. 
Enfin, $\Phi_{\kappa_G}(\rho)=\Phi_G(\rho)$ a trois feuilles, deux extérieurement connexes, une extérieurement non connexe, toutes trois étant de structure interne grossière.
\end{exm}


\begin{prop}\label{Phik regulier} Pour toute représentation connective $\rho$, le feuilletage $\Phi_\kappa (\rho)$ est régulier.
\end{prop}
\noindent \textbf{Preuve}. La structure interne de $\Phi_\kappa (\rho)$ étant (par définition de cette structure interne et de la structure fonctorielle $\kappa$) engendrée par des parties extérieurement connexes, elle est nécessairement plus fine que la structure externe.
\begin{flushright}$\square$\end{flushright} 
\pagebreak[3]

La proposition suivante découle immédiatement des définitions.

\begin{prop} Soient $(\gamma_0,\gamma_1)$ un couple de structures connectives fonctorielles, tel que $\gamma_0 \subset\gamma_1$, et soit $\Phi=\Phi_{(\gamma_0,\gamma_1)}$ le foncteur $\mathbf{RC}\to\mathbf{FC}$ associé. 

Si $\rho$ est une représentation distincte alors, pour $Z=\Phi(\rho)$, on a
\[ 
\kappa_0(Z)=\bigcup_{i\in\{0,1\}}\bigcup_{a\in A_i}{( \gamma_i(B)\cap \mathcal{P}_{\rho(a)})},
 \]
de sorte que si $\gamma_1\supset\kappa$ et que $a$ est un point connecté de $ob(\rho)$, alors $\rho(a)$ est une composante connexe de $(Z,\kappa_0(Z))$, c'est-à-dire une feuille de $Z$.


Si $ob(\rho)$ est intègre, alors la première structure connective fonctorielle $\gamma_0$ de $\Phi_{(\gamma_0,\gamma_1)}$ n'a trivialement aucune incidence : si $\gamma'_0$ est une autre structure connective fonctorielle telle que $\gamma'_0 \subset\gamma_1$, alors $\Phi_{(\gamma_0,\gamma_1)}(\rho)=\Phi_{(\gamma'_0,\gamma_1)}(\rho)$.


\end{prop}

\begin{cor}\label{cor feuilles representation} Si $\rho$ est une représentation distincte d'un objet intègre, les feuilles de $\Phi_\kappa (\rho)$ sont les parties de $sp(\rho)$ de la forme $\rho(a)$ :
\[ 
\mathcal{F}(\Phi_\kappa(\rho))=\{\rho(a), a\in ob(\rho)\}.
 \]
\end{cor}

\subsection{Foncteurs $\mathbf{FC} \rightarrow \mathbf{RCD}$}

\`{A} tout feuilletage on veut associer une représentation de son espace des feuilles.  Puisque nous avons vu qu'il y avait deux notions distinctes possibles pour la structure de l'espace des feuilles, à savoir l'espace induit et l'espace quotient, nous sommes conduits à définir deux foncteurs différents : le foncteur \og représentation induite\fg\, $\mathcal{R}^\downarrow$ et le foncteur \og représentation quotient\fg\, $\mathcal{R}^\uparrow$.

\subsubsection{Représentation induite $\mathcal{R}^\downarrow$ d'un feuilletage}

\paragraph{Définition de $\mathcal{R}^\downarrow$  sur les objets}

\`{A} tout feuilletage $Z$,  $\mathcal{R}^\downarrow$  associe la représentation $\mathcal{R}^\downarrow(Z):Z^\downarrow \rightsquigarrow Z_1$ définie pour toute feuille $F\in \mathcal{F}(Z)=\vert Z^\downarrow \vert$ par
\[ 
\mathcal{R}^\downarrow(Z)(F)=F\subset\vert Z_1\vert.
 \]
Ceci est bien une représentation connective puisque, par définition, un ensemble de feuilles est connexe dans $\mathcal{P}_{Z_1}$ si et seulement si son union est connexe dans $Z_1$, et que cette dernière propriété caractérise la structure connective de l'espace induit $Z^\downarrow$.  On vérifie en outre immédiatement que la représentation $\mathcal{R}^\downarrow(Z)$ est claire (si un ensemble de feuilles est non connexe dans $Z^\downarrow$, alors leur union est également non connexe dans l'espace externe du feuilletage, donc dans l'espace de la représentation), et distincte (deux point différents, c'est-à-dire deux feuilles différentes, sont représentées par deux composantes connexes internes nécessairement disjointes).

\begin{prop}\label{objet integre si feuilletage regulier} Si le feuilletage $\mathcal{Z}$ est régulier, l'objet de la représentation $\mathcal{R}^\downarrow \mathcal{Z}$ est intègre.
\end{prop}

\noindent \textbf{Preuve}. Toute feuille étant extérieurement connexe, elle constitue un singleton connexe de $ob(\mathcal{R}^\downarrow \mathcal{Z})$.
\begin{flushright}$\square$\end{flushright} 
\pagebreak[3]

\paragraph{Définition de $\mathcal{R}^\downarrow$  sur les flèches}

$\mathcal{R}^\downarrow$ est défini sur les flèches de $\mathbf{FC}$ en associant à tout morphisme de feuilletage $\phi:Z\to Z'$ le morphisme de représentations $(\phi_0,\phi_1)$, où $\phi_0: Z^\downarrow \to Z'^\downarrow$  est défini pour toute feuille $F\in \mathcal{F}(Z)$ par : $\phi_0(F)$ est celle des composantes connexes de l'espace interne $(\vert Z'\vert,\kappa_0(Z')$ qui contient le $\kappa_0(Z')$-connexe $\phi^\mathcal{P}(F)$, et où $\phi_1:Z_1\to Z'_1$ est le morphisme connectif qui en tant qu'application ensembliste coïncide avec $\phi$.

\subsubsection{Représentation quotient $R^\uparrow$}

\paragraph{Définition de $R^\uparrow$ sur les objets} 
\`{A} tout feuilletage $Z$ est associé la représentation 
$R^\uparrow(Z):\mathcal{F}^\uparrow(Z)\rightsquigarrow (\vert Z\vert,[\kappa_0(Z)\cup \kappa_1(Z)])$
de l'espace quotient des feuilles $\mathcal{F}^\uparrow(Z)$ dans l'espace $(\vert Z\vert,[\kappa_0(Z)\cup \kappa_1(Z)])$ qui à toute feuille $F$ (qu'elle soit ou non connexe dans $\mathcal{F}^\uparrow(Z)$) associe elle-même en tant que partie (nécessairement connexe) de $(\vert Z\vert,[\kappa_0(Z)\cup \kappa_1(Z)])$. 
Ceci définit bien une représentation car si un ensemble  de feuilles est un connexe \og de base\fg\, dans $\mathcal{F}^\uparrow(Z)$, c'est-à-dire qu'il s'agit d'un ensemble de feuilles dont l'union contient un connexe $K$ de $Z_1$ rencontrant toutes ces feuilles, alors cette union est connexe dans $(\vert Z\vert,[\kappa_0(Z)\cup \kappa_1(Z)])$ du fait de la  $\kappa_0(Z)$-connexité de chaque feuille et de la $\kappa_1(Z)$-connexité de $K$.

\begin{exc} Préciser ce que doit être l'action du foncteur $R^\uparrow$ sur les morphismes de feuilletages.
\end{exc}

\begin{rmq}
La représentation quotient associée à un feuilletage n'est pas nécessairement distincte (une feuille non connexe peut \og devenir connexe\fg\, dans l'espace de représentation). 
\end{rmq}

%
%
%

\subsection{Une adjonction}

La composition des divers foncteurs définis plus haut entre catégories de feuilletages et catégories de représentations conduisent à des endofoncteurs présentant diverses propriétés intéressantes. Par exemple, notons
\[ 
\rho^{\downarrow}_{G} = \mathcal{R}^\downarrow(\Phi_{G}(\rho) )
 \]
la représentation associée à une représentation $\rho$ par l'endofoncteur $\mathcal{R}^\downarrow \circ\Phi_{G}$. On a alors la proposition suivante.

\begin{prop} Si $\rho$ est une représentation claire et distincte, alors le couple d'applications $(\alpha,\beta)$ défini par 
$\alpha(a)=\rho(a)\in ob(\rho^{\downarrow}_{G})$ et $\beta=Id_{sp(\rho)}$ constitue un isomorphisme entre les représentations  $\rho$ et
$\rho^{\downarrow}_{G}$.
\end{prop}

\noindent \textbf{Preuve}. Soit $\rho$ une telle représentation. Le feuilletage associé $\Phi_G(\rho)$ a pour structure interne celle engendrée par les parties des $\rho(a)$ lorsque $a$ décrit $ob(\rho)$. La représentation étant distincte, les $\rho(a)$ sont deux à deux disjoints, de sorte que les composantes connexes de la structure interne du feuilletage sont précisément les $\rho(a)$, et l'ensemble des feuilles s'identifie à l'objet de $\rho$. Si un ensemble $K\subset ob(\rho)$ est connexe, l'ensemble correspondant de feuilles, $\{\rho(a), a\in K\}$, est connexe dans l'espace induit des feuilles $\mathcal{F}^{\downarrow}(\Phi_G(\rho))=ob(\rho^{\downarrow}_{G})$ puisque $\rho$ étant une représentation, on a $\rho^\mathcal{P}(K)$ connexe dans $\mathcal{P}^*(sp(\rho))$, autrement dit $^\mu{\rho^\mathcal{P}(K)}=\bigcup_{a\in K}\rho(a)$ est connexe pour la structure externe du feuilletage $\Phi_G(\rho)$. Si au contraire un ensemble $A\subset ob(\rho)$ est non connexe, la représentation $\rho$ étant claire, on aura $\bigcup_{a\in K}\rho(a)$ non connexe dans $sp(\rho)$, autrement dit non connexe pour la structure externe du feuilletage $\Phi_G(\rho)$, de sorte que l'ensemble des feuilles $\{\rho(a), a\in A\}$ sera non connexe dans l'espace $\mathcal{F}^{\downarrow}(\Phi_G(\rho))=ob(\rho^{\downarrow}_{G})$. Ainsi, $\rho$ et $\rho^{\downarrow}_{G}$ ont-ils des objets isomorphes. Ils ont également même espace, et le couple indiqué constitue alors trivialement un isomorphisme entre ces deux représentations.
\begin{flushright}$\square$\end{flushright} 
\pagebreak[3]

Beaucoup d'autres propriétés restent à explorer. On se contentera dans cette section de prouver l'existence d'une adjonction entre les foncteurs $\mathcal{R}^\downarrow$ et $\Phi_\kappa$ lorsqu'ils sont restreints à certaines catégories de feuilletages et de représentations.
Pour cela, nous ferons appel aux trois lemmes suivants (lemme \ref{lemme beta est un morphisme de feuilletages } à lemme \ref{lemme existence de alpha donnant representation}).

\begin{lm}\label{lemme beta est un morphisme de feuilletages } Soit $Z$ un feuilletage \emph{régulier}, $\rho$ une représentation quelconque, et $(\alpha,\beta) :\mathcal{R}^\downarrow Z \rightarrow \rho $ un morphisme de représentations. Alors $\beta$ est un morphisme de feuilletages $Z\rightarrow \Phi_\kappa (\rho)$.
\end{lm}
\noindent \textbf{Preuve}. Par définition d'un morphisme de représentations, $\beta$ est un morphisme connectif $sp(\mathcal{R}^\downarrow Z) \rightarrow sp(\rho)$, autrement dit un morphisme connectif $Z_1\rightarrow(\Phi_\kappa(\rho))_1$.

D'autre part, en appliquant le foncteur $\Phi_\kappa$ au morphisme $(\alpha,\beta)$ (proposition \ref{prop Phi foncteur}), on en déduit que $\beta$ est un morphisme de feuilletages $\Phi_\kappa(\mathcal{R}^\downarrow Z) \rightarrow \Phi_\kappa(\rho)$, donc en particulier un morphisme pour les structures internes

$(\Phi_\kappa(\mathcal{R}^\downarrow Z))_0 \rightarrow (\Phi_\kappa(\rho))_0$.

Mais $Z$ étant régulier, $\kappa_0(Z)\subset \kappa_0(\Phi_\kappa(\mathcal{R}^\downarrow Z))$. En effet, tout connexe intérieur est trivialement inclus dans une composante connexe intérieur et, par la régularité de $Z$, est aussi un connexe extérieur, de sorte que, par définition de la structure $\kappa_0(\Phi_\kappa(\mathcal{R}^\downarrow Z))$, se trouve bien appartenir à celle-ci.

Finalement, on a à la fois $\beta:Z_1\rightarrow(\Phi_\kappa(\rho))_1$ et $\beta:Z_0\rightarrow(\Phi_\kappa(\rho))_0$, autrement dit $\beta$ est bien un morphisme $Z\rightarrow \Phi_\kappa(\rho)$.
\begin{flushright}$\square$\end{flushright} 
\pagebreak[3]

\begin{lm}\label{lemme beta determine alpha} Soient $\mathcal{Z}$ un feuilletage connectif, $\rho$ une représentation connective \emph{distincte} et $(\alpha,\beta):\mathcal{R}^\downarrow(Z) \rightarrow \rho$ un morphisme de représentations. Alors la connaissance de $\beta$ détermine celle de $\alpha$. Autrement dit, si $(\alpha',\beta):\mathcal{R}^\downarrow(Z) \rightarrow \rho$ est également un morphisme de représentations, on a nécessai\-rement $\alpha=\alpha'$.
\end{lm}
\noindent \textbf{Preuve}. 
Par définition, $\alpha$ est un morphisme connectif de 
$ob(\mathcal{R}^\downarrow Z)=\mathcal{F}^\downarrow(Z)=Z^\downarrow$ 
dans $ob(\rho)$. 
Soit $F\in ob(\mathcal{R}^\downarrow(Z))$, 
autrement dit une composante connexe de $Z_0=(\vert Z\vert, \kappa_0(Z))$. 
Par définition d'un morphisme de représentations, on a l'inclusion
\[ 
\beta^\mathcal{P}(\mathcal{R}^\downarrow Z (F))\subset \rho(\alpha(F)).
 \]
Or, $\mathcal{R}^\downarrow Z (F)=F \subset \vert Z \vert$, d'où 
$\beta^\mathcal{P}(F)\subset\rho(\alpha(F))$. 
\pagebreak[3]
La représentation $\rho$ étant distincte, il n'y a au plus qu'un point $a$ de $ob(\rho)$ pouvant vérifier $\beta^\mathcal{P}(F)\subset\rho(a)$, d'où l'unicité annoncée. \begin{flushright}$\square$\end{flushright} 
\pagebreak[3]

\begin{lm}\label{lemme existence de alpha donnant representation}
Soit $Z$ un feuilletage, $\rho$ une représentation claire et distincte, d'objet $ob(\rho)$ intègre, et soit $\beta:Z\rightarrow \Phi_\kappa(\rho)$ un morphisme de feuilletages. Alors il existe un et un seul morphisme connectif $\alpha:\mathcal{F}^\downarrow Z \rightarrow ob(\rho)$ tel que $(\alpha,\beta)$ soit un morphisme de représentations $\mathcal{R}^\downarrow Z\rightarrow\rho$.
\end{lm}
\noindent \textbf{Preuve}. S'il existe, le morphisme $\alpha$ est unique d'après le lemme \ref{lemme beta determine alpha}. Précisons l'application ensembliste $\alpha:\mathcal{F}Z\rightarrow\vert ob(\rho)\vert$  dont, nécessairement, il s'agit. Pour $F\in \mathcal{F}Z$, on a $\beta^\mathcal{P}(F)\in \kappa_0(\Phi_\kappa(\rho))$, puisque $\beta$ préserve aussi les morphismes internes. 

Notons $\overline{\beta^\mathcal{P}(F)}$ la composante $\kappa_0(\Phi_\kappa(\rho)$-connexe contenant $\beta^\mathcal{P}(F)$. Alors $\overline{\beta^\mathcal{P}(F)}\in\mathcal{F}(\Phi_\kappa(\rho))$. D'après le corolaire \ref{cor feuilles representation}, il existe alors un élément unique $a_F\in ob(\rho)$ tel que $\overline{\beta^\mathcal{P}(F)}=\rho(a_F)$. L'application $\alpha$ est donc définie par $\alpha(F)=a_F$. Autrement dit, 

\[ 
\alpha(F)=a 
\Leftrightarrow \beta^\mathcal{P}(F) \subset \rho(A) 
\Leftrightarrow 
\beta^\mathcal{P}(F)\subset \overline{\beta^\mathcal{P}(F)}= \rho(A).
 \]
Il s'agit de prouver que l'application $\alpha$ ainsi définie est un morphisme connectif $\mathcal{F}^\downarrow Z \rightarrow ob(\rho)$, et que le couple $(\alpha,\beta)$ est bien un morphisme de repré\-senta\-tions. 

Soit donc $\mathcal{L}$ un ensemble $\kappa(Z^\downarrow)$-connexe de feuilles. Par définition de $Z^\downarrow$, on a $\bigcup_{F\in\mathcal{L}} F \in \kappa_1(Z)$, donc l'ensemble $W=\bigcup_{F\in\mathcal{L}} \beta^\mathcal{P}(F)$ vérifie $W \in \kappa_1(\Phi_\kappa(\rho))$.

Posons  $\mathcal{A}=\alpha^\mathcal{P}(\mathcal{L})=\{a\in ob(\rho), \exists F\in\mathcal{L}, \rho(a) \supset \beta^\mathcal{P}(F)\}$. 
On veut montrer que $\mathcal{A}$ est une partie connexe de $ob(\rho)$. 
Or, $\rho$ étant claire, il suffit pour cela de prouver que $^\mu \rho (\mathcal{A})= 
\bigcup_{F\in\mathcal{L}}\overline{\beta^\mathcal{P}(F)}$ est connexe dans $sp(\rho)$. 

Par définition, les $\overline{\beta^\mathcal{P}(F)}$ sont $\kappa_0(\Phi_\kappa(\rho))$-connexes. Mais, le feuilletage $\Phi_\kappa(\rho)$ étant régulier (proposition \ref{Phik regulier}), les $\overline{\beta^\mathcal{P}(F)}$ sont également $\kappa_1(\Phi_\kappa(\rho))$-connexes. Il en découle que

\[ 
\bigcup_{F\in\mathcal{L}}\overline{\beta^\mathcal{P}(F)}
=
\bigcup_{F\in\mathcal{L}}(\overline{\beta^\mathcal{P}(F)}\cup W)
 \]
\noindent est l'union de $\kappa_1(\Phi_\kappa(\rho))$-connexes d'intersection non vide . Ainsi, $\bigcup_{F\in\mathcal{L}}\overline{\beta^\mathcal{P}(F)}$ est $\kappa_1(\Phi_\kappa(\rho))$-connexe, autrement dit $\kappa_1(sp(\rho))$-connexe, de sorte que $\mathcal{A}$ est connexe dans $ob(\rho)$. 

Reste à vérifier que $\beta^\mathcal{P} \circ \mathcal{R}^\downarrow Z \subset \rho \circ \alpha$, mais c'est là une conséquence immédiate de la construction même de $\alpha$.
\begin{flushright}$\square$\end{flushright} 
\pagebreak[3]

Soit $\mathbf{FR}$ la sous-catégorie pleine de  $\mathbf{FC}$  constituée des feuilletages connectifs réguliers, et soit $\mathbf{RIO}$ la sous-catégorie pleine de $\mathbf{RCD}$ constituée des représentations claires et distinctes dont l'objet est intègre. Reprenons les notations  $\mathcal{R}^\downarrow$ et $\Phi_\kappa$ employées précédemment, mais pour désigner cette fois les restrictions de ces foncteurs à $\mathbf{FR}$ et à $\mathbf{RIO}$. 

D'après la proposition \ref{objet integre si feuilletage regulier}, on obtient bien de cette manière un foncteur $\mathcal{R}^\downarrow:\mathbf{FR}\rightarrow \mathbf{RIO}$. Et d'après la proposition \ref{Phik regulier}, on obtient de même un foncteur $\Phi_\kappa:\mathbf{RIO}\rightarrow \mathbf{FR}$.

Soit maintenant $Z$ un feuilletage régulier, et $\rho$ une représentation claire et distincte d'un objet intègre. \`{A} tout morphisme de représentation $(\alpha,\beta):\mathcal{R}^\downarrow Z \rightarrow \rho$, on associe, d'après le lemme \ref{lemme beta est un morphisme de feuilletages }, le morphisme de feuilletages $\beta: Z\rightarrow \Phi_\kappa (\rho)$. Réciproquement, à tout morphisme de feuilletage $\beta: Z\rightarrow \Phi_\kappa (\rho)$, on associe d'après le lemme \ref{lemme existence de alpha donnant representation}, un unique morphisme de représentations $(\alpha,\beta):\mathcal{R}^\downarrow Z \rightarrow \rho$. On a ainsi construit des applications réciproques, donc bijectives, entre $Hom_{\mathbf{RIO}}(\mathcal{R}^\downarrow Z, \rho)$ et $Hom_{\mathbf{FR}}(Z, \Phi_\kappa(\rho))$, et il est clair que ces bijections sont naturelles par rapport à $Z$ et $\rho$. On peut ainsi énoncer :

Pour tout feuilletage régulier

\begin{thm}\label{Adjonction entre feuilletages et representations} Le  foncteur $\mathcal{R}^\downarrow:\mathbf{FR}\rightarrow \mathbf{RIO}$ est adjoint à gauche du foncteur $\Phi_\kappa:\mathbf{RIO}\rightarrow \mathbf{FR}$ :
\[ \mathcal{R}^\downarrow\dashv \Phi_\kappa \] 
\end{thm}


\chapter{Dynamiques catégoriques ensemblistes} \label{chapitre dynamiques categoriques ens}

\begin{flushright}

\textit{Il crée ainsi divers avenirs, divers temps}\\

\textit{qui prolifèrent aussi et bifurquent.}

\mbox{}

Jorge Luis Borgès\\

\end{flushright}

 \section{Introduction}
 
La conception classique des systèmes dynamiques est celle de systèmes déterministes dont la définition diffère selon qu'il est question de systèmes continus ou discrets, découpage qui d'ailleurs se prolonge dans la subdivision entre systèmes réversibles et systèmes irréversibles. Ces quatre types de systèmes dynamiques correspondent de fait à quatre modélisations distinctes du temps, fondées respectivement sur $\mathbf{N}$ et $\mathbf{Z}$ pour les systèmes discrets, sur $\mathbf{R}_+$ et $\mathbf{R}$ pour les systèmes continus. En outre, ces modélisations ne sont généralement pas suffisamment explicites pour que l'on y distingue parfaitement les notions d'instants et d'écoulements temporels (les durées), en témoigne en un sens la référence permanente au mythique \og instant $0$\fg, tandis que le rôle exact joué par le fait que les monoïdes cités sont ordonnés reste souvent difficile à définir. 

Souhaitant appliquer les notions considérées dans les chapitres précédents aux dynamiques qui semblent y apparaître assez spontanément --- les exemples \og naturels\fg\, de représentations connectives d'espaces connectifs infinis et de feuilletages connectifs comportant une infinité de feuilles viennent en effet souvent des systèmes dynamiques, notamment de la mécanique (fibration de Hopf, tores de Liouville-Arnold, ...) --- il était logique de chercher à reprendre à la base la question des modélisations temporelles. La base dont il s'agit, qui n'est sans doute pas aussi profonde qu'on pourrait le souhaiter --- nous n'avons en particulier pas su prendre d'emblée en compte le point de vue relativiste sur l'espace-temps, qui supposerait en quelque sorte de ne pas présupposer une source (le temps) et un but (l'espace) aux foncteurs dynamiques, mais de les faire émerger de tels foncteurs --- consiste en ce fait élémentaire : un système dynamique transforme l'addition des durées en composition des transitions d'états.  Autrement dit, l'essentiel, dans un système dynamique et pour toute notion de flot, est que les transformations des états du système soient réglées sur la composition des écoulements temporels. Ainsi, dans la définition des systèmes dynamiques, la notion d'instant n'intervient pas, seule celle de durée (ou d'écoulement temporel) est un jeu, et plus précisément la composition de ces écoulements.
C'est pourquoi les équation différentielles à coefficients variables n'entrent pas (ou pas directement) dans la théorie classique des systèmes dynamiques\footnote{La théorie des systèmes dynamiques non-autonomes (\textit{Nonautonomous Dynamical
Systems}) vise précisément à élargir le cadre classique en ce sens.}.

Cette relégation au second plan des instants y place simultanément la relation qui habituellement les ordonne. Même si le modèle proposé plus loin n'est pas, \textit{a priori},  relativiste, il est difficile d'ignorer que, depuis 1905, les travaux d'Albert Einstein ont remis définitivement en question l'idée d'un temps totalement ordonné. Du reste, l'idée même de relation d'ordre est exclue des conceptions les plus antiques de temps cyclique, conceptions qui font régulièrement retour comme en témoigne la pensée de Nietzsche. Plus près de nous, de nombreux travaux en informatique font appel à des temporalités arborescentes, dans lesquelles les instants sont partiellement ordonnés.  C'est notamment le cas, depuis le début des années 1980, de la logique dite CTL (\textit{Computation tree logic}) \cite{Ben-Ari:1983, Emerson_Halpern:1985}. 

Quelles sont, dès lors, les pistes qui s'offrent à nous pour une refondation des modèles temporels ? Dans des travaux récents, Claudio Mazzola et Marco Giunti \cite{Giunti_Mazzola:2010} discutent divers types d'irréversibilités dans les systèmes déterministes fondés sur une temporalité définie par un monoïde. 
Dans sa thèse, \cite{Mazzola:2010}, C. Mazzola soutient que ce cadre serait le plus général possible pour les systèmes dynamiques, et trouve un argument catégorique pour soutenir ce point de vue. 
Il apparaît pourtant que le typage des écoulements temporels, donc des restrictions posées à leur composabilité, élargit considérablement la portée de ce type de modèle. Du reste, dès 1965, Mme Andrée Ehresmann-Bastiani \cite{MmeEhresmann:1965} a considéré des \textit{systèmes guidables} (déterministes) où la temporalité était donnée par une catégorie topologique\footnote{i.e. munie d'une structure topologique} agissant sur un ensemble d'états, ces systèmes étant d'ailleurs eux-même une généralisation de certains \textit{polysystèmes dynamiques} \cite{Bushaw:1963}. 
Dans ce qui suit, nous reprenons cette idée et en explorons diverses conséquences qui nous semblent nouvelle, notamment avec les notions de solutions existentielles et essentielles d'une dynamique, d'interprétation d'une dynamique par une autre, également avec la composition de certains foncteurs reliant la catégorie des petites catégories et la catégorie que nous définissons comme catégorie des dynamiques. \`{A} noter qu'une grande partie de l'intérêt selon nous de ces considérations vient de ce qu'elles s'appliquent à des dynamiques non-déterministes.

\section{Ecoulements catégoriques}

\begin{df} Un \emph{système catégorique d'écoulements} $\mathbf{E}$ est une petite catégorie. Les flèches de $\mathbf{E}$ sont appelés les $\mathbf{E}$-\emph{écoulements}. Les objets de $\mathbf{E}$ sont les \emph{modes}\footnote{On pourrait aussi les appeler les \emph{types} d'écoulements, à condition de ne pas confondre cette notion avec celle de \emph{type d'une dynamique}, un tel type désignant une catégorie $\mathbf{E}$ (voir la définition \ref{definition dynamique de type E} page \pageref{definition dynamique de type E}).} d'écoulements. 
\end{df}


\begin{exm}\label{foncteur MC}  Comme il est bien connu, la catégorie des monoïdes (avec élément neutre) peut être plongée dans celle des petites catégories.

Bien que cette question puisse paraître triviale, plusieurs points méritent d'être précisés à ce sujet. Tout d'abord, il y a une opposition entre les conventions usuelles concernant la composition des flèches dans une catégorie, où c'est la convention employée pour la composition des applications qui domine, et la convention issue de la concaténation, convention qui convient davantage aux monoïdes puisque c'est celle employée pour les monoïdes libres. Par ailleurs, il y a deux choix possibles pour l'unique objet de la petite catégorie que l'on souhaite associer à un monoïde : soit on prend un objet abstrait arbitraire $\bullet$, et dans ce cas les flèches sont simplement définies par la manière dont elles se composent, soit on prend pour objet l'ensemble des éléments du monoïde, et dans ce cas les flèches peuvent s'interpréter comme des actions sur le monoïde, plus précisément comme des actions à droite ou à gauche selon que la loi de composition du monoïde s'interprète dans l'ordre opposé ou dans le même ordre que la composition des morphismes correspondants.

Par exemple, on définit un foncteur \textit{covariant}, appelons-le $MC$, de la catégorie des monoïdes dans celle des petites catégories,  en posant :

\begin{itemize}
\item pour tout monoïde $(M,*)$, $MC(M,*)$ est la catégorie ayant pour unique objet l'ensemble $M$, pour flèches les applications de la forme $*f:M\rightarrow M$, définies pour tout $m\in M$ par $*f (m) = m*f$, avec $f\in M$, de sorte que, notant simplement $f$ l'action $*f$, la loi de composition $\circ$ des morphismes de cette catégorie s'écrit $f\circ g=g*f$,
\item pour tout morphisme de monoïdes $\mu:(M,*)\to (N,.)$, $MC(\mu)$ est le foncteur covariant\footnote{$MC(\mu)$ est bien un foncteur covariant puisque
$MC(\mu)(f\circ g)=\mu(g*f)=\mu(g)* \mu (f) = MC(\mu)(f) \circ MC(\mu)(g)$. } de la catégorie $MC(M,*)$ dans $MC(N,.)$ défini sur les (uniques) objets par  $MC(\mu)(M)=N$ et sur les flèches $e\in M$ par $MC(\mu)(e)=\mu(e)\in N$.
\end{itemize}

Par ce plongement, les systèmes présentant un unique mode d'écoulements s'identifient aux monoïdes. Ainsi, munis de l'addition, $\mathbf{N}$, $\mathbf{Z}$, $\mathbf{R}_+$, $\mathbf{R}$, constituent autant de systèmes commutatifs d'écoulements; le monoïde libre à deux lettres est un exemple de système non commutatif d'écoulements; les groupes cycliques forment des systèmes d'écoulements non ordonnés, etc.

\end{exm}

\begin{notations}\label{notations fg} On notera souvent $gf$ la composée $f\circ g$ de deux écoulements composables. Avec cette notation, le plongement de la catégorie des monoïdes dans la catégorie des petites catégories est définie par la formule $gf=g*f$.
\end{notations}

\begin{exm}\label{ecoulements par preordre}  Tout préodre définit une petite catégorie. En particulier, tout ensemble partiellement ordonné définit un système catégorique d'écoulements temporels. Par exemple, la catégorie définie par $(\mathbf{R},\leqslant)$ a pour écoulements les couples  de réels $(r,s)$ avec $r\leqslant s$. On notera $(r\leqslant s)$ \label{notations r inf s} un tel écoulement. On a alors
\[  
(s\leqslant t)\circ (r\leqslant s)=(r\leqslant s)(s\leqslant t)=(r\leqslant  t)
\]
\end{exm}

\begin{rmq} Lorsqu'un monoïde est ordonné, ce qui est le cas des quatre structures temporelles classiques, il donne lieu à deux catégories d'écou\-lements généralement fort différentes l'une de l'autre selon que l'on considère l'aspect monoïdal ou l'aspect ordonné. En effet, les élé\-ments d'un monoïde définissent les \textit{écoulements} d'un certain système catégorique, tandis que les élé\-ments d'un ensemble ordonné définissent les \textit{modes} d'écoulements d'un autre système catégorique. Par exemple, pour $\mathbf{E}=(\mathbf{R_+},+)$ on obtient des écoulements toujours composables, puisqu'il n'y a qu'un seul mode, tandis que pour $\mathbf{E}=(\mathbf{R_+},\leqslant)$, les écoulements sont typés, les types, c'est-à-dire les modes d'écoulements, pouvant  être identifiés à des instants\footnote{Voir plus loin la definition \ref{def instants essentiels} des \emph{instants essentiels}.}. Les aspects \textit{monoïde} et \textit{ordonné} d'un monoïde ordonné sont en quelque sorte orthogonaux l'un à l'autre. Philosophiquement, toute réflexion sur le temps devrait certainement tenir compte de cette distinction fondamentale, qui prolonge en un sens la différence entre les écoulements et les instants. 
\end{rmq}

\begin{exm} \label{exemple ecoulements RN} On constitue un système d'écoulements à deux modes en reliant les deux monoïdes $(\mathbf{R}_+,+)$ et $(\mathbf{N},+)$  par des flèches $f_r:\mathbf{R}_+ \rightarrow \mathbf{N}$ de la forme $ f_r(x)= E(x+r)$, où $E$ désigne la partie entière, la composition de trois flèches $s:\mathbf{R}_+\rightarrow\mathbf{R}_+$, $f_r:\mathbf{R}_+ \rightarrow \mathbf{N}$ et $n:\mathbf{N}\rightarrow\mathbf{N}$ étant définie par
\[ 
n\circ f_r \circ s = f_{s+r+n}: x \mapsto E(x+s+r)+n
\]
\end{exm}

\begin{exm}\label{exemple ecoulements MNMN} Plus généralement, étant donné un monoïde $(M,*)$ et une suite, finie ou non, d'objets abstraits deux à deux distincts $T_1$, ... $T_n$, ... on définit un système catégorique d'écoulements en se donnant, pour chaque indice $k$, un sous-monoïde $M_k$ de $M$, et pour chaque indice $k$ tel que $k+1$ soit également un indice de la suite, un sous-monoïde $N_k$ de $M$ contenant $M_k$ et $M_{k+1}$ : les flèches de $T_k$ dans lui-même s'identifient à $M_k$, la composition étant donnée par $f\circ g=g*f$, les flèches de $T_k$ dans $T_{k+1}$ s'identifient à $N_k$, toutes les autres flèches étant engendrées par  composition, définie de même par la loi $*$ de $M$. Notons\label{notation exm categorie ecoulements par monoides} $\mathbf{E}[(M,*); M_1, N_1, M_2,...M_n,...]$ la catégorie obtenue. L'exemple \ref{exemple ecoulements RN} correspond alors à  \[\mathbf{E}[(\mathbf{R},+); \mathbf{R}_+, \mathbf{R}_+, \mathbf{N}].\]
\end{exm}

\begin{exm} \label{ecoulements RZ12} On constitue un système d'écoulements à deux modes en reliant les deux monoïdes $(\mathbf{R}_+,+)$ et $(\mathbf{Z}/{12\mathbf{Z}},+)$ par des flèches $\mathbf{R}_+ \rightarrow \mathbf{Z}/{12\mathbf{Z}}$ de la forme $ f_r:x \mapsto E(x+r) \mod 12 $, la composition de trois flèches $s:\mathbf{R}_+\rightarrow\mathbf{R}_+$, $f_r:\mathbf{R}_+ \rightarrow \mathbf{Z}/{12\mathbf{Z}}$ et $n:\mathbf{Z}/{12\mathbf{Z}}\rightarrow\mathbf{Z}/{12\mathbf{Z}}$ étant définie par
\[ 
n\circ f_r \circ s = f_{s+r+n}: x \mapsto E(x+s+r)+n \mod 12
\]

\end{exm}


\begin{df}[régularité, irréversibilité]
Un écoulement $f$ est dit régulier à droite si pour tous écoulements $g,h$ composables à gauche avec $f$, on a $g\circ f=h\circ f \Rightarrow g=h$. Un système de temporalités est dit régulier à droite si tout écoulement l'est. On définit de même les écoulements et les systèmes de temporalités réguliers à gauche. Un écoulement inversible, c'est-à-dire un \textit{iso},  est aussi dit \emph{réversible}. Dans le cas contraire, on le dit \emph{irréversible}.
\end{df}

\section{Transitions}

\begin{notations}\label{notations transitions}

Pour tout ensemble $A$, on note $\mathcal{P}A$ l'ensemble des parties de $A$. 
Pour toute application $f:A\rightarrow B$, on note $f^\mathcal{P}$ l'application de $\mathcal{P}A$ dans $\mathcal{P}B$ définie pour toute partie $U$ de $A$ par $f^\mathcal{P}(U)=\{f(u), u\in U\}$. 
Pour toute application $f:A\rightarrow \mathcal{P}B$, on note ${^\mu}f$ l'application de $\mathcal{P}A$ dans $\mathcal{P}B$ définie pour toute partie $U$ de $A$ par ${^\mu}f(U)=\bigcup_{u\in U} f(u)$.
\end{notations}

\begin{df}\label{def categorie transitions} On appelle \emph{catégorie des transitions} la catégorie $\mathbf{P}$ 
\begin{itemize}
\item dont les objets sont les ensembles,
\item telle que, pour tout couple d'ensembles $(A,B)$, les flèches de $A$ vers $B$, appelées \emph{transitions} de $A$ vers $B$, sont les applications de $A$ vers $\mathcal{P}B$,
\item et telle que la composée de deux transitions $f\in \mathbf{P}(A,B)$ et $g\in \mathbf{P}(B,C)$, notée $g\odot f$, est donnée par 
\[ g\odot f={^\mu}g \circ f. \]
\end{itemize}
Les flèches de $\mathbf{P}$ sont appelées les \emph{transitions}.
\end{df}

\begin{notations}\label{notation fleches transitions}
On notera $f:A\rightsquigarrow B$ pour exprimer que $f$ est une transition de $A$ vers $B$.
\end{notations}

\begin{rmq} La catégorie $\mathbf{P}$ est isomorphe à la catégorie des relations ensemblistes. Elle admet l'ensemble vide pour objet terminal (objet nul), et le produit cartésien y coïncide avec l'union disjointe. Par ailleurs, les flèches de $\mathbf{P}$ peuvent être interprétées comme des distributeurs très particuliers (entre les catégories discrètes associées aux ensembles considérés). 
\end{rmq}

\begin{rmq} Au lieu d'ensembles, on pourrait considérer les objets d'un topos. Nous ne développerons pas ce point de vue ici.
\end{rmq}


\begin{df} Une transition $p:A\rightsquigarrow B$ est dite \emph{quasi-déterministe} si, pour tout $a\in A$, l'ensemble $p(a)$ a au plus un élé\-ment. Elle est dite \emph{complète} si $p(a)$ n'est vide pour aucun $a\in A$. Elle est dite \emph{déterministe} si elle est à la fois quasi-déterministe et complète. 
\end{df}

\begin{notations}\label{notations transitions deterministes} Une transition quasi-déterministe $p:A\rightsquigarrow B$ s'identifie à une fonction de $A$ vers $B$, et une transition déterministe à une application de $A$ vers $B$. Dans ce cas, commettant un léger abus d'écriture, on notera souvent encore $p$ la fonction correspondante, et donc $p(a)$ l'unique élé\-ment, s'il existe, de la transition $p$ appliquée à $a$.
\end{notations}

\begin{df}\label{definition transition inversible} Une transition déterministe $p:A \rightarrow B$ est dite \emph{réversible} si, en tant qu'application  $p:A \rightarrow B$, elle est bijective.
\end{df}

\section{Dynamiques de type $\mathbf{E}$}

\subsection{Définitions}


\begin{df}[$\mathbf{E}$-dynamiques]\label{definition dynamique de type E} Étant donnée un système d'écoulements $\mathbf{E}$, une \emph{dynamique de type $\mathbf{E}$}, ou \emph{$\mathbf{E}$-dynamique}, est un foncteur de $\mathbf{E}$ dans $\mathbf{P}$. 
\end{df}

\begin{notations}\label{notation image par une dynamique} L'image par une $\mathbf{E}$-dynamique $\alpha$ d'un mode d'écoulement $T\in \dot{\mathbf{E}}$ sera le plus souvent notée $T^\alpha$.
De même, l'image par $\alpha$ d'un écoulement $d\in\vec{\mathbf{E}}$ sera noté $d^\alpha$ plutôt que $\alpha(d)$. 
\end{notations}


\begin{df}[Etats et transitions d'une dynamique] L'ensemble des \emph{états} de la $\mathbf{E}$-dynamique $\alpha$ est
\[ St_\alpha= \bigcup_{T\in \dot{\mathbf{E}}} {T^\alpha}.\] 
Pour tout écoulement $d\in \vec{\mathbf{E}}$,  $d^\alpha$ est la transition  associée à $d$ par $\alpha$.
\end{df}

L'image par une transition $d^\alpha$ d'un état $e\in \mathrm{dom}(d)^\alpha$ sera également désignée comme l'action de $d$ sur $e$ dans (ou pour) la dynamique $\alpha$, ou encore, s'il n'y pas d'ambiguïté sur la dynamique en jeu, comme l'action dynamique de $d$ sur $e$.


Certaines constructions\footnote{Notamment celles de la section \ref{section categories dynamiques}.} nécessitent qu'on se limite aux dynamiques pour lesquelles un seul mode d'écoulements peut agir sur chaque état de la dynamique.  D'où la définition suivante.

\begin{df}[Dynamiques propres]\label{definition dynamiques propres} Une dynamique $\alpha: \mathbf{E}\rightarrow\mathbf{P}$ est dite \emph{propre} si pour tout couple $(S,T)\in \vert E\vert_0$ de types temporels, on a l'implication
\[ S \neq T \Longrightarrow S^\alpha \cap T^\alpha=\emptyset.\]
\end{df}


\begin{prop}\label{prop preodre sur les etats} Étant donnée une dynamique propre $\alpha: \mathbf{E}\rightarrow\mathbf{P}$, la relation $\preceq_\alpha$ définie sur les états de $\alpha$  par
\[ 
\forall a,b \in St_\alpha, a \preceq_\alpha b \Leftrightarrow \exists e\in\vec{\mathbf{E}}, b\in e^\alpha (a)
 \]
est une relation de pré-ordre.
\end{prop}
\noindent \textbf{Preuve}. La relation $\preceq_\alpha$ est réflexive car tout état est laissé invariant par l'action dynamique de l'identité du mode de temporalité (unique) à l'image duquel il appartient.  Cette relation est en outre transitive car, la dynamique ayant été supposée propre, si deux écoulements peuvent agir (dans la dynamique) successivement sur un état c'est nécessairement que ces écoulements se composent (dans la catégorie des écoulements).
\begin{flushright}$\square$\end{flushright} 


\begin{df}[Orbite]
Soit  $\alpha$ une $\mathbf{E}$-dynamique, et $a\in St_\alpha$ un état de $\alpha$. L'\emph{orbite} de $a$ est l'ensemble 
\[Orb_a
=
\bigcup_{f\in \vec{\mathbf{E}},a\in \mathrm{dom}(f^\alpha)} f^\alpha(a).\]
\end{df}

\begin{rmq} Pour une dynamique propre, la relation de pré-ordre $a\preceq_\alpha b$ peut être définie par $b\in Orb_a$. Par contre, pour une dynamique impropre, cette dernière relation n'est pas nécessairement transitive : il peut arriver qu'un état $c$ appartenant à l'orbite d'un état $b$ de l'orbite de $a$ n'appartienne pas à l'orbite de $a$. 
\end{rmq}


\begin{df}[complétude, déterminisme, réversibilité]
Soit  $\alpha$ une $\mathbf{E}$-dynamique. Si, pour tout écoulement $f:T\rightarrow T'$ de $\mathbf{E}$ la transition $f^\alpha$ est déterministe (resp. quasi-déterministe, complète, inversible), on dit que la dynamique $\alpha$ elle-même est \emph{déterministe} (resp. \emph{quasi-déterministe}, \emph{complète}, \emph{inversible}). 
\end{df}

\subsection{Exemples}


\begin{exm}[Dynamiques sur un groupe]\label{exemple : dynamiques sur un groupe} Pour tout groupe $(G,*)$, une $G$-dynamique est une \emph{action} du groupe $G$ sur un ensemble d'états $A$.  En effet, $G$ s'identifie à une catégorie avec un unique objet $*$, l'ensemble des états de la dynamique est donc de la forme $A=*^\alpha$, tandis que pour tout $g\in G$, on a $(g^{-1})^\alpha \odot g^\alpha=g^\alpha\odot(g^{-1})^\alpha=id_A$, d'où l'on déduit que $g^\alpha$ est déterministe (et, en fait, inversible). Ainsi,

\begin{itemize}
\item une $(\mathbf{Z},+)$-dynamique $\alpha$ consiste en la donnée d'un ensemble d'états $A$ et d'une bijection $f=1^\alpha:A\to A$, les autres transitions de $\alpha$ vérifiant $n^\alpha=f^n$;
\item une dynamique sur $(\mathbf{R},+)$ est né\-cessai\-rement déterministe\footnote{Rappelons qu'à ce stade il n'est pas question de la notion d'\emph{instant} : ce sont les écoulements qui décrivent ici $(\mathbf{R},+)$.}, et elle est entièrement déterminée par la connaissance des bijections $f^r:A\to A$  pour $r$ décrivant un intervalle ouvert non vide quelconque. L'exemple classique est le flot d'une équation différentielle autonome sur une variété différentielle admettant pour toutes conditions de Cauchy une solution unique définie sur $\mathbf{R}$ tout entier. A noter que, dans le cadre ensembliste où nous  nous plaçons pour le moment, les trajectoires $r\mapsto f^r(x)$ ne sont pas né\-cessai\-rement différentiables, ni même continues.
\end{itemize}

\end{exm}


\begin{exm} Une $(\mathbf{N},+)$-dynamique consiste en la donnée d'un ensemble d'états $A$ et d'une transition $f:A\to \mathcal{P}A$. Cette application détermine en effet la suite de transitions définie par $f^0=id_A$ et $f^{n+1}=f\odot f^{n}$. 

Les jeux de plateau (go, échecs, dames,...) peuvent être vus comme de telles dynamiques, évidemment non déterministes puisque jeu il y a\footnote{Ces dynamiques seront également incomplètes, dès lors que l'on considère des espaces d'états incluant des situations ne pouvant exister dans ces jeux.}, sur $(\mathbf{N}, +)$.

Par exemple, pour le go, on prend pour ensemble d'états \[A=(\{-1,0,1\}^{(\{1,...19\}^2)})^3\times \{-1,1\},\] et pour tout $e\in A$, $1^\alpha(e)$ sera l'ensemble des états possibles, y compris la mémoire des deux dernières configurations, après qu'on aura appliqué la règle suivante : \emph{le joueur $1$ ou $-1$, selon ce qu'en dit $e$, aura posé sa pierre sur le \emph{goban} en respectant les trois règles fondamentales du go, y compris la règle du \emph{ko}\footnote{La règle du \textit{ko} rend né\-cessai\-re de garder en mémoire les deux configurations précédentes.}. }

\end{exm}

\begin{exm} Pour une $(\mathbf{N},+)$-dynamique quasi-déterministe, on peut remplacer $f$ par une fonction $A\to A$. Pour une dynamique déterministe, $f$ est une application.
\end{exm}

\begin{exm} Une dynamique $\alpha$ sur $(\mathbf{R}_+,+)$ consiste en la donnée d'un ensemble d'états $A$ et, pour tout réel positif $r$, d'une application $f^r=r^\alpha:A\to \mathcal{P}A$, de sorte que $f^r\odot f^s=f^{s+r}$. En voici trois exemples admettant $A=\mathbf{R}$ pour ensemble des états.

\textbf{1.} Pour tout $r\geqslant 0$, et tout $t\in\mathbf{R}$, on pose
\[ f^r(t)=[t,t+r]. \]

\textbf{2.} Pour tout $r\geqslant 0$, on définit $f^r$ par
\begin{itemize}
\item pour $x>0$, $f^r(x)=]max(0,x-r), x+r]$,
\item pour $x=0$, $f^r(0)=\{0\}$,
\item pour $x<0$, $f^r(x)=[x-r, min(x+r,0)[$.
\end{itemize}

\textbf{3.} Dans le même genre, on prend pour $f^r$ la transition définie par
\[ f^r(x)= [max(E(x), x-r), min (x+r, E(x)+1)[,\] où $E(x)$ désigne la partie entière de $x$.

\end{exm}


\begin{exm} Pour tout monoïde (unitaire) $(M,*)$, une $M$-dynamique \emph{détermi\-niste} s'identifie à une \emph{action} de $M$ sur l'ensemble des états $A=*^\alpha$ de cette dynamique.
\end{exm}


\begin{exm} Une dynamique $\alpha$ sur $(\mathbf{N},\leqslant)$ consiste en la donnée d'une suite d'ensembles d'états $(A_n)_{n\geqslant 0}$ et d'une suite \emph{quelconque} de transitions $f_n:A_n \rightsquigarrow A_{n+1}$. En effet, tout écoulement de $(\mathbf{N},\leqslant)$ se décompose de façon unique en écoulements élé\-mentaires $(k,k+1)$ mutuellement indépendants. En particulier, \textit{toute} suite $(u_n)_{n\geqslant 0}$ dans un ensemble quelconque $A$ peut être considérée 
\begin{itemize}
\item soit comme une dynamique déterministe avec pour ensembles d'états les singletons $A_n=\{u_n\}$,
\item soit comme une dynamique quasi-déterministe avec pour ensembles d'états $A_n=A$ et pour transitions les fonctions $f_n$ dont le domaine de définition est réduit à $\{u_n\}\subset A$, avec $f_n(u_n)=u_{n+1}$.
\end{itemize}

\end{exm}

\begin{exm} Une dynamiques sur $(\mathbf{Z},\leqslant)$ consiste en la donnée d'une suite  \textit{quelconque} indexée par $n\in\mathbf{Z}$ de transitions $f_n:A_n \rightsquigarrow A_{n+1}$.
\end{exm}

\begin{exm}\label{dynamique sur R ordonne} Une dynamiques $\gamma$ sur $(\mathbf{R},\leqslant)$ consiste en la donnée d'une famille d'ensembles $\Gamma_r=r^\gamma$ indexée par $r\in\mathbf{R}$, et d'une famille de transitions $\gamma_{(r\leqslant s)}=(r\leqslant s)^\gamma:\Gamma_r \rightsquigarrow \Gamma_s$, indexée par les couples $(r,s)$ décrivant le demi-plan $r\leqslant s$ de $\mathbf{R}^2$, vérifiant
\[ 
\gamma_{(s\leqslant t)}\odot \gamma_{(r\leqslant s)}=\gamma_{(r\leqslant s)}\gamma_{(s\leqslant t)}=\gamma_{(r\leqslant t)}.
\]

Remarquons que, $\epsilon>0$ étant donné, une telle dynamique est déterminée par la donnée des transitions de la forme $\gamma_{(r\leqslant s)}$ où $s$ décrit $\mathbf{R}$ et $r\in ]s-\epsilon,s[$.

A noter que, pour tout ensemble $\Gamma$, \textit{tout} chemin $g:\mathbf{R}\to \Gamma$ peut être considéré comme une dynamique déterministe en prenant $\Gamma_r=\{g(r)\}$(ou comme une dynamique quasi-déterministe avec $\Gamma$ pour ensemble constant d'états).

\end{exm}
 
\begin{exm} On prend $\mathbf{E}=\mathbf{E}[(\mathbf{R},+); \mathbf{R}_+, \mathbf{R}_+, \mathbf{N}]$ (voir les exemples \ref{exemple ecoulements RN} et \ref{exemple ecoulements MNMN} page \pageref{exemple ecoulements RN}).
Une $\mathbf{E}$-dynamique consiste alors en la donnée de deux ensembles d'états $A$ et $B$ et
\begin{itemize}
\item pour tout réel positif $r$, d'une transition $f^r:A \rightsquigarrow A$,
\item pour tout réel positif $t$, d'une transition $g_t:A \rightsquigarrow B$,
\item d'une transition $h: B\rightsquigarrow B$,
\end{itemize}
de sorte que les relations suivantes soient satisfaites: 

\[ \forall r, s \in \mathbf{R}_+, f^r\odot f^s = f^{s+r},  \]
\[ \forall r, t \in \mathbf{R}_+, \forall n\in\mathbf{N},  h^n\odot g_t\odot f^r = g_{r+t+n}.  \]

Puisque $g_t=g_0\odot f^t$, la donnée de $g_0$ détermine tous les $g_t$, et la définition d'une $\mathbf{E}$-dynamique se ramène à la donnée des $f^r$, de $g_0$ et de $h$ vérifiant

\[ \forall r, s \in \mathbf{R}_+, f^r\odot f^s = f^{s+r},  \]
\[h\odot g_0=g_0\odot f^1.\]

Par exemple, on définit une telle dynamique en prenant $A=B=\mathbf{R}$, et

\textbf{1.} $f^r(x)=x+r$, $g_0(x)=]-\infty,x]$ et $h(x)=x+1$.

\textbf{2.} $f^r(x)=x+r$, $g_0(x)=]-\infty,x]$ et $h(x)=[x-1, x+1]$.

\textbf{3.} $f^r(x)=x+r$, $g_0(x)=0$ et $h(x)=[-2 \vert x\vert, \vert x\vert]$.

\textbf{4.} $f^r(x)=[x-r,x+r]$, $g_0(x)=0$ et $h(x)=2x$.

\textbf{5.} $f^r(x)=[x-r,x+r]$, $g_0(x)=0$ et $h(x)=[-2 \vert x\vert, \vert x\vert]$.

\textbf{6.} $f^r(x)=[x-r,x+r]$, $g_0(x)=\mathbf{R}$ et $h(x)=2x$.

\end{exm}



\subsection{Instants essentiels, existentiels, etc...}

Jusqu'à présent, nous avons considéré les dynamiques sans faire appel à la notion d'instant. Dans cette section, nous associons des instants aux systèmes catégoriques d'écoulements : l'idée est que les instants d'un tel système $\mathbf{E}$ sont les états de dynamiques déterministes spécifiques, en particulier celles que nous appellerons respectivement la 
\emph{dynamique essentielle} et la \emph{dynamique existentielle} de $\mathbf{E}$, conduisant respectivement à la définition des instants essentiels et des instants existentiels. Plus généralement, nous définissons les $\alpha-$instants de $\mathbf{E}$ pour toute dynamique déterministe $\alpha$ sur $\mathbf{E}$.

\begin{df} [$\alpha$-instants] Soit $\alpha$ une dynamique propre déterministe sur un système catégorique d'écoulements $\mathbf{E}$. On appelle \emph{$\alpha$-instants} de $\mathbf{E}$ les états de $\alpha$. 
\end{df}

\begin{rmq} Comme pour tout dynamique propre, l'ensemble des $\alpha$-instants est naturellement muni d'un pré-ordre (voir la proposition \ref{prop preodre sur les etats} page \pageref{prop preodre sur les etats}).
\end{rmq}

\begin{df}[Dynamique et instants essentiels]\label{def instants essentiels} La \emph{dynamique essentielle} associée à un système d'écoulements $\mathbf{E}$ est la $\mathbf{E}$-dynamique déterministe $\zeta_\mathbf{E} =\zeta$ ainsi définie :
\begin{itemize}
\item pour tout mode d'écoulement $T\in \dot{\mathbf{E}}$, $T^\zeta=\{T\}$,
\item pour tout écoulement $(f:S\rightarrow T)\in \vec{\mathbf{E}}$, $f^\zeta$ est défini par $f^\zeta (S)=\{T\}$.
\end{itemize}
Les $\zeta$-instants de $\mathbf{E}$, autrement dit les modes d'écoulement de $\mathbf{E}$, sont appelés les \emph{instants essentiels} de $\mathbf{E}$.
\end{df}

Dans la définition suivante, on note\footnote{voir notations \ref{notations fleches et objets} page \pageref{notations fleches et objets}.} $\{\rightarrow S\}$ la classe des flèches de but $S$ dans la catégorie considérée, où $S$ désigne un objet de ladite catégorie.

\begin{df} [Dynamique et instants existentiels] La \emph{dynamique existentielle} associée à un système d'écoulements $\mathbf{E}$ est la $\mathbf{E}$-dynamique déterministe $\xi_\mathbf{E} =\xi$ ainsi définie :
\begin{itemize}
\item pour tout mode d'écoulement $T\in \dot{\mathbf{E}}$, 
$T^\xi=\{\rightarrow T\}$,
\item pour tout écoulement $(f:S\rightarrow T)\in \vec{\mathbf{E}}$, $f^\xi$ est définie pour tout état $a\in S^\xi$ par $f^\xi (a)=\{f\circ a\}$.
\end{itemize}
Les $\xi$-instants de $\mathbf{E}$ sont appelés les \emph{instants existentiels} de $\mathbf{E}$.
\end{df}

\begin{exm} Soit $\mathbf{E}=(\mathbf{R}_+,+)$. 

\begin{itemize}
\item La dynamique $\zeta_\mathbf{E}$ admet un unique instant essentiel, notons-le $\bullet$, et l'identité pour unique transition : pour tout écoulement $r\in \mathbf{R}_+$, \[ r^\zeta (\bullet)= \bullet.\]
\item Les instants existentiels s'identifient aux réels positifs $s\in \mathbf{R}_+$, et la transition associée à l'écoulement $r$ agit sur tout instant $s$ selon \[  r^\xi(s)= s+r.\]
\item \label{instants reels} On définit les \emph{instants réels} pour $(\mathbf{R}_+,+)$ comme les états de la dynamique déterministe $\rho$ définie par  
\begin{itemize}
\item ses états : ${\mathbf{R}_+}^\rho=\mathbf{R}$,
\item  ses transitions : pour tout écoulement  $r\geqslant 0$, la transition $r^\rho$
    est l'application $\mathbf{R}\rightarrow\mathbf{R}$
     telle que $r^\rho(t)=t+r$ pour tout $t\in\mathbf{R}$. 
\end{itemize} 
\end{itemize}

\end{exm}

\begin{exm}
Soit $\mathbf{E}=(\mathbf{R}_+,\leqslant)$.
\begin{itemize}
\item Les instants essentiels sont les $r\in\mathbf{R}_+$, et la transition associée à l'écoulement\footnote{Voir la notation de l'exemple \ref{notations r inf s} page \pageref{notations r inf s}.} $(r\leqslant s)$ n'agit que sur le seul instant $r$, selon \[(r\leqslant s)^\zeta (r)=s.\]
\item Les instants existentiels de mode $r\in \mathbf{R}_+$ sont de la forme $(t\leqslant r)$ et constituent un ensemble  $[0,r] \times \{r\}\simeq [0,r]$. La transition associée à l'écoulement $(r\leqslant s)$ est définie pour tout $t\in[0,r]$ par \[  (r\leqslant s)^\xi (t\leqslant r)= (t\leqslant s),\] ce que nous noterons $(r\leqslant s)^\xi : [0,r]\ni t \mapsto t\in [0,s]$.
\end{itemize}
\end{exm} 

\begin{exm} [Instants essentiels et existentiels de ${(\mathbf{R},\leqslant)}$]
\label{instants de R ordonne}
Soit $\mathbf{E}=(\mathbf{R},\leqslant)$.
\begin{itemize}
\item Les instants essentiels sont les $r\in\mathbf{R}$, et la transition associée à l'écoulement\footnote{Voir la notation de l'exemple \ref{notations r inf s} page \pageref{notations r inf s}.} $(r\leqslant s)$ n'agit que sur le seul instant $r$, selon \[(r\leqslant s)^\zeta (r)=s.\]
\item Les instants existentiels de mode $r\in \mathbf{R}$ sont de la forme $(t\leqslant r)$ et constituent un ensemble  $]-\infty,r] \times \{r\}\simeq ]-\infty,r]$. La transition associée à l'écoulement $(r\leqslant s)$ est définie pour tout $t\in]-\infty,r]$ par \[  (r\leqslant s)^\xi (t\leqslant r)= (t\leqslant s),\] ce que nous noterons $(r\leqslant s)^\xi : ]-\infty,r]\ni t \mapsto t\in ]-\infty,s]$.
\end{itemize}
\end{exm}

\subsection{$\mathbf{E}$-dynamorphismes}

On se propose de constituer une catégorie dont les objets seront les $\mathbf{E}-$dynamiques. Plusieurs idées peuvent venir à l'esprit pour la définition des morphismes :  transformations naturelles entre dynamiques considérées comme des foncteurs, endofoncteurs de la catégorie des transitions $\mathbf{P}$ ou bien foncteurs entre catégories images des dynamiques conduisant à des triangles commutatifs. Mais c'est une définition encore différente que nous choisissons, fondée comme on va voir sur l'\textit{inclusion} plutôt que sur l'égalité, en particulier parce que nous avons en vue une définition de la notion de \textit{solutions} d'une dynamique. En outre, nous rencontrerons une situation analogue avec l'inclusion d'une représentation connective dans une autre.

\begin{df}[Inclusion des transitions] Étant donnés deux ensembles $G$ et $D$, et $u:G\leadsto D$ et $v:G\leadsto D$ deux transitions, on dit que $u$ est incluse dans $v$, et l'on note
\[ u\subset v \]
pour exprimer le fait que, pour tout élé\-ment $g\in G$, on ait
\[ u(g) \subset v(g). \]
\end{df}

\begin{df}[$\mathbf{E}$-dynamorphismes]\label{def E dynamo} Soit  $\alpha$ et $\beta$ deux $\mathbf{E}$-dynamiques.
Un  \emph{$\mathbf{E}$-dynamorphisme} $\delta : \alpha \looparrowright\beta$ consiste en la donnée, pour tout mode d'écoulement $S\in \dot{\mathbf{E}}$, d'une transition $\delta_S:S^\alpha \rightsquigarrow S^\beta$ de sorte que pour tout écoulement $(f:S\rightarrow T)\in \vec{\mathbf{E}}$ on ait
\[ 
{\delta_T}\odot f^\alpha \subset\, {f}^\beta \odot \delta_S.
 \]
\end{df}

Autrement dit, $\delta $ est un $\mathbf{E}$-dynamorphisme de $\alpha$ vers $\beta$ si, pour tout écoulement $f:S\rightarrow T$ et  pour tout état $s\in S^\alpha$, on a l'inclusion
\[ 
 ^\mu{\delta_T}(f^\alpha(s)) \subset\, ^\mu{f}^\beta (\delta_S(s)).
 \]
 
On vérifie sans difficulté que,  munie des $\mathbf{E}$-dynamor\-phismes, la classe des $\mathbf{E}$-dynamiques constitue une catégorie. On la notera $\mathbf{E-Dy}$. \label{categorie E-Dy}

\begin{rmq}
Cette catégorie admet la dynamique vide comme objet nul.
\end{rmq}

\begin{df} Un $\mathbf{E}$-dynamorphisme $\delta:\alpha \looparrowright\beta$ est dit \emph{complet} (resp. \emph{quasi-déterministe}, resp. \emph{déterministe}) si, pour tout mode d'écoulement $T$ et tout état $t\in T^\alpha$, $\delta_T(t)$ est une partie non vide (resp. a au plus un élé\-ment, resp. est un singleton) de  $T^\beta$.
\end{df}

 \begin{rmq}[Topos des $\mathbf{E}$-dynamiques déterministes] \label{topos des E-dynamiques deterministes} Étant donnée une petite catégorie quelconque $\mathbf{E}$, la classe des $\mathbf{E}$-dynamiques déterministes et la classe des $\mathbf{E}$-dynamorphismes \emph{déterministes} entre ces dynamiques constituent un topos. En effet, dans ce cas, l'inclusion qui préside à la définition des dynamorphismes devient une égalité, et la catégorie ainsi définie coïncide avec la catégorie des préfaisceaux sur (la catégorie opposée à) $\mathbf{E}$. En particulier, lorsque $\mathbf{E}$ est un monoïde, il s'agit du topos des actions de $\mathbf{E}$. Si $\mathbf{E}$ est un groupe, les objets de ce topos sont toutes les $\mathbf{E}$-dynamiques. 
\end{rmq}
 
\subsection{Solutions d'une dynamique}

\begin{df}[solutions d'une dynamique] \'{E}tant donnée une $\mathbf{E}$-dynamique propre déterministe $\tau$, une \emph{$\tau$-solution} $\sigma$ d'une $\mathbf{E}$-dynamique $\alpha$ est un $\mathbf{E}$-dyna\-mor\-phisme quasi-déterministe $\sigma:\tau \looparrowright \alpha$. Si $\sigma$ est un dynamorphisme complet, on dit que la solution est \emph{complète}. 
Si $\tau=\zeta$, on dit que $\sigma$ est une \emph{solution essentielle} de $\alpha$, et si $\tau=\xi$, on dit que c'est une \emph{solution existentielle}.   
\end{df}

\begin{rmq} Si une \emph{$\tau$-solution} $\sigma$ d'une $\mathbf{E}$-dynamique $\alpha$ est vide à un $\tau$-instant $s\in A^\tau$, alors elle est également vide à tout instant ultérieur $t\succeq_\tau s$. En effet, en notant $e:A\rightarrow B$ un écoulement tel que $t=s^\delta(s)$, on a, par définition d'un dynamorphisme 
\[ 
\sigma^B(e^\tau(s))\subset e^\alpha(\sigma^A(s)),
 \]
\noindent d'où
\[ 
\sigma^B(t)=\emptyset.
\]

\end{rmq}

\begin{exm}[Solutions vides d'une dynamique $\alpha$] Pour toute dynamique propre déterministe $\tau$, le dynamorphisme vide $\tau \looparrowright \alpha$ constitue une solution (évidem\-ment non complète), dite solution vide. 
\end{exm}

\begin{exm}[solutions essentielles et existentielles en temps continu irréversible] On prend $\mathbf{E}=(\mathbf{R}_+,+)$. Soit $\alpha$ une dynamique de type $\mathbf{E}$, et soit $S={\mathbf{R}_+}^\alpha$ l'ensemble de ses états. 

Une solution existentielle de $\alpha$ consiste en la donnée d'une fonction $f : \mathbf{R}_+ \rightarrow S$ telle que pour tous réels positifs $t$ et $r$, on ait $f(t+r)\in r^\alpha (f(t))$. Cette solution est complète si $f$ est une application définie sur tout $\mathbf{R}_+$.

Une solution essentielle de $\alpha$ consiste en la donnée d'un état $s_0\in S$ tel que pour tout écoulement $r\in \mathbf{R}_+$, $s_0\in r^\alpha (s_0)$.  En particulier, une telle solution essentielle détermine une solution existentielle stationnaire, définie par $f(t)=s_0$ pour tout $t\in \mathbf{R}_+$.

Une solution réelle de $\alpha$ est une $\rho$-solution de $\alpha$ (voir l'exemple \ref{instants reels}), autrement dit c'est une fonction $f$ de $\mathbf{R}$ dans l'ensemble des états de $\alpha$ telle que pour tout $\rho$-instant $t\in\mathbf{R}$ et tout écoulement $r\geqslant 0$, on ait
\[ 
f(t+r)\in {r^\alpha}(f(t)).
 \]
\noindent On retrouve ainsi la notion habituelle de solution définie sur $\mathbf{R}$ d'une dynamique, sans que celle-ci soit pour autant né\-cessai\-rement déterministe, puisque le mode d'écoulement n'est pas le groupe $\mathbf{R}$ mais le semi-groupe $\mathbf{R}_+$. Soulignons que la nette distinction entre écoulements et instants est indispensable à l'établissement d'un tel point de vue.
\end{exm}

\begin{exm} 


Pour $\mathbf{E}=(\mathbf{R}_+,+)$, 
on considère la $\mathbf{E}$-dynamique $\alpha$ définie par $(\mathbf{R}_+)^\alpha=\mathbf{R}^2$ et pour $M\in \mathbf{R}^2$, et $r\in \mathbf{R}_+$, $r^\alpha(M)$ est le disque fermé de centre $M$ et de rayon $r$ :
\[  r^\alpha(M)=\mathcal{D}(M,r).\]

Alors les solutions existentielles complètes de $\alpha$ sont les applications $\sigma:\mathbf{R}_+ \rightarrow \mathbf{R}^2$ telles pour tout instant $s\in\mathbf{R}_+$ et tout écoulement $r\in\mathbf{R}_+$, on a 
$\sigma(s+r)\in \mathcal{D}(\sigma(s),r)$, autrement dit 
\[d(\sigma(s),\sigma(s+r))\leqslant r.
 \]

\noindent Ce sont donc les mouvements plans $\sigma$ de \og vitesse\fg\, inférieure ou égale à 1.

 \end{exm}

\begin{exm}[Solution existentielle canonique de la dynamique essentielle]. \'{E}tant donnés un système catégorique d'écoulements $\mathbf{E}$, sa dynamique essentielle $\zeta$ et sa dynamique existentielle $\xi$, on définit une solution existentielle complète $Z$ de $\zeta$ en posant, pour tout mode d'écoulement $S\in \dot{\mathbf{E}}$, et tout instant existentiel $t\in S^\xi$,
\[ 
{Z}_S(t)=S.
\]
\noindent \`{A} toute solution essentielle $\sigma$ d'une $\mathbf{E}$-dynamique $\alpha$ se trouve alors canoniquememt associée une solution existentielle de la même dynamique, à savoir $\sigma\circ Z$. Les solutions existentielles de cette forme pourront être appelées les solutions existentielles essentielles. Dans le cas des systèmes d'écoulement possédant un seul mode, les solutions existentielles essentielles sont les solutions stationnaires.
\end{exm}

\begin{exm} Soit $\alpha_{(r\leqslant s)}:A_r \rightsquigarrow A_s$ une $\mathbf{E}$-dynamique avec $\mathbf{E}={(\mathbf{R},\leqslant)}$ (voir l'exemple \ref{dynamique sur R ordonne} page \pageref{dynamique sur R ordonne} et l'exemple \ref{instants de R ordonne} page \pageref{instants de R ordonne}). 

Une solution essentielle de  $\alpha$ est la donnée d'une fonction $f:\mathbf{R}\longrightarrow \bigcup_{r\in \mathbf{R}}{A_r}$ telle que, pour tout $r\in \mathbf{R}$, $f(r)\in A_r$ et pour tout $s\geqslant r$, $f(s)\in \alpha_{(r\leqslant s)}(f(r))$.

Une solution existentielle de $\alpha$ est une fonction $g:{\{(r,s)\in\mathbf{R}^2, r\leqslant s\}}\longrightarrow \bigcup_{t\in \mathbf{R}}{A_t}$ telle que, pour tout $(r\leqslant s)$, $g(r,s)\in A_s$ et pour tout $(r\leqslant s\leqslant t)$, $g(r,t)\in (s,t)^\alpha (g(r,s))$.  En particulier, si $\alpha$ est une dynamique déterministe, et si $g$ en est une solution complète, alors $g$ est une application définie sur ${\{(r,s)\in\mathbf{R}^2, r\leqslant s\}}$ qui vérifie
\[ 
g(r,t)= (s,t)^\alpha (g(r,s))
 \]
\noindent pour tout triplet $(r,s,t)$ tel que $r\leqslant s \leqslant t$.

Par exemple, la solution existentielle canonique $Z$ de la dynamique essentielle $\zeta$ de ${(\mathbf{R},\leqslant)}$ est définie pour tout $(r\leqslant s)$ par $Z(r,s)=s$, et l'on a bien $t=(s,t)^\zeta (s)$ pour tout triplet $(r,s,t)$ tel que $r\leqslant s \leqslant t$.

A noter qu'une solution complète $g$ d'une $\mathbf{E}$-dynamique déterministe $\alpha$ est entièrement déterminée par les valeurs qu'elle prend sur la diagonale $\{(r,r), r\in \mathbf{R}\}$, puisque

\[ g(r,s)=(r,s)^\alpha g(r,r). \]

\noindent Par contre, les $g(r,r)$ ne sont pas, en général, déterminées par d'autres valeurs. En munissant l'ensemble ${\{(r,s)\in\mathbf{R}^2, r\leqslant s\}}$ du préordre $(r_1,s_1)\preceq (r_2,s_2) \Leftrightarrow s_1\leqslant s_2$, on peut interpréter ce qui précède de la façon suivante : chaque instant réel $r$ est constitué de tous les instants existentiels de la forme $(q,r)$ avec $q\leqslant r$, où $q$ représente en quelque sorte la trace à l'instant $r$ de l'instant antérieur $q$, et une dynamique déterministe ne détermine en fait à chaque instant que ce qui dans cet instant est une trace du passé, la pointe de présent pur $(r,r)$ de l'instant $r$ étant quant à lui le lieu du libre choix de nouvelles conditions initiales pour une telle dynamique.

\end{exm}


\begin{df}[trajectoire d'une solution]\label{trajectoire d'une solution} La \emph{trajectoire} d'une $\tau$-solution $\sigma$ d'une $\mathbf{E}$-dynamique $\alpha$ est l'ensemble $\{\sigma_{T}(t), T\in \vert\mathbf{E}\vert_0, t\in{T^\tau}\}$, où l'on a identifié la transition quasi-déterministe $\sigma_{T}$ à une fonction ${T^\tau}\rightarrow{T^\alpha}$. La trajectoire d'une $\tau$-solution sera appelée \emph{$\tau$-trajectoire} de la dynamique. En particulier, les \emph{trajectoires existentielles} (resp. essentielles) sont les trajectoires des solutions existentielles (resp. essentielles).
\end{df}


En général, les trajectoires d'une dynamique ne coïncident pas avec ses orbites.


\section{Catégorie des dynamiques : $\mathbf{Dy}$}

\subsection{Définition de $\mathbf{Dy}$}

Nous allons voir que l'association à toute petite catégorie de sa dynamique existentielle, ou de sa dynamique essentielle, est une opération fonctorielle. Un tel constat n'a de sens que si l'on commence par définir la catégorie de toutes les dynamiques catégoriques ensemblistes, indépendamment du choix d'une catégorie d'écoulements. L'objet de la présente section est de définir une telle catégorie, que nous noterons $\mathbf{Dy}$. Entre autres choses, cela nous permettra de définir les $\tau$-solutions d'une dynamique $\alpha$, avec $\tau$ et $\alpha$ de types éventuellement différents.

\begin{df}\label{categorie Dy} La catégorie  $\mathbf{Dy}$ des dynamiques catégoriques ensemblistes a
\begin{itemize}
\item pour objets : les $\alpha$ tels qu'il existe une petite catégorie $\mathbf{E}$, né\-cessai\-rement unique, telle que $\alpha$ soit une $\mathbf{E}$-dynamique,
\item pour morphismes d'une $\mathbf{E}$-dynamique $\alpha$ vers une $\mathbf{F}$-dynamique $\beta$ les couples $(\Delta,\delta)$ où
\begin{itemize}
\item  $\Delta$ est un foncteur $\mathbf{E}\rightarrow \mathbf{F}$,
\item $\delta$ consiste en la donnée, pour tout $\mathbf{E}$-mode d'écoulement $T$, d'une transition $\delta_T: T^\alpha \rightsquigarrow (\Delta T)^\beta$,
\end{itemize}
\noindent tels que pour tout $\mathbf{E}$-écoulement $e:S\rightarrow T$ on ait
\[ 
\delta_T\odot e^\alpha \subset (\Delta e)^\beta \odot \delta_S
 \]
\noindent autrement dit, tel que pour tout état $s\in S^\alpha$, on ait
\[ 
^\mu{\delta_T}(e^\alpha(s))\subset\, {^\mu{(\Delta e)}}^\beta (\delta_S (s)).
\]
\end{itemize}
Ces morphismes seront appelés des \emph{dynamorphismes}.
\end{df}

\begin{rmq} Étant donnée  $\alpha$ une dynamique quelconque, on notera parfois $\mathbf{E}_\alpha$ sa catégorie source, de sorte que $\alpha\mapsto \mathbf{E}_\alpha$ est la partie objet du foncteur d'oubli $\mathbf{Dy}\longrightarrow \mathbf{Cat}$, dont l'action sur les flèches est donnée par $(\Delta,\delta)\mapsto \delta$.
\end{rmq}

\begin{rmq} La catégorie des $\mathbf{E}$-dynamiques se plonge naturellement dans $\mathbf{Dy}$, en identifiant tout $\mathbf{E}$-dynamorphisme $\delta$ au dynamorphisme $(Id_\mathbf{E}, \delta)$.
\end{rmq}

\begin{df}\label{divers styles de dynamo} Un dynamorphisme $(\Delta,\delta)$ est dit fidèle si le foncteur $\Delta$ est fidèle. Il est dit quasi-déterministe (resp. complet, resp. déterministe) si toutes les transitions $\delta_T$ sont quasi-déterministes (resp. complètes, resp. déterministes).
\end{df}

\begin{df}[$\lambda$-solutions d'une dynamique] Soit $\alpha$ une dynamique de type $\mathbf{A}$, et $\lambda$ une dynamique déterministe de type $\mathbf{L}$. On appelle $\lambda$-solution de $\alpha$ tout dynamorphisme quasi-déterministe de $\lambda$ vers $\alpha$. Selon les propriétés de ce dynamorphisme, la solution sera dite complète, fidèle, etc... Si $\lambda$ est la dynamique essentielle (resp. existentielle) de $\mathbf{L}$, une $\lambda$-solution sera aussi dite $\mathbf{L}$-solution essentielle (resp. existentielle).
\end{df}

\begin{exm} Soit $\mathbf{F}$ un système catégorique d'écoulements, $\beta$ une dynamique de type $\mathbf{F}$, $\mathbf{E}$ une sous-catégorie de $\mathbf{F}$ et  $\alpha$ la restriction de $\beta$ à $\mathbf{E}$. On obtient alors un dynamorphisme déterministe $(\Delta,\delta)$  de $\alpha$ vers $\beta$ en prenant pour $\Delta$ l'injection de $\mathbf{E}$ dans $\mathbf{F}$ et pour $\delta$ l'identité sur chaque ensemble d'états $T^\alpha$ de $\alpha$.
\end{exm}

\begin{exm} On prend $\mathbf{E}=(\mathbf{N},+)$ et $\mathbf{F}=(\mathbf{R_+},+)$. Une $\mathbf{E}$-dynamique $\alpha$ est caractérisée par la donnée d'un ensemble d'états $A$ et une transition $1^\alpha=u:A\rightsquigarrow A$, de sorte que l'image par $\alpha$ de l'écoulement $n$ est $u^n=u\odot ... \odot u$. Une $\mathbf{F}$-dynamique $\beta$ est caractérisée par la donnée d'un ensemble d'états $B $ et pour chaque $r\geqslant 0$ d'une transition $r^\beta:B \rightsquigarrow B $ de sorte que $s^\beta \odot r^\beta =(r+s)^\beta$. Un dynamorphisme $(\Delta,\delta):\alpha \looparrowright \beta$ consiste alors en un foncteur $\Delta:(\mathbf{N},+)\rightarrow (\mathbf{R_+},+)$, c'est-à-dire un morphisme de monoïdes de la forme $\Delta(n)=n\tau$ avec $\tau\geqslant 0$ constante réelle, et en une transition $\delta:A \rightsquigarrow B $, de sorte que
\[ 
\delta\odot u \subset \tau^\beta\odot \delta,
 \]

\noindent autrement dit, pour tout état $a\in A$,
\[ 
^\mu\delta(u(a))\subset ^\mu{\tau}^\beta(\delta(a)).
\]
 
Cette condition est en effet suffisante, puisqu'on a 
\[\delta\odot u \subset \tau^\beta\odot \delta\linebreak
\Rightarrow \delta\odot u \odot u \subset \tau^\beta\odot \delta\odot u  \subset \tau^\beta\odot \tau^\beta\odot \delta\linebreak
\]

\noindent soit $\delta \odot u^2 \subset (2\tau)^\beta \odot \delta$, 
puis par récurrence $\delta \odot u^n \subset (n\tau)^\beta \odot \delta$ 
pour tout écoulement $n\in \dot{\mathbf{N}}$, ce qui est la condition requise par définition d'un dynamorphisme.

Dans le cas particulier où $\alpha$ est la dynamique existentielle de $(\mathbf{N},+)$, où $\beta$ est une dynamique déterministe et où le dynamorphisme $\delta$ est lui-même déterministe, la condition précédente s'écrit
\[ 
\delta_{n+1}=\tau^\beta(\delta_n),
 \]
\noindent la suite $(\delta_n)_{n\in\mathbf{N}}$ d'états appartenant à $B $ constituant alors, si $\tau>0$, une solution discrète de la dynamique continue $\beta$.
 
\end{exm}

\begin{exm} On prend $\mathbf{E}=(\mathbf{N},\leqslant)$ et $\mathbf{F}=(\mathbf{N},+)$.
Une $\mathbf{E}$-dynamique $\alpha$ est définie par la donnée d'une famille d'ensembles $(A_n)_{n\in\mathbf{N}}$
et, pour tout entier $n$, d'une transition $u_n$
$:A_n\rightsquigarrow A_{n+1}$, les autres transitions se déduisant de celles-ci par composition. 
Une  $\mathbf{F}$-dynamique $\beta$ est définie par un ensemble d'états $B $ et une transition 
$f:B \rightsquigarrow B $ d'où se déduisent par itération les autres transitions. Un dynamorphisme $(\Delta,\delta)$ de $\alpha$ vers $\beta$ consiste alors en un foncteur $\Delta:{(\mathbf{N},\leqslant)}\longrightarrow {(\mathbf{N},+)}$ --- caractérisé par la valeur entière $\Delta_n$ attribuée par $\Delta$ à chaque flèche élé\-mentaire $(n\leqslant n+1)$ de $(\mathbf{N},\leqslant)$ --- et, pour chaque entier $n$, d'une transition $\delta_n:A_n\rightsquigarrow B $, tels que pour tout $n\in\mathbf{N}$ et tout $a\in A_n$,
\[ 
{^\mu}{\delta_{n+1}}(u_n(a))\subset{^\mu}(f^{\Delta_n})(\delta_n(a)).
 \]

En particulier, si les dynamiques considérées sont déterministes, les dynamorphismes déterministes $(\Delta,\delta)$ entre elles sont caractérisés par la  donnée d'une suite d'entiers positifs $N_n$ et d'une famille d'applications $\delta_n:A_n\rightarrow B $ telles que
\[ 
{\delta_{n+1}}\circ u_n= f^{N_n}\circ\delta_n.
 \]
Par exemple, si $\alpha$ désigne la $\mathbf{E}$-dynamique essentielle et $\beta$ la $\mathbf{F}$-dynamique existentielle, de sorte que $u_n:A_n=\{n\}\rightarrow\{n+1\}$ et $f:B=\mathbf{N}\ni n\mapsto n+1=f(n)$, les dynamorphismes déterministes $\alpha\looparrowright\beta$ s'identifient aux suites croissantes (au sens large) d'entiers positifs $(\beta_n)_{n\in\mathbf{N}}$. On constate ainsi que la $(\mathbf{N},+)$-dynamique existentielle admet une infinité (continue) de $(\mathbf{N},\leqslant)$-solutions essentielles.
\end{exm}

\subsection{Foncteur \og dynamique essentielle\fg\, (resp. existentielle)}

On sait que la dynamique essentielle et la dynamique existentielle d'une petite catégorie sont des dynamiques déterministes propres. Le théorème suivant précise la nature fonctorielle de ces notions.

Notons $\mathbf{DyP}$\label{DyP} la sous-catégorie pleine de $\mathbf{Dy}$ constituée des dynamiques propres, et $\mathbf{DyPd}$ la sous-catégorie de $\mathbf{DyP}$ obtenue en ne conservant que les dynamorphismes déterministes. 

\begin{thm}
L'association à tout petite catégorie de sa dynamique essentielle (resp. existentielle) se prolonge en un foncteur $\zeta$ (resp. $\xi$) de la catégorie $\mathbf{Cat}$ des petites catégories dans la catégorie $\mathbf{DyPd}$. La solution existentielle canonique de la dynamique essentielle définit une transformation naturelle $Z$ de $\xi$ dans $\zeta$.
\end{thm}

\noindent \textbf{Preuve}.
Soit, en effet, $\Delta:\mathbf{E}\rightarrow \mathbf{F}$ un foncteur entre deux petites catégories. La dynamique déterministe essentielle $\zeta_\mathbf{E}$ associe à tout $\mathbf{E}$-écoulement $e:S\rightarrow T$ l'unique application entre singletons $\{S\}\rightarrow\{T\}$. 

On définit alors $\zeta_\Delta$ comme le dynamorphisme \textit{déterministe} $(\Delta,\dot{\delta})$, avec $\dot{\delta}=\dot{\Delta}_{\vert\{\bullet\}}$, autrement dit : $\dot{\delta}$, appliquée à un $\mathbf{E}$-mode d'écoulement quelconque $S$, désigne l'unique application 
entre singletons $\dot{\Delta}_{\vert\{S\}}:\{S\}\rightarrow\{\Delta S\}$. On a alors en effet $\dot{\delta}_T(e^{\zeta_\mathbf{E}}(S))=\dot{\delta}_T(T)=\Delta T$ tandis que $(\Delta e)^{\zeta_\mathbf{F}}(\dot{\delta}_S(S))=(\Delta e)^{\zeta_\mathbf{F}}(\Delta S)=\Delta T$, d'où $\dot{\delta}_T(e^{\zeta_\mathbf{E}}(S))=(\Delta e)^{\zeta_\mathbf{F}}(\dot{\delta}_S(S))$ pour tout $e$, ce qui,
 dans le cas déterministe,  est par définition la condition qui doit être vérifiée par $\dot{\delta}$. La fonctorialité de $\zeta:\mathbf{Cat}\rightarrow\mathbf{Dy}$ ainsi défini se vérifie alors sans aucune difficulté.

De même, la dynamique déterministe existentielle $\xi_\mathbf{E}$ associe à tout $\mathbf{E}$-écoulement $e:S\rightarrow T$ l'application\footnote{Voir notations \ref{notations fleches et objets} page \pageref{notations fleches et objets}.}  $e^{\xi_\mathbf{E}}:\{\rightarrow S\}\rightarrow\{\rightarrow T\}$ définie par $e^{\xi_\mathbf{E}}(s)=e\circ s$. On obtient alors un dynamorphisme \textit{déterministe} $(\Delta,\vec{\delta}):\xi_\mathbf{E}\looparrowright\xi_\mathbf{F}$ en posant
\[ 
\vec{\delta}=\vec{\Delta}
_{\vert \{ \rightarrow\bullet\}},
 \]
\noindent autrement dit en prenant pour tout $\mathbf{E}$-mode d'écoulement $S$, $\vec{\delta}_S:\{\rightarrow S\}\rightarrow\{\rightarrow \Delta S\}$ tel que $\vec{\delta}_S(s)=\Delta s$. On a alors $\vec{\delta}_T(e^{\xi_\mathbf{E}}(s))=\vec{\delta}_T(e\circ s)=\Delta(e\circ s)$ tandis que $(\Delta e)^{\xi_\mathbf{F}}(\vec{\delta}_S(s))=(\Delta e)^{\xi_\mathbf{F}}(\Delta s)=\Delta e\circ \Delta s$, d'où $\vec{\delta}_T(e^{\xi_\mathbf{E}}(s))=(\Delta e)^{\xi_\mathbf{F}}(\vec{\delta}_S(s))$ pour tout $e$ et pour tout $s\in \xi_\mathbf{E}(\mathrm{dom}(e))$, ce qui, dans le cas déterministe,  est par définition la condition qui doit être vérifiée par $\vec{\delta}$. On pose donc $\xi_\Delta=(\Delta,\vec{\delta})$, et la fonctorialité d'un tel $\xi:\mathbf{Cat}\rightarrow\mathbf{Dy}$ se vérifie alors sans plus de difficulté que pour le cas essentiel.

Vérifions enfin que la solution existentielle canonique $Z_\mathbf{E}$ de la dynamique essentielle $\zeta_\mathbf{E}$ d'une petite catégorie quelconque $\mathbf{E}$ définit bien une transformation naturelle du foncteur $\xi$ sur le foncteur $\zeta$. $Z_\mathbf{E}$ est un dynamorphisme $\xi_\mathbf{E}\looparrowright \zeta_\mathbf{E}$ qui, avec, des notations analogues aux précédentes, peut s'écrire

\[ Z_\mathbf{E}=(Id_\mathbf{E}, \mathrm{cod}_{\vert\{\rightarrow\bullet\}}).\]

Il s'agit alors de vérifier que pour tout foncteur $\Delta:\mathbf{E}\rightarrow \mathbf{F}$, on a $Z_\mathbf{F}\circ \xi_\Delta=\zeta_\Delta\circ Z_\mathbf{E}$, où

\[ \xi_\Delta=(\Delta,\vec{\Delta}_{\vert \{\rightarrow\bullet\}} )
\textrm{ et }
\zeta_\Delta=(\Delta,\dot{\Delta}_{\vert \{\bullet\}} ).
\]

Or, on a 
\begin{eqnarray*}
Z_\mathbf{F}\circ \xi_\Delta 
& = 
&(Id_\mathbf{F}, \mathrm{cod}_{\vert\{\rightarrow\bullet\}}) 
\circ 
(\Delta,\vec{\Delta}_{\vert \{\rightarrow\bullet\}} ) \\
& = 
&(Id_\mathbf{F} \circ \Delta,\mathrm{cod}_{\vert\{\rightarrow\bullet\}} \circ \vec{\Delta}_{\vert \{\rightarrow\bullet\}} ) \\
& =
& (\Delta, (\dot{\Delta}\circ \mathrm{cod})_{\vert\{\rightarrow\bullet\}})\\
& =
& (\Delta,\dot{\Delta}_{\vert \{\bullet}\} ) \circ (Id_\mathbf{E}, \mathrm{cod}_{\vert\{\rightarrow\bullet\}}).\\
\end{eqnarray*}

\begin{flushright}$\square$\end{flushright}

\subsection{Catégories dynamiques}\label{section categories dynamiques}

Soit $\alpha:\mathbf{E}\rightarrow\mathbf{P}$ une dynamique propre sur une petite catégorie $\mathbf{E}$. 

\begin{df}\label{definition categorie dynamique d'une dynamique}
 On appelle \emph{catégorie dynamique} de $\alpha$, \emph{catégorie des transitions} de $\alpha$,  ou encore \emph{catégories des états} de $\alpha$,  la petite catégorie $\mathcal{TC}_\alpha$ dont les objets sont les états de $\alpha$ et telle que les flèches d'un objet $a\in \vert{\mathcal{TC}_\alpha}\vert_0 =St_\alpha$ vers un objet $b\in\vert{\mathcal{TC}_\alpha}\vert_0$ sont les triplets $(a,f,b)$ avec $f\in\vec{\mathbf{E}}$ telle que $b\in f^\alpha(a)$, la composition des flèches de $\mathcal{TC}_\alpha$ étant définie par
$(a,f,b)(b,g,c)=(a,fg,c)$.
\end{df}

\begin{rmq} Pour cette définition, il est nécessaire en effet de supposer que la dynamique $\alpha$ est propre, faute de quoi $b\in\mathrm{cod}(f^\alpha)\cap \mathrm{dom}(g^\alpha)\nRightarrow \mathrm{cod}(f)= \mathrm{dom}(g)$
\end{rmq}

\begin{rmq} On pourrait interpréter la catégorie ${\mathcal{TC}_\alpha}$ comme une catégorie des \og chemins\fg\, entre les états de la dynamique $\alpha$, à condition de préciser qu'il s'agit alors de chemins en quelque sorte \og quantiques \fg, en ce sens qu'un  écoulement $(a,f,b)$ ne contient pas en général d'information sur les états \og intermédiaires\fg\, entre $a$ et $b$ (pour autant que la notion ait un sens) : on ne sait pas \og par où c'est passé\fg.
\end{rmq}

\begin{prop} L'association à toute dynamique propre de sa catégorie dynamique se prolonge en un foncteur $\mathbf{DyPd}\rightarrow \mathbf{Cat}$.
\end{prop}

\noindent \textbf{Preuve}. L'image par $\mathcal{TC}$ d'un dynamorphisme déterministe quelconque $(\Delta,\delta):\alpha\rightarrow \beta$ est donnée par $\mathcal{TC}_{(\Delta,\delta)}=D$, où $D$ est le foncteur de la catégorie $\mathcal{TC}_\alpha$ dans la catégorie $\mathcal{TC}_\beta$ ainsi défini :
\begin{itemize}
\item  pour tout objet $s\in \vert\mathcal{TC}_\alpha\vert_0=st(\alpha)$, on pose $D(s)=\delta_S(s)$, où $S$ désigne l'objet --- nécessairement unique puisque $\alpha$ est une dynamique \textit{propre} --- de $\mathbf{E}_\alpha$ tel que $s\in S^\alpha$,
\item pour toute flèche $(s,f,t)\in\vert \mathcal{TC}_\alpha\vert_1$, on pose $Df=(Ds,\Delta f, Dt)$.
\end{itemize}

Vérifions qu'on a bien défini ainsi un $\mathcal{TC}_\beta$-écoulement  $Df:Ds\rightarrow Dt$.
Puisque $(s,f,t)$ est un morphisme de $\mathcal{TC}_\alpha$, on a $t\in f^\alpha(t)$, 
mais $(\Delta,\delta)$ étant un $\mathbf{DyPd}$-morphisme $\alpha\rightarrow\beta$, on a également $\delta_T^\mathcal{P}(f^\alpha(s))\subset \Delta f^\beta(\delta_S(s))$, de sorte qu'on a bien, en particulier, $\delta_T(t)\in (\Delta f)^\beta (\delta_S(s))$, autrement dit $Dt\in(\Delta f)^\beta(Ds)$ et $Df$ est bien défini.

La fonctorialité de $D:\mathcal{TC}_\alpha\rightarrow\mathcal{TC}_\beta$, puis de $\mathcal{TC}:\mathbf{DyPd}\rightarrow \mathbf{Cat}$ se vérifient alors immédiatement.

\begin{flushright}$\square$\end{flushright} 

\paragraph{Toute petite catégorie est dynamique}

Quelles sont les petites catégories qui s'identifient à la catégorie dynamique de $\alpha$ pour une certaine dynamique $\alpha$ ? La proposition suivante montre notamment que c'est le cas de toute petite catégorie.

\begin{prop}\label{prop TCozeta iso 1} Toute petite catégorie est isomorphe à la catégories des transitions de sa dynamique essentielle. Plus précisément, l'endofoncteur $TC\circ\zeta:\mathbf{Cat}\rightarrow \mathbf{Cat}$ est naturellement isomorphe à l'endofoncteur identité $Id_\mathbf{Cat}$.
\end{prop}

\noindent \textbf{Preuve}.
On vérifie sans difficulté que 

\begin{itemize}
\item pour tout $\mathbf{Cat}$-objet $\mathbf{E}$, $TC\circ\zeta (\mathbf{E})$ est obtenue en gardant les mêmes objets que $\mathbf{E}$ et y remplaçant les flèches $f:S\rightarrow T$ par les triplets de la forme $(S,f,T)$,
\item pour tout $\mathbf{Cat}$-morphisme $D:\mathbf{E}\rightarrow \mathbf{F}$, $TC\circ\zeta (D)$ est le foncteur défini
\begin{itemize}
\item pour tout objet $S$ de $TC\circ\zeta (\mathbf{E})$, par $TC\circ\zeta (D)(S)=DS$,
\item pour toute flèche $(S,f,T):S\rightarrow T$, par $TC\circ\zeta (D)(S,f,T)=(DS,Df,DT)$.
\end{itemize}
\end{itemize}

On obtient alors un isomorphisme naturel de $Id_\mathbf{Cat}$ vers $TC\circ\zeta$ en associant à toute petite catégorie $\mathbf{E}$ le foncteur $\mathbf{E}\rightarrow TC\circ\zeta (\mathbf{E})$ qui à tout $\mathbf{E}$-morphisme $f:S\rightarrow T$ associe $(S,f,T):S\rightarrow T$.

\begin{flushright}$\square$\end{flushright} 

\paragraph{Essentialisation}\label{Essentialisation}

Notons $Ess=\zeta\circ\mathcal{TC}$ l'endofoncteur de $\mathbf{DyPd}$ que nous appellerons foncteur d'essentialisation. \`{A} toute dynamique catégorique propre $\alpha:\mathbf{E}\rightarrow\mathbf{P}$, il associe une dynamique déterministe $Ess_\alpha$, dite \emph{essentialisée} de $\alpha$, définie par 
\begin{itemize}
\item pour tout type $s$ de ${\mathcal{TC}_\alpha}$-écoulement, autrement dit pour tout état $s\in St_\alpha$ de $\alpha$, $s^{Ess_\alpha}=\{s\}$,
\item pour tout écoulement $((s,f,t):s\rightarrow t)\in \vec{\mathcal{TC}_\alpha}$, $(s,f,t)^{Ess_\alpha}$ est défini par $(s,f,t)^{Ess_\alpha} (s)=\{t\}$.
\end{itemize}
\begin{rmq} La dynamique essentielle $Ess_\alpha$ étant déterministe, on pourra aussi noter, avec l'abus d'écriture usuel pour les dynamiques déterministes indiquée dans la notation \ref{notations transitions deterministes} page \pageref{notations transitions deterministes} : $(s,f,t)^{Ess_\alpha} (s)=t$.
\end{rmq}

\begin{prop} $Ess^2$ est naturellement isomorphe à $Ess$.
\end{prop}
\noindent \textbf{Preuve}.
$Ess\circ Ess=(\zeta \circ TC)\circ (\zeta \circ TC)
  = \zeta \circ (TC \circ \zeta )\circ TC$,
 d'où, d'après la proposition \ref{prop TCozeta iso 1}, $ Ess\circ Ess\simeq Ess$.
\begin{flushright}$\square$\end{flushright} 

Comme, en général, $Ess_\alpha$ n'est pas isomorphe à $\alpha$, on en déduit que la connaissance de $Ess_\alpha$ ne suffit pas en général à caractériser $\alpha$ (même à isomorphisme près), puisque $Ess_\alpha$ et  $\alpha$ ont même essentialisée.

Par contre, il existe une transformation naturelle $\mathcal{AV}:Ess \rightarrow Id_{\mathbf{DyPd}}$ qui \og recouvre\fg\, toute dynamique propre $\alpha$ par son essentialisée, au sens où tout état de $\alpha$ est dans l'image d'une des transitions qui constituent le dynamorphisme $\mathcal{AV}_\alpha$. 
On définit en effet une telle transformation naturelle en posant 
$\mathcal{AV}_\alpha =(AV_\alpha, av_\alpha)$, 
avec 
\begin{itemize}
\item $AV_\alpha:TC_\alpha\rightarrow\mathbf{E}_\alpha$ le foncteur qui 
\begin{itemize}
\item à tout objet $s\in\vert TC_\alpha\vert_0=st(\alpha)$ associe l'unique\footnote{Rappelons qu'on ne travaille ici qu'avec des dynamiques propres.} $S\in\vert\mathbf{E}_\alpha  \vert_0$ tel que $s\in S^\alpha$,
\item à toute flèche $(s,f,t)\in \vert TC_\alpha\vert_1$, associe $(f:AV_\alpha(s)=S\rightarrow T=AV_\alpha(t))\in \vert\mathbf{E}_\alpha  \vert_1$,
\end{itemize}
\item $av_\alpha$ défini pour tout $s\in st(\alpha)$ par $av_\alpha s:\{s\}\rightarrow (AV_\alpha(s))^\alpha$ telle que $av_\alpha s (s)=s$.
\end{itemize}
Puisque $av_\alpha s(s)=s$, tout état de $\alpha$ est dans l'image d'une des transitions. Par contre, on ne peut pas en déduire que $\mathcal{AV}_\alpha$ soit épique, comme le montre le contre-exemple de la dynamique vide sur une petite catégorie non vide $\mathbf{E}$. 

\begin{exc} Vérifier que ce qui précède définit effectivement une transformation naturelle $\mathcal{AV}:Ess \rightarrow Id_{\mathbf{DyPd}}$.
\end{exc} 

\begin{exc} \`{A} quelles conditions sur $\alpha$ le dynamorphisme $\mathcal{AV}_\alpha$ est-il épique dans $\mathbf{DyPd}$ ?
\end{exc} 

\paragraph{Verticalisation}

On pose $Vert=TC\circ \xi : \mathbf{Cat}\rightarrow \mathbf{Cat}$, et on appelle \og verticalisation \fg\, l'action de ce foncteur.

\begin{exc} Quelle est l'image par l'endofoncteur $Vert$ de la catégorie définie par le monoïde $(\mathbf{R},+)$. Justifier l'appellation \og verticalisation\fg\, donnée à ce foncteur. Déterminer la catégorie $Vert^2 (\mathbf{R},+)$.
\end{exc}

\subsection{$\mathbf{E}$-interprétations}

L'essentialisation des dynamiques permet en particulier, en affinant les temporalités, de \og rendre compte \fg\, de n'importe laquelle d'entre elles par une dynamique \textit{déterministe}.

Dans la présente section, on introduit le concept d'interprétation d'une dynamique par une autre, dans laquelle c'est l'espace des états qui est affiné. L'idée intuitive est qu'une interprétation d'une \og dynamique observée \fg\, $\beta$ par une dynamique $\alpha$ doit permettre de \og rendre compte\fg\, ou d'\og expliquer \fg la première par la seconde.

Il s'avère que plusieurs définitions non équivalentes, plus ou moins restrictives, peuvent être posées en ce sens. En particulier, on peut définir une interprétation de $\beta$ par $\alpha$ comme un dynamorphisme $\alpha\looparrowright \beta$ vérifiant certaines conditions --- nous parlerons dans ce cas d'interprétation \textit{entrante} ---, ou bien comme un dynamorphisme $\beta\looparrowright \alpha$ vérifiant certaines autres conditions, et nous parlerons alors d'interprétation \textit{sortante}. En appliquant des conditions restrictives, on peut faire en sorte que les deux notions deviennent équivalentes : on parlera dans ce cas d'\textit{interprétation mixte}. D'un autre coté, on peut aussi chercher à englober les deux notions dans une notion plus générale, dans laquelle une interprétation ne serait plus nécessairement un dynamorphisme, et dans ce cas on parlera d'\textit{interprétation généralisée}.

Les définitions qui suivent posent des jalons en ce sens.

\begin{df}
[Interprétation entrante d'une $\mathbf{E}$-dynamique] 
\label{def E interpretation entrante 1}
Soit $\beta:\mathbf{E}\longrightarrow\mathbf{P}$ une $\mathbf{E}$-dynamique.
On appelle $\mathbf{E}$-interprétation entrante de $\beta$ 
tout couple $(\alpha,\psi)$ constitué d'une $\mathbf{E}$-dynamique $\alpha$ 
et d'un dynamorphisme quasi-dé\-ter\-ministe
\footnote{Voir la définition \ref{divers styles de dynamo} page \pageref{divers styles de dynamo}.}
$\psi:\alpha \longrightarrow \beta$ tel que 
pour tout $\mathbf{E}$-écoulement $f:S\rightarrow T$, on ait
\[ 
\psi_T^\mathcal{P}\circ {{^\mu f}^\alpha} \circ \psi_S^{-1}=f^\beta
\]
\noindent autrement dit\footnote{Avec un trans-typage implicite selon que l'on considère $\psi_T$ comme une \emph{fonction} définie sur $S^\alpha$ ou comme une transition.}
\[ 
\psi_T \odot {{f}^\alpha} \odot \psi_S^{-1}=f^\beta
\]
\end{df}

La définition ci-dessus exprime en effet l'idée que $\psi$ permet d'interpréter les états et les transitions \og observés \fg\, de $\beta$ comme correspondant \og en réalité \fg\, aux états et transitions de $\alpha$. On demande que $\psi$ soit quasi-dé\-ter\-mi\-niste car deux états observés distincts ne doivent pas pouvoir s'interpréter par un même état de la dynamique $\alpha$, tandis que certains états de $\alpha$ peuvent bien ne pas être \og observés \fg\, dans la dynamique $\beta$.

Remarquons que, par définition d'un dynamorphisme, on a nécessairement l'inclusion  $\psi_T \odot {{f}^\alpha} \odot \psi_S^{-1}\subset f^\beta$, de sorte qu'il aurait suffit dans la définition précédente de demander l'inclusion réciproque.

\begin{prop}\label{prop CNS interpretation entrante} Étant donné un dynamorphisme quasi-déterministe $\psi:\alpha\rightarrow \beta$, avec $\mathbf{E}_\alpha=\mathbf{E}_\beta=\mathbf{E}$, le couple $(\alpha,\psi)$ est une $\mathbf{E}$-interprétation entrante de $\beta$ si et seulement si les deux conditions suivantes sont satisfaites
\begin{itemize}
\item pour tout $f:S\rightarrow T$ dans $\vec{\mathbf{E}}$, $\psi_T^\mathcal{P}\circ f^\alpha = f^\beta \circ \psi_S$,
\item pour tout $S\in\dot{\mathbf{E}}$, $\psi_S$ est une fonction surjective.
\end{itemize}
\end{prop}

\noindent \textbf{Preuve}. 

1. On suppose que $(\alpha,\psi)$ est une $\mathbf{E}$-interprétation entrante de $\beta$. Pour tout $S$ tel que $S^\beta$ soit non vide, et tout $s\in S^\beta$, on a $\psi_T^\mathcal{P}\circ {{^\mu Id_S}^\alpha} \circ \psi_S^{-1}(s)=Id_S^\beta(s)$, soit $\psi_T^\mathcal{P}(\psi_S^{-1}(s))=\{s\}$, donc $\psi_S^{-1}(s)\neq \emptyset$, d'où la surjectivité de $\psi_S$. 

En outre, pour tout $r\in S^\alpha$ tel que $\psi_S$ soit défini en $r$ (sinon la conclusion est immédiate), et tout $t\in f^\beta(\psi_S(r))$, on a $t\in \psi_T^\mathcal{P}\circ {{^\mu f}^\alpha} \circ \psi_S^{-1}(\psi_S(r))=\psi_T^\mathcal{P}\circ {{^\mu f}^\alpha} (\{r\})=\psi_T^\mathcal{P}\circ {{ f}^\alpha} (r)$.

2. Réciproquement, si $\psi_S$ est une fonction surjective alors tout élément $s\in S^\beta$ s'écrit $\{s\}=\psi_S^\mathcal{P}(\psi_S^{-1}(s))$, de sorte que l'égalité $\psi_T^\mathcal{P}\circ f^\alpha=f^\beta\circ\psi_S$ implique que  $\psi_T^\mathcal{P}\circ {^\mu f}^\alpha (\psi_S^{-1}(s))={^\mu f}^\beta\circ\psi_S^\mathcal{P}(\psi_S^{-1}(s))$, d'où $\psi_T^\mathcal{P}\circ {^\mu f}^\alpha (\psi_S^{-1}(s))={f}^\beta(s)$, et $(\alpha,\psi)$ est bien une interprétation entrante de $\beta$.

\begin{flushright}$\square$\end{flushright} 

\begin{rmq} On trouve facilement des contre-exemples illustrant qu'aucune des deux conditions de la proposition \ref{prop CNS interpretation entrante} n'est à elle seule suffisante à assurer que $\psi_S$ constitue une interprétation entrante.

Par exemple, soit $\mathbf{E}=\{ S \rightarrow T\}$ une catégorie ayant deux objets $S$ et $T$ avec une unique flèche non identité $f:S\rightarrow T$, et soit 
\begin{itemize}
\item $\alpha$ la $\mathbf{E}$-dynamique définie par $S^\alpha=\{s_0\}$, $T^\alpha=\{t_0\}$, $f^\alpha(s_0)=\{t_0\}$,
\item $\beta$ la $\mathbf{E}$-dynamique définie par $S^\beta=\{s_0, s_1\}$, $T^\beta=\{t_0,t_1\}$, $f^\beta(s_0)=f^\beta(s_1)=\{t_0\}$,
\item $\psi:\alpha\looparrowright \beta$ défini par $\psi_S(s_0)=s_0$ et $\psi_T(t_0)=t_0$,
\end{itemize}
\noindent où $s_0$, $t_0$ et $t_1$ désignent les états de ces dynamiques. Alors la première condition de la proposition \ref{prop CNS interpretation entrante} est satisfaite, mais $\psi_S$ et $\psi_T$ ne sont pas surjectives, et l'on n'a donc pas ici une interprétation de $\beta$ puisque tous ses états ne sont pas \og sauvés\fg.

De même, toujours sur  $\mathbf{E}=\{ S \rightarrow T\}$, considérons 
 
\begin{itemize}
\item $\beta$ la $\mathbf{E}$-dynamique définie par $S^\beta=\{s_0\}$, $T^\beta=\{t_0,t_1\}$, $f^\beta(s_0)=\{t_0,t_1\}$,
\item $\alpha$ la $\mathbf{E}$-dynamique définie par $S^\alpha=\{s_0\}$, $T^\alpha=\{t_0,t_1\}$, $f^\alpha(s_0)=\{t_0\}$,
\item $\psi:\alpha\looparrowright \beta$ défini par $\psi_S(s_0)=s_0$ et $\psi_T(t_i)=t_i$ pour $i\in\{0,1\}$.
\end{itemize}

\noindent  Cette fois $\psi_S$ est surjective pour tout $S$, mais cela ne suffit pas car une interprétation entrante doit en quelque sorte rendre compte de \textit{toutes} les transitions \textit{possibles} : sur ce contre-exemple $\alpha$ \og n'explique pas\fg\, comment $\beta$ peut passer de $s_0$ à $t_1$. 
\end{rmq}

\begin{df} Une $\mathbf{E}$-\emph{interprétation sortante} de $\beta$ est un couple $(\varphi,\alpha)$ où
 $\alpha$ est une $\mathbf{E}$-dynamique et $\varphi$ est un $\mathbf{E}$-dynamorphisme 
 $\varphi:\beta\looparrowright \alpha$ vérifiant les conditions suivantes pour tout $S\in\dot{\mathbf{E}}$ et tout $(f:S\rightarrow T)\in\vec{\mathbf{E}}$ :
\begin{itemize}
\item $\forall s\in S^\beta, \varphi_S(s)\neq\emptyset$, 
\item $\forall s_0, s_1\in S^\beta, s_0\neq s_1 \Rightarrow\varphi_S(s_0)\cap\varphi_S(s_1)=\emptyset$.
\end{itemize}

Si en outre $\varphi$ vérifie ${f}^\alpha\odot \varphi_S={{\varphi}_T}\odot f^\beta$, l'interprétation sortante $(\varphi,\alpha)$ est dite \emph{régulière}.

\end{df}

Pour toute $\mathbf{E}$-interprétation entrante $(\alpha,\psi)$ de $\beta$, on note $\tilde{\psi}$ la famille de transitions $(\tilde{\psi}_S:S^\beta\rightsquigarrow S^\alpha)_{S\in\dot{\mathbf{E}}}$ définie par
\[\tilde{\psi}_S=\psi_S^{-1}.\] 

Pour toute interprétation sortante $(\varphi,\alpha)$ de $\beta$, on note $\tilde{\varphi}$ la famille de transitions $(\tilde{\varphi}_S:S^\alpha\rightsquigarrow S^\beta)_{S\in\dot{\mathbf{E}}}$ définie par
\begin{itemize}
\item $\tilde{\varphi}_S(r)=\emptyset$ si $r\notin {^\mu \varphi}_S(S^\beta)$
\item $\tilde{\varphi}_S(r)=\{s\}$ si $r\in {\varphi}_S(s)$ (il existe au plus un tel $s$),
\end{itemize}

\begin{df} Soient $(\alpha,\psi)$ et $(\varphi,\alpha)$ deux $\mathbf{E}$-interprétations, respectivement entrante et sortante, de $\beta$. Si, pour tout $S\in\dot{\mathbf{E}}$, $\psi=\tilde{\varphi}$ ou, de façon équivalente, $\varphi=\tilde{\psi}$, on dit que ces deux interprétations sont \emph{associées}. Lorsqu'une interprétation admet une interprétation associée, elle est dite \emph{mixte}. Une interprétation mixte est dite régulière si elle est sortante régulière ou associée à une interprétation sortante régulière.
\end{df}

\begin{prop} Soit $(\alpha,\psi)$ une $\mathbf{E}$-interprétation entrante de $\beta$. Si $\varphi=\tilde{\psi}$  définit un dynamorphisme $\varphi:\beta\looparrowright\alpha$, alors $(\alpha,\psi)$ est mixte, d'interprétation sortante associée $(\varphi,\alpha)$.

Soit $(\varphi,\alpha)$ une $\mathbf{E}$-interprétation sortante de $\beta$. 
Si $\psi=\tilde{\varphi}$  définit un dynamorphisme $\psi:\alpha\looparrowright\beta$, alors $(\varphi,\alpha)$ est mixte, d'interprétation entrante associée $(\alpha,\psi)$.
\end{prop}

\noindent \textbf{Preuve}. La première assertion découle simplement de ce que $\psi$ est quasi-déterministe, de sorte que $\varphi_S(s_1)\cap\varphi_S(s_2)=\emptyset$ pour $s_1\neq s_2$ deux éléments de $S^\beta$ avec $S\in\dot{\mathbf{E}}$ quelconque, et qu'en outre $\psi_S$ est surjective, de sorte que $\varphi_S(s)\neq \emptyset$ pour tout $s\in S^\beta$.

Vérifions la seconde assertion. 
Pour tout $\mathbf{E}$-écoulement $f:S\rightarrow T$, $\varphi$ étant un dynamorphisme,
on a $\varphi_T\odot f^\beta\subset f^\alpha\odot \varphi_S$, 
d'où $\psi_T\odot \varphi_T\odot f^\beta\subset\psi_T\odot f^\alpha\odot \varphi_S$. 
Or, $\psi_T\odot\varphi_T=Id_{T^\beta}$. 
D'où $f^\beta\subset\psi_T\odot f^\alpha\odot \varphi_S$.
D'un autre coté, $\psi$ étant également un dynamorphisme, on a
$\psi_T\odot f^\alpha\subset f^\beta\odot\psi_S$, d'où
$\psi_T\odot f^\alpha\odot \varphi_S\subset f^\beta\odot\psi_S\odot \varphi_S$. En utilisant l'égalité $\psi_S\odot\varphi_S=Id_{S^\beta}$, on en déduit
$\psi_T\odot f^\alpha\odot \varphi_S\subset f^\beta$, d'où finalement
$\psi_T\odot f^\alpha\odot \varphi_S = f^\beta$ ce qui, d'après la définition \ref{def E interpretation entrante 1}, devait être démontré.
\begin{flushright}$\square$\end{flushright} 

Dans les cinq exemples suivants (\ref{exm inter mixte non reg} à \ref{exm inter gen}), on prend pour $\mathbf{E}$ la catégorie ayant deux objets $S$ et $T$ avec une unique flèche non identité : $f:S\rightarrow T$, et pour dynamique \og à interpréter \fg la dynamique $\beta$ définie par $S^\beta=\{p\}$, $T^\beta=\{q\}$ et $f^\beta(p)=\{q\}$.

\begin{exm}[Une interprétation mixte non régulière]\label{exm inter mixte non reg}
On définit $\alpha$ par
$S^\alpha=\{x\}$, $T^\alpha=\{y,y'\}$ et $f^\alpha(x)=\{y,y'\}$.
On définit $\psi$ par $\psi_S(x)=\{p\}$, $\psi_T(y)=\{q\}$ et $\psi_T(y')=\emptyset$.

L'interprétation entrante $(\alpha,\psi)$ est mixte, mais n'est pas régulière à cause du $y'$. On remarque que dans la dynamique $\alpha$, l'écoulement $f$ peut conduire de l'état $x$ à l'état $y'$ qui reste insu de $\beta$ : cela n'est pas en contradiction avec le fonctionnement de $\beta$, car on rappelle que $f^\beta(p)=\{q\}$ ne signifie pas que le passage de $p$ à $q$ soit obligatoire, certaines solutions étant vides en  en $T$.
\end{exm}

\begin{exm}[Une interprétation entrante (non mixte)]\label{exm inter entr}
On définit $\alpha$ par
$S^\alpha=\{x\}$, $T^\alpha=\{y,y',y''\}$ et $f^\alpha(x)=\{y',y''\}$.
On définit $\psi$ par $\psi_S(x)=\{p\}$, $\psi_T(y)=\psi_T(y')=\{q\}$ et $\psi_T(y'')=\emptyset$.

On remarque que dans cette interprétation, $q$ n'est pas nécessairement interprété comme \og conséquence \fg\, de $p$, car il peut aussi se faire qu'un état $y$ ait surgi on ne sait d'où.
\end{exm}

\begin{exm}[Une interprétation sortante (non mixte) régulière ]\label{exm inter sort reg}
On pose
$S^\alpha=\{x,x'\}$, $T^\alpha=\{y,y',y''\}$, $f^\alpha(x)=\{y,y'\}$ et $f^\alpha (x')=y'$. et comme interprétation sortante on prend $\varphi$ avec $\varphi_S(p)=\{x\}$, $\varphi_T(q)=\{y,y'\}$.

L'interprétation n'est pas mixte à cause de l'état $x'$, \og insu\fg\, de $\beta$, mais qui évolue vers un état $y'$ correspondant à un état de $\beta$.
\end{exm}

\begin{exm}[Une interprétation sortante (non mixte) non régulière ]\label{exm inter sort non reg}
On pose comme dans l'exemple précédent
$S^\alpha=\{x,x'\}$, $T^\alpha=\{y,y',y''\}$, $f^\alpha(x)=\{y,y'\}$ et $f^\alpha (x')=y'$. 

L'interprétation sortante est défini par $\varphi_S(p)=\{x\}$, et $\varphi_T(q)=\{y'\}$.

Comme dans l'exemple précédent, l'interprétation n'est pas mixte à cause du $x'$ qui \og entre dans le jeu\fg\. En outre, elle est non régulière à cause du $y$ qui, lui, en sort.
\end{exm}

\begin{exm}[Une interprétation généralisée]\label{exm inter gen}
Si, dans les données de l'exemple précédent, on remplace uniquement $\varphi_T(q)=\{y'\}$ par $\varphi_T(q)=\{y',y''\}$, on obtient une famille $\varphi=(\varphi_S,\varphi_T$ de transitions qui, avec $\alpha$, ne constitue une interprétation au sens d'aucune des définitions précédentes. Néanmoins, on peut décrire la configuration obtenue en assemblant des morceaux des situations précédentes, et l'on pourrait par conséquent songer à élargir la notion d'interprétation pour y inclure ce type d'exemples. Nous ne le ferons pas ici.
\end{exm}

\begin{exc} On dit qu'une interprétation entrante $(\alpha, \psi)$ ou une interprétation sortante $(\varphi, \alpha)$, est déterministe (resp. non déterministe, etc.) si $\alpha$ l'est. Montrer que toute $\mathbf{E}$-dynamique déterministe admet une $\mathbf{E}$-interprétation mixte régulière non déterministe, et que toute  $\mathbf{E}$-dynamique non déterministe admet une $\mathbf{E}$-interprétation mixte régulière déterministe.
\end{exc}

\subsection{Interprétations trans-temporelles}

On cherche à présent à définir l'interprétation d'une dynamique par une autre n'ayant pas nécessairement même type temporel. Par exemple, on veut pouvoir interpréter une dynamique discrète, définie sur $(\mathbf{N},+)$ par une dynamique continue\footnote{Par contre, on ne cherchera pas à interpréter une dynamique continue par une dynamique discrète : au sens où nous l'entendons, une interprétation n'est pas une approximation, elle doit rendre compte de façon complète de la dynamique interprétée.}, définie sur $(\mathbf{R}_+,+)$. Un dynamorphisme entre deux dynamiques admet une partie fonctorielle qui, si on voulait définir les interprétations par des dynamorphismes \og entrant \fg, autrement dit avec pour but la dynamique à interpréter, serait sur l'exemple considéré introuvable, puisqu'il n'y a qu'un seul morphisme de monoïde $(\mathbf{R}_+,+)\rightarrow(\mathbf{N},+)$, autrement dit qu'un seul foncteur entre les catégories correspondantes.

C'est pourquoi nous posons la définition suivante en généralisant la notion de $\mathbf{E}$-interprétation sortante d'une $\mathbf{E}$-dynamique $\beta$.

\begin{df} Étant données deux dynamiques $\beta: \mathbf{E}\rightarrow \mathbf{P}$ et $\alpha:\mathbf{F}\rightarrow \mathbf{P}$, et $(\Phi,\varphi):\beta \looparrowright \alpha$ un dynamorphisme. On dit que $((\Phi,\varphi),\alpha)$ constitue une interprétation (sortante) de $\beta$ si
\begin{itemize}
\item $\Phi$ est un foncteur fidèle\footnote{C'est-à-dire injectif sur les flèches, donc sur les objets.},
\item pour tout $S\in\dot{\mathbf{E}}$ et tout $(f:S\rightarrow T)\in\vec{\mathbf{E}}$ :
\begin{itemize}
\item $\forall s\in S^\beta, \varphi_S(s)\neq\emptyset$, 
\item $\forall s_0, s_1\in S^\beta, s_0\neq s_1 \Rightarrow\varphi_S(s_0)\cap\varphi_S(s_1)=\emptyset$.
\end{itemize}
\end{itemize}
On dit que l'interprétation est régulière si pour tout $\mathbf{E}$-écoulement $f:S\rightarrow T$, on a $\varphi_T\odot f^\beta ={(\Phi f)}^\alpha \odot \varphi_S$.
\end{df}

\begin{exm}[Interprétation continue d'une dynamique discrète]
En voici un exemple élémentaire. On prend $\mathbf{E}=(\mathbf{N},+)$, $st(\beta)=\mathbf{R}$, $1^\beta(x)=x+1$, $\mathbf{F}=(\mathbf{R},+)$, $st(\alpha)=\mathbf{R}$ et $r^\alpha (x)=x+\alpha$. On obtient ainsi une interprétation régulière de la dynamique discrète $\beta$ par la dynamique continue $\alpha$. 
\end{exm}

%

\chapter{Dynamiques catégoriques connectives}

Le but de ce chapitre est essentiellement de \emph{définir} les dynamiques caté\-go\-ri\-ques connectives, à savoir les dynamiques catégoriques dont le type de temporalité et l'ensemble des états sont chacun munis d'une structure connective. \`{A} cette fin, on commence par définir les catégories connectives. La définition des dynamiques catégoriques connectives est alors donnée, ainsi que celle du feuilletage associé, d'où découle immédiatement la notion d'ordre connectif d'une dynamique catégorique connective.


%
%
%

\section{Catégories connectives}

\subsection{Définition}

Pour tout couple $(A,B)$ d'ensembles de flèches d'une catégorie, on note $B\circ A$ l'ensemble des flèches de la forme $b\circ a$ avec $b\in B$ et $a\in A$ :

\begin{equation}
B\circ A = \{b\circ a, a\in A, b\in B, \mathrm{dom}(b)=\mathrm{cod}(a)\}.
\end{equation}

\begin{df}[Catégories connectives]
Une catégorie connective $\mathbf{C}$ est une petite catégorie $\vert \mathbf{C} \vert$ dont l'ensemble des flèches $\vec{\mathbf{C}}$ est muni d'une structure connective $\kappa_{\mathbf{C}}$ telle que, pour tout couple d'ensembles connexes de flèches $(A,B)$, l'ensemble de flèches $B\circ A$ soit connexe:
\begin{equation}\label{structures connectives categorique}
\forall (A,B)\in {\kappa_{\mathbf{C}}}^2, B\circ A \in {\kappa_{\mathbf{C}}}
\end{equation}
\end{df}

Étant donnée une catégorie connective $\mathbf{C}$, on notera, conformément à la définition précédente, 
\begin{itemize}
\item $\mathbf{C}_0$ ou $\dot{\mathbf{C}}$ l'ensemble de ses objets,
\item $\mathbf{C}_1$ ou $\vec{\mathbf{C}}$ l'ensemble de ses flèches,
\item $\vert \mathbf{C} \vert=(\mathbf{C}_0, \mathbf{C}_1, \circ, \mathrm{dom}, \mathrm{cod}, id)$ la catégorie sous-jacente,
\item $\kappa_{\mathbf{C}}$ sa structure connective, définie sur $\mathbf{C}_1$.
\end{itemize}

Pour distinguer les structures connectives sur l'\emph{ensemble} $\mathbf{C}_1$ des structures connectives sur la \emph{catégorie} $\vert \mathbf{C} \vert$, on dira parfois de ces dernières qu'elles sont des \emph{structures connectives catégoriques}. Une structure connective sur l'ensemble $\vec{\mathbf{C}}$ est donc catégorique si elle vérifie en outre la relation {(\ref{structures connectives categorique})}.

\begin{df}\label{def cat des categories connectives} La \emph{catégorie des catégories connectives} est définie en prenant pour morphismes entre catégories connectives les \emph{foncteurs connectifs}, c'est-à-dire les foncteurs qui transforment tout ensemble connexe de flèches en ensemble connexe de flèches.
\end{df}

\subsection{Treillis de structures connectives catégoriques}

\begin{prop}[Treillis de structures] L'ensemble des structures connectives dont peut être munie une petite catégorie donnée constitue un treillis complet pour l'inclusion. La structure minimale est la structure \emph{totalement discrète}, qui n'admet pour connexe que la partie vide. La structure maximale est la structure grossière ou indiscrète, pour laquelle tout ensemble de flèches est connexe. La borne inférieure d'une famille de structures connectives sur la catégorie considérée est leur intersection. 

Les structures connectives intègres vérifient les mêmes propriétés, à ceci près que la structure minimale est dans ce cas la structure \emph{discrète intègre}, qui n'admet pour connexes non vides que les singletons.
\end{prop}

\noindent \textbf{Preuve}. La vérification est immédiate, en particulier le fait que l'intersection d'une famille non vide quelconque de structures connectives catégoriques sur une catégorie est encore une structure connective catégorique sur cette catégorie, et de même pour les structures intègres.
\begin{flushright}$\square$\end{flushright}

\begin{cor} Tout ensemble $\mathcal{F}\subset\mathcal{P}(\vec{\mathbf{E}})$ 
d'ensembles de flèches est contenu dans une structure connective catégorique minimale, notée\label{notation double crochet structure} $\llbracket \mathcal{F} \rrbracket_0$. 
De même, il existe une structure connective catégorique intègre minimale contenant $\mathcal{F}$, notée\label{notation double crochet structure integre} $\llbracket \mathcal{F} \rrbracket$, 
et l'on a
\[ 
\llbracket \mathcal{F} \rrbracket = \llbracket \mathcal{F}\cup \{\{f\}, f\in \vec{\mathbf{E}}\} \rrbracket_0 .
\]
\end{cor}

\begin{exc} Donner une construction de ${\llbracket \mathcal{F} \rrbracket}_0$ et de $\llbracket \mathcal{F} \rrbracket$ à partir de $\mathcal{F}\subset\mathcal{P}(\vec{\mathbf{E}})$. Dans ce but, on pourra noter\label{notation ens compositions fleches} $^\circ\mathcal{F}$ l'ensemble des ensembles de flèches défini par
\begin{equation}
^\circ\mathcal{F}=\bigcup_{n\geqslant 1}\{A_1\circ ... \circ A_n, (A_1,...,A_n)\in \mathcal{F}^n\}.
\end{equation}
L'obtention d'une telle construction permettra de désigner ${\llbracket \mathcal{F} \rrbracket}_0$ comme la structure connective catégorique \emph{engendrée} par $\mathcal{F}$
\end{exc}

%
%
%
%

\subsection{Ordre brunnien d'une petite catégorie}

Étant donnée $\mathbf{E}$ une petite catégorie. On lui associe  une catégorie connective  $\mathbf{E}_\mathcal{B}$, dite brunnienne intègre, en prenant pour structure connective sur $\mathbf{E}_1$ la plus fine des structures catégoriques intègres qui rende $\mathbf{E}_1$ connexe :
\begin{equation}
\kappa(\mathbf{E}_\mathcal{B})=\llbracket \{\mathbf{E}_1\} \rrbracket.
\end{equation}

\begin{df}[Ordre brunnien d'une petite catégorie] On appelle \emph{ordre brunnien} d'une petite catégorie $\mathbf{E}$ l'ordre connectif de la catégorie brunnienne intègre associée, c'est-à-dire l'ordre connectif de l'espace connectif $(\mathbf{E}_1, \kappa(\mathbf{E}_\mathcal{B}))$.
\end{df}

On remarquera que l'ordre brunnien d'un groupe (vu comme petite catégorie) est $1$, que celui de $(\mathbf{N},+)$ est $\aleph_0$, que celui de $(\mathbf{R_+},+)$ est $\aleph_1$.

\begin{exm} [Composantes connexes d'une catégorie] On munit ordinairement la classe des objets d'une catégorie $\mathbf{E}$ d'une structure connective (en fait une relation d'équivalence), de la manière suivante : deux objets sont équivalents (ou connectés) si et seulement si il existe un morphisme de l'un vers l'autre. La façon la plus économique de retrouver cette structure en terme de (petites) catégories connectives est sans doute de poser
\[ 
\kappa(\mathbf{E})
={[\{\{id_A,id_B\}, 
(A,B)\in{\dot{\mathbf{E}}}^2, 
\exists f\in\vec{\mathbf{E}}, 
\mathrm{dom}(f)=A, \mathrm{cod}(f)=B\}]}_0
\]
\noindent un ensemble d'objets de $\mathbf{E}$ étant alors considéré comme connexe si l'ensemble de leurs identités l'est.
\end{exm}
 
\subsection{Monoïdes connectifs}

Le foncteur $MC$ (voir l'exemple \ref{foncteur MC} page \pageref{foncteur MC}) permet de définir les monoïdes connectifs : ce sont les monoïdes qui, vus comme petites catégories, sont munis d'une structure connective. La définition suivante explicite et développe cette idée.

\begin{df} [monoïdes connectifs, (semi-)groupes (semi-)connectifs] 
Un \emph{monoïde connectif} est un quadruplet $(M,*,\varepsilon,\mathcal{M})$ tel que 
\begin{itemize}
\item $(M,*)$ est un monoïde d'élé\-ment neutre $\varepsilon$,
\item $(M,\mathcal{M})$ est un espace connectif,
\item $\forall (A,B)\in \mathcal{M}^2, A*B\in \mathcal{M}$,
 où  $A*B=\{a*b, (a,b)\in {A}\times {B}\}$. 
\end{itemize}

Un \emph{semi-groupe connectif} est un monoïde connectif qui est un semi-groupe, c'est-à-dire dont la loi est régulière : $\forall a, b, c \in M, [(a*b=a*c) \, \mathrm{ou} \, (b*a=c*a)]\Rightarrow b=c$.

Un \emph{groupe semi-connectif} est un monoïde connectif $(G,*,\mathcal{G})$ tel que $(G,*)$ soit un groupe.

Un \emph{groupe connectif} est un groupe semi-connectif tel qu'en outre l'inversion transforme toute partie connexe en partie connexe.

\end{df}

\begin{rmq}\label{remarque tensoriel connectif} Contrairement à l'usage qui prévaut par exemple dans la dé\-fi\-ni\-tion des groupes topologiques, la définition des monoïdes connectifs n'implique pas que la loi de composition $M\times M\rightarrow M$ soit connective. Par exemple, muni de l'addition et de sa structure connective usuelle, l'ensemble $\mathbf{N}$ des entiers naturels constitue un monoïde connectif, bien que l'addition ne soit pas connective sur le carré cartésien de cet espace, puisque par exemple l'ensemble de couples $\{(0,0),(1,1)\}$ est une partie connexe de ce carré cartésien (ses deux projections sont connexes dans $\mathbf{N}$) mais son image par l'addition est l'ensemble non connexe $\{0,2\}$. Par contre, l'addition est connective sur le carré \emph{tensoriel} connectif de $\mathbf{N}$, dont $\{(0,0),(1,1)\}$ n'est pas une partie connexe. Cette notion de produit tensoriel des espaces connectifs est développée dans la section 4 de notre article \cite{Dugowson:201012}. On retiendra en particulier qu'une application de plusieurs variables est partiellement connective par rapport à chacune de ces variables si et seulement si elle est connective sur le produit tensoriel des espaces connectifs décrits par ces variables, et que muni du produit tensoriel la catégorie des espaces connectifs intègres est une catégorie monoïdale fermée.
\end{rmq}
\begin{prop} Muni d'une structure connective \emph{intègre}, un monoïde $(M,.)$ est un monoïde connectif si et seulement si pour tout $x\in M$, la translation à gauche $y\mapsto x.y$ et la translation à droite $y\mapsto y.x$ sont des applications connectives de $M$ dans lui-même, autrement dit si pour tout $x\in M$ et toute partie connexe $K$ de $M$, on a $x.K$ et $K.x$ connexes.
\end{prop}
\noindent \textbf{Preuve}. L'espace étant supposé intègre, la condition est évidemment né\-cessai\-re. Montrons qu'elle est suffisante. Soient donc $F$ et $G$ deux parties connexes de $M$, que l'on peut supposer non vides, la connexité de $FG$ étant triviale si l'une de ces parties est vide. Soit $g_0\in G$, on a 
\[ 
F.G
= \bigcup_{f\in F}{f.G}
= \bigcup_{f\in F}{(f.G \cup {Fg_0})}
 \]

 Les deux parties $f.G$ et $F.g_0$ sont connexes et d'intersection non vide, donc leur union est connexe. Les parties  $f.G \cup {F.g_0}$ sont donc connexes, mais comme elles sont elles aussi d'intersection non vide, leur union est encore {connexe.
 
\begin{flushright}$\square$\end{flushright} }


\section{Dynamiques catégoriques connectives}

\begin{df}\label{definition dynamique categorique connective} Étant donnée un système catégorique d'écou\-lements $\mathbf{E}$, une $\mathbf{E}$-dynamique connective $(\alpha, \kappa_0, \kappa_1)$, simplement notée $\alpha$ s'il n'y a pas d'ambiguïté, est la donnée
\begin{itemize}
\item d'une $\mathbf{E}$-dynamique catégorique $\alpha:\mathbf{E}\rightarrow \mathbf{P}$,
\item d'une structure connective $\kappa_0$ sur la catégorie $\mathbf{E}$,
\item d'une structure connective $\kappa_1$ sur l'ensemble $St_\alpha=\bigcup_{T\in \mathbf{E}_0}{T^\alpha}$ des états de $\alpha$.
\end{itemize}
\end{df}

A partir de cette définition, on pourra définir et étudier différentes catégories de dynamiques connectives. On pourra en particulier étudier ce qui, dans les notions relatives aux dynamiques catégoriques ensemblistes considérées au chapitre 3, peut être transposé de façon naturelle dans le contexte connectif où nous nous plaçons ici. 

Nous ne procéderons pas dans le présent ouvrage à ces explorations, notre objectif ayant été ici essentiellement de poser la \textit{définition} des dynamiques catégoriques connectives. On se contentera pour le moment de définir l'ordre connectif d'une dynamique connective.

\begin{df} Étant donnée $(\alpha, \kappa_0,\kappa_1)$ une dynamique catégorique connective définie sur une petite catégorie $E$, le feuilletage connectif $\mathcal{F}_\alpha$ associé est défini par
\begin{itemize}
\item ses points, qui sont les états de la dynamique : $\vert \mathcal{F}_\alpha\vert=St_\alpha$,
\item sa structure interne, qui est la structure connective engendrée par les 
$K^\alpha (a) = \{b\in St_\alpha, \exists f\in K, b\in f^\alpha (a)\}$ avec $K\in \kappa_0$  et $a\in St_\alpha$,
\item sa structure externe, qui est $\kappa_1$.
\end{itemize}
\end{df}

\begin{rmq}\label{dynamique avec feuilletage interne non integre} La structure interne du feuilletage associé à une dynamique catégorique connective n'est pas nécessairement intègre. En effet, si pour un état $a\in T^\alpha$ le seul écoulement temporel $f$ et le seul état $b$ tels que $f^\alpha(b)=\{a\}$ sont $b=a$ et $f=Id_T$, alors il suffit que le singleton $\{Id_T\}$ ne soit pas dans la structure temporelle $\kappa_0$ pour que $\{a\}$ n'appartienne pas à la structure interne du feuilletage. C'est par exemple ce qui se produit si l'on muni le monoïde $(\mathbf{R}_+,+)$ de la structure connective telle que les connexes sont les intervalles \emph{qui ne contiennent pas} $0$, et que l'on considère une dynamique $\alpha$ telle que tous les états $a$ situés en dehors d'une famille donnée d'orbites vérifient $f^\alpha(a)=\emptyset$ dès que $f>0$.
\end{rmq}

\begin{df} L'\emph{ordre connectif} d'une dynamique connective $\alpha$ est l'ordre connectif de son feuilletage $\mathcal{F}_\alpha$, c'est-à-dire l'ordre connectif de l'espace connectif induit des feuilles $\mathcal{F_\alpha}^\downarrow$.

\end{df}


\noindent \textbf{Exemple 1}. Une rotation irrationnelle sur un cercle définit une dynamique discrète. Le monoïde $(\mathbf{N},+)$ étant muni d'une structure connective qui en fasse un espace connexe (on peut par exemple prendre la structure engendrée par les $\{n,n+1\}$, ou bien la structure brunnienne intègre, etc.), et le cercle de la structure connective associée à la topologie usuelle, la dynamique en question admet pour espace de feuilles induit un ensemble qui a la puissance du continu et dont la structure connective est la structure brunnienne non intègre.

\noindent \textbf{Exemple 2}. Le pendule sphérique linéaire détermine, sur un sous-espace d'iso-énergie $S^3$ de son espace de phase, une dynamique connective d'ordre 1 ($S^3$ étant munie de la structure de séparation standard). En effet, les feuilles de cette dynamiques sont les cercles constituant une fibration de Hopf, et l'on sait que ces cercles sont deux à deux entrelacés. Même chose pour le papillon de Lorenz, dont les orbites périodiques sont deux à deux entrelacées (et dont les orbites non périodiques ne sont pas davantage séparables par un plan topologique).

\noindent \textbf{Exemple 3}. En 1996, dans \cite{Ghrist:1995, GH:1996, GHS:1997}, Ghrist et Holmes présentent une équation différentielle dont l'ordre connectif  est  au moins $\aleph_0$, puisqu'on peut former tout entrelacs avec ses trajectoires périodiques, de sorte que, d'après le théorème de Brunn-Debrunner-Kanenobu, toute structure connective finie se trouve représentée dans le flot.

\chapter*{Notations}


\paragraph{Catégories, flèches et objets}

Pour tout objet $S$ d'une catégorie, la flèche identité de $S$ est notée $Id_S$. En particulier, pour toute catégorie $\mathbf{E}$, l'endofoncteur identité de $\mathbf{E}$ est noté $Id_\mathbf{E}$, tandis que $\mathbf{1}_\mathbf{E}$ désigne, s'il existe, l'objet terminal de $\mathbf{E}$. En outre, 

\begin{notations}\label{notations fleches et objets} 
\mbox{}
\begin{itemize}
\item la classe des objets de $\mathbf{E}$ est notée $\dot{\mathbf{E}}$ ou bien $\vert\mathbf{E}\vert_0$,
\item la classe des flèches de $\mathbf{E}$ est notée $\vec{\mathbf{E}}$ ou $\vert\mathbf{E}\vert_1$.
\item Si $F$ désigne un foncteur, on notera parfois $\dot{F}$ la partie objet de $F$, et $\vec{F}$ son action sur les flèches.
\item La source de $f\in \vec{\mathbf{E}}$ est $\mathrm{dom}(f)$, son but est $\mathrm{cod}(f)$,
\item On note parfois $fg$ la composée $g\circ f$ de deux flèches (voir à ce propos la notation \ref{notations fg} page \pageref{notations fg}),
\item Si $T$ est un objet de $\mathbf{E}$, on note $\{\rightarrow T\}$ la classe des flèches de  $\mathbf{E}$ de but $T$ : $\{\rightarrow T\}=\{f\in \vec{\mathbf{E}},\mathrm{cod}(f)=T\}$.
\end{itemize}
\end{notations}

Outre les catégories usuelles, telle que la catégorie des petites catégories $\mathbf{Cat}$, 
les principales catégories spécifiques considérées dans cet ouvrage sont :

\begin{itemize}
\item $\mathbf{Cnc}$ et $\mathbf{Cnct}$ : notation \ref{notation categories espaces connectifs} page \pageref{notation categories espaces connectifs},
\item $ \mathbf{RC}$ et $ \mathbf{RCD}$ : définition \ref{categorie des representations et rcd} page \pageref{categorie des representations et rcd},
\item $\mathbf{FC}$ : définition \ref{categorie FC} page \pageref{categorie FC},
\item $\mathbf{P}$ : définition \ref{def categorie transitions} page \pageref{def categorie transitions},
\item $\mathbf{E-Dy}$ : proposition \ref{categorie E-Dy} page \pageref{categorie E-Dy},
\item  $\mathbf{Dy}$ : définition \ref{categorie Dy} page \pageref{categorie Dy},
\item $\mathbf{DyP}$, $\mathbf{DyPd}$ : page \pageref{DyP}, section \ref{DyP}.
\end{itemize}

Certaines flèches sont notées de façon spécifique. Ainsi les \textit{transitions} (voir la définition \ref{def categorie transitions} page \pageref{def categorie transitions}) sont-elles notées $\rightsquigarrow$, tandis que les dynamorphismes (voir la définition \ref{def E dynamo} page \pageref{def E dynamo} et la définition \ref{categorie Dy} page \pageref{categorie Dy}) sont notés $\looparrowright$.

Nous notons également de façon spécifique la composition des représentations et des transitions, à l'aide du symbole $\odot$ (voir la section \ref{composition des representations} page \pageref{composition des representations} pour les représentations, et la définition \ref{def categorie transitions} page \pageref{def categorie transitions} pour les transitions).

\mbox{}

\noindent \textbf{Remarque}. La catégorie des dynamiques catégoriques connectives, $\mathbf{CaCoDy}$,  n'est pas définie dans cet ouvrage, où l'on se contente d'en définir les objets.

\paragraph{Ensemble de parties}

Outre la notation standard $\mathcal{P}$ pour désigner l'ensemble des parties, et les notations $f^\mathcal{P}$, $^\mu f$ définies dans la notation \ref{notations transitions} page \pageref{notations transitions}, on utilise également les notations suivantes, lorsque $\mathcal{K}X$ désigne un ensemble de parties de $X$ :

\begin{notations}[parties non vides, parties non triviales]

On pose :

\begin{itemize}
\item $\mathcal{K}^* X$ :  ensemble des élements de $\mathcal{K}X$ qui sont des parties non vides de $X$,
\item $\mathcal{K}^\bullet X$ : ensemble des élements de $\mathcal{K}X$ qui sont des parties de cardinal $\geqslant 2$ de $X$,
\end{itemize}

En particulier,

\begin{itemize}
\item $\mathcal{P}^* X$ désigne l'ensemble des parties non vides de $X$,
\item $\mathcal{P}^\bullet X$ désigne l'ensemble des parties de $X$ de cardinal $\geqslant 2$.
\end{itemize}

\end{notations}

Les notation $[\mathcal{A}]_0$ et $[\mathcal{A}]$ sont rappelées page \pageref{notation structure connective engendree} (section \ref{notation structure connective engendree}).

De même, les notations $\llbracket \mathcal{F} \rrbracket_0$ et $\llbracket \mathcal{F} \rrbracket$ sont introduites dans le corolaire \ref{notation double crochet structure}  page \pageref{notation double crochet structure}.

\paragraph{Autres notations}

Les autres notations sont introduites dans le cours du texte, par exemple :

\begin{itemize}
\item notation  $\hat{A}$ dans la formule \ref{notation A chapeau} page \pageref{notation A chapeau},
\item notation $(r\leqslant s)$ dans l'exemple \ref{notations r inf s} page \pageref{notations r inf s}, 
\item notation $\mathbf{E}[(M,*); M_1, N_1, M_2,...M_n,...]$ dans l'exemple \ref{notation exm categorie ecoulements par monoides} page \pageref{notation exm categorie ecoulements par monoides},
\item notation \ref{notation fleches transitions} page \pageref{notations transitions}, notation \ref{notations transitions deterministes} page \pageref{notations transitions deterministes},  \ref{notation image par une dynamique} page \pageref{notation image par une dynamique}, etc...
\end{itemize}

%

%
%

%
%
\bibliographystyle{plain}

\begin{thebibliography}{}

\end{thebibliography}


\begin{thebibliography}{10}

\bibitem{joycats:2005}
Jiri Adamek, Horst Herrlich, and George Strecker.
\newblock {\em Abstract and {C}oncrete {C}ategories: {T}he {J}oy of {C}ats}.
\newblock John Wiley and Sons, 1990.
\newblock On-line edition, 18th January 2005 :
  http://katmat.math.uni-bremen.de/acc.

\bibitem{Ben-Ari:1983}
M.~Ben-Ari, A.~Pnueli, and Z.~Manna.
\newblock The temporal logic of branching time.
\newblock {\em Acta informatica}, 20(3):207--226, 1983.

\bibitem{Borger:1981}
Reinhard B\"orger.
\newblock {\em Kategorielle Beschreibungen von Zusammenhangsbegriffen}.
\newblock PhD thesis, Fernuniversit\"at, Hagen, 1981.

\bibitem{Borger:1983}
Reinhard B\"orger.
\newblock Connectivity spaces and component categories.
\newblock In {\em Categorical topology, International Conference on Categorical
  Topology (1983)}, Berlin, 1984. Heldermann.

\bibitem{Brunn:1892a}
Hermann Brunn.
\newblock Ueber verkettung.
\newblock {\em Sitzungsberichte der Bayerische Akad. Wiss., MathPhys. Klasse},
  22:77--99, 1892.

\bibitem{Brunn:1982b}
Hermann Brunn.
\newblock De l'encha\^\i nement (\emph{Errata}).
\newblock {\em Ornicar ?}, 26-27:289--292, 1982.

\bibitem{Brunn:1982a}
Hermann Brunn.
\newblock De l'encha\^\i nement (french translation by {C}. {L}\'eger and {M}.
  {T}urnheim).
\newblock {\em Ornicar ?}, 25:207--223, 1982.

\bibitem{Bushaw:1963}
D.~Bushaw.
\newblock Dynamical polysystems and optimization.
\newblock {\em Contributions to differential equations}, 2:351, 1963.

\bibitem{Debrunner:19600416}
Hans Debrunner.
\newblock {Links of Brunnian type}.
\newblock {\em Duke Math. J.}, 28:17--23, 1961.

\bibitem{Debrunner:1964}
Hans Debrunner.
\newblock {\"Uber den Zerfall von Verkettungen.}
\newblock {\em Mathematische Zeitschrift}, 85:154--168, 1964.
\newblock {http://www.digizeitschriften.de}.

\bibitem{Dugowson:2003}
Stéphane Dugowson.
\newblock Espaces connectifs et espaces de partage, 2003.
\newblock Unpublished.

\bibitem{Dugowson:2007c}
Stéphane Dugowson.
\newblock Les fronti\`eres dialectiques.
\newblock {\em Mathematics and Social Sciences}, 177:87--152, 2007.

\bibitem{Dugowson:201012}
Stéphane Dugowson.
\newblock On connectivity spaces.
\newblock {\em Cahiers de {T}opologie et {G}éométrie {D}ifférentielle
  {C}atégoriques}, LI(4):282--315, 2010.
\newblock {http://hal.archives-ouvertes.fr/hal-00446998/fr/}.

\bibitem{MmeEhresmann:1965}
Mme A.~Bastiani Ehresmann.
\newblock Sur le problème général d'optimisation.
\newblock In {\em Indentification, Optimisation et Stabilité (actes du congrès
  d'automatique théorique, Paris 1965)}. Dunod, 1967.

\bibitem{Emerson_Halpern:1985}
E.A. Emerson and J.Y. Halpern.
\newblock Decision procedures and expressiveness in the temporal logic of
  branching time* 1.
\newblock {\em Journal of computer and system sciences}, 30(1):1--24, 1985.

\bibitem{Ghrist:1995}
Robert~W. Ghrist.
\newblock Flows on ${S}^3$ supporting all links as orbits.
\newblock {\em Electron. Res. Announc. Amer. Math. Soc}, 1(2):91--97, 1995.

\bibitem{GH:1996}
Robert~W. Ghrist and Philip~J. Holmes.
\newblock {An ODE whose solutions contain all knots and links.}
\newblock {\em Int. J. Bifurcation Chaos Appl. Sci. Eng.}, 6(5):779--800, 1996.

\bibitem{GHS:1997}
Robert~W. Ghrist, Philip~J. Holmes, and Michael~C. Sullivan.
\newblock {\em {Knots and links in three-dimensional flows.}}
\newblock {Lecture Notes in Mathematics. 1654. Berlin: Springer}, 1997.

\bibitem{Giunti_Mazzola:2010}
M.~Giunti and C.~Mazzola.
\newblock {Dynamical systems on monoids: toward a general theory of
  deterministic systems and motion}, {Q}uinto congresso nazionale di sistemica,
  {A}ssociazione italiana per la ricerca sui sistemi, {F}ermo,
  \begin{scriptsize}
  {https://www.alophis.unica.it/files/Dynamical\%20Systems\%20on\%20Monoids.pdf}
  \end{scriptsize}.
\newblock october 2010.

\bibitem{PIZ:2012}
P.~Iglezias-Zemmour.
\newblock {\em Diffeology}.
\newblock à paraître.
\newblock \\{http://math.huji.ac.il/$\sim$piz/documents/Diffeology.pdf}.

\bibitem{Kanenobu:198504}
Taizo Kanenobu.
\newblock {Satellite links with Brunnian properties.}
\newblock {\em Arch. Math.}, 44(4):369--372, 1985.

\bibitem{Kanenobu:1986}
Taizo Kanenobu.
\newblock {Hyperbolic links with Brunnian properties.}
\newblock {\em J. Math. Soc. Japan}, 38:295--308, 1986.

\bibitem{Matheron.Serra:1988}
Georges Matheron and Jean Serra.
\newblock {\em Image analysis and Mathematical morphology}, volume~2.
\newblock Academic Press, London, 1988.

\bibitem{Matheron.Serra:1988k}
Georges Matheron and Jean Serra.
\newblock Strong filters and connectivity.
\newblock In {\em Image analysis and Mathematical morphology}, volume~2, pages
  141--157. Academic Press, London, 1988.

\bibitem{Mazzola:2010}
Claudio Mazzola.
\newblock {\em Temporal becoming and the Algebra of time}.
\newblock PhD thesis, Università degli studi di Cagliari, 2010.
\newblock \\http://veprints.unica.it/609/1/PhD\_ClaudioMazzola.pdf.

\bibitem{Muscat_Buhagiar:2006}
Joseph Muscat and David Buhagiar.
\newblock Connective space.
\newblock {\em Mem. Fac. Sci. Eng. Shimane Univ. (Series B: Mathematical
  Science)}, (39):1--13, 2006.

\bibitem{Rolfsen:1976}
Dale Rolfsen.
\newblock {\em Knots and links}.
\newblock {P}ublish or {P}erish, Inc., Houston, 1976, sec. ed. 1990.

\end{thebibliography}
%
%







\tableofcontents

\end{document}